\documentclass[12pt,reqno]{amsart}
\usepackage{latexsym}
\usepackage{amssymb}
\newtheorem{theorem}{Theorem}
\newtheorem{proposition}{Proposition}
\newtheorem{lemma}{Lemma}
\newtheorem{corollary}{Corollary}

\def\Hom{\operatorname{Hom}}
\newcommand{\sq}{\mathrm{sq}}
\newcommand{\R}{\mathbb{R}}
\newcommand{\NN}{\mathbb{N}}

\newcommand{\C}{\mathbb{C}}
\renewcommand{\labelenumi}{(\theenumi)}

\newcommand{\F}{\mathcal{F}}
\newenvironment{proofof}[1]{\par\noindent{\em {\it Proof of #1.}}}{\hfill $\square$\par\medskip}
\begin{document}
\title[Characters, Fuchsian groups, and random walks]{Character theory of symmetric groups, subgroup growth
of Fuchsian groups, and random walks}
\maketitle
\begin{center}
\textsc{Thomas W.
M\"uller}\hspace{2.5mm}and\hspace{2.5mm}\textsc{Jan-Christoph
Schlage-Puchta}
\end{center}
\renewcommand{\theenumi}{\roman{enumi}}
\renewcommand{\labelenumi}{(\theenumi)}
\section{Introduction}
\label{Sec:Intro}
The purpose of this paper is three-fold. On the one hand -- and that was its original motivation -- we
establish an asymptotic estimate for the subgroup growth of Fuchsian groups, that is, groups $\Gamma$ of the form
\begin{multline}
\label{Eq:GammaIntro}
\Gamma = \Big\langle x_1, \ldots, x_r, y_1, \ldots, y_s, u_1, v_1, \ldots, u_t, v_t\,\big|\\
x_1^{a_1}=\cdots=x_r^{a_r} = x_1\cdots x_r y_1^{e_1}\cdots y_s^{e_s} [u_1, v_1]\cdots
[u_t, v_t] = 1\Big\rangle
\end{multline}
with integers $r,s,t\geq0$ and $e_1,\ldots, e_s\geq2$, and $a_1,\ldots,a_r\in\NN\cup\{\infty\}$.\\[3mm]
{\bf Theorem A.} {\em Let $\Gamma$ be a Fuchsian group such that
\begin{equation}
\label{Eq:ThmACond}
\alpha(\Gamma):= \sum_i\big(1-\frac{1}{a_i}\big)\hspace{.6mm}+\hspace{.6mm}\sum_j
\frac{2}{e_j}\hspace{.6mm}+\hspace{.6mm}2(t-1) > 0,
\end{equation}
and let
\[
\mu(\Gamma) = \sum_i\big(1-\frac{1}{a_i}\big)\hspace{.6mm}+\hspace{.6mm}s\hspace{.6mm}+\hspace{.6mm}2(t-1)
\]
be the hyperbolic measure of $\Gamma$. Then the number $s_n(\Gamma)$ of index $n$ subgroups in
$\Gamma$ satisfies an asymptotic expansion
\begin{equation}
\label{Eq:AsympIntro}
s_n(\Gamma) \approx \delta L_\Gamma (n!)^{\mu(\Gamma)}
\Phi_\Gamma(n)\bigg\{1\hspace{.6mm}+\hspace{.6mm}\sum_{\nu=1}^\infty
a_\nu(\Gamma)\hspace{.3mm}n^{-\nu/m_\Gamma}\bigg\},\quad(n\rightarrow\infty).
\end{equation}
Here,
\begin{eqnarray*}
\delta  &=& \begin{cases} 2, &\forall i: a_i \text{\em\ finite and odd}, \forall j: e_j\text{\em\ even}\\[1mm]
1, & \text{\em otherwise,}
\end{cases}\\
L_\Gamma & = & (2\pi)^{-1/2-\sum_i(1-1/a_i)} (a_1\cdots a_r)^{-1/2} \exp
\left(-\sum_{\underset{2|a_i}{i}}\frac{1}{2a_i}\right),\\[1mm]
\Phi_\Gamma(n) & = & n^{3/2-\sum_i(1-1/a_i)} \exp\left(\sum_{i=1}^r\sum_{\underset{t<a_i}{t|a_i}}
\frac{n^{t/a_i}}{t}\right),\\[1mm]
m_\Gamma &=& [a_1,a_2,\ldots,a_r],
\end{eqnarray*}
and the $a_\nu(\Gamma)$ are explicitly computable constants depending only on $\Gamma$.}\\[3mm]
On the other hand, the proof of Theorem~A requires the representation theoretic approach initiated in \cite{MuPuSurf},
and large parts of the present paper are concerned with various statistical aspects of symmetric groups, and rather subtle
estimates for values and multiplicities of their characters, which are also of independent interest.
In particular, we show the following.\\[3mm]
{\bf Theorem B.}{\em\ Let $\varepsilon>0$ and an integer $q$ be given, $n$ sufficiently large, and let $\chi$ be an
irreducible character of $S_n$.\vspace{-3mm}
\begin{enumerate}
\item We have $|\chi(\mathbf{c})|\leq \big(\chi(1)\big)^{1-\delta}$ with
\[
\delta = \left(\big(1-1/(\log n)\big)^{-1}\frac{12\log n}{\log (n/f)} + 18\right)^{-1},
\]
where $\mathbf{c}$ is any conjugacy class of $S_n$ with $f$ fixed points.\\[-3mm]
\item We have
\[
\sum_{\pi^q=1}|\chi(\pi)| \leq \big(\chi(1)\big)^{\frac{1}{q}+\varepsilon}\sum_{\pi^q=1} 1.
\]
\item Let $m_\chi^{(q)}$ be the multiplicity of $\chi$ in the $q$-th root number
function of $S_n$. Then
\[
m_\chi^{(q)}\leq \big(\chi(1)\big)^{1-2/q+\varepsilon}.
\]
\end{enumerate}
}
All these bounds are essentially best possible. In characteristic zero, estimates
for character values and multiplicities are among the least understood topics
in the representation theory of symmetric groups. In recent times,
additional interest in this circle of problems was sparked by the theory
of random walks on finite groups.
In this context, Theorem B enables us to prove the following.\\[3mm]
{\bf Theorem C.}{\em\ Let $\mathbf{c}$ be a conjugacy class in $S_n$. Denote
by $t_c(\mathbf{c})$ the least even integer such that $t_c(\mathbf{c})$
elements chosen at random from $\mathbf{c}$ have, with probability $\geq 1-\frac{1}{n},$
no common fixed point, and let $t_s(\mathbf{c})$ be the mixing time for
the random walk generated by $\mathbf{c}$. Then, for $n\geq 4000,$ we have}
\[
t_c(\mathbf{c})\leq t_s(\mathbf{c})\leq 10 t_c(\mathbf{c}).
\]
Theorem C establishes in full generality a conjecture of Roichman;
cf. \cite[Conj. 6.6]{Roich}. For special choices of $\mathbf{c}$,
Roichman's conjecture had already been known to hold: Diaconis and Shahshahani \cite{DiaSha}
established it for transpositions, Roichman \cite{Roich} generalized their
result to conjugacy classes having at least $c n$ fixed points, and
Fomin and Lulov \cite{FoLu} established a character bound
implying Theorem C for conjugacy classes having only cycles of the same length.\\[3mm]
In the remainder of this introduction, we describe in more detail the organization of
our paper, and the contents of individual sections. After a short introduction to random
walks on finite groups generated by conjugacy classes and their connection with character
theory, Section~\ref{Sec:Random} establishes Theorem~C and the first part of Theorem~B.
This includes the proof of a variety of preliminary character estimates for arbitrary conjugacy classes,
which are used throughout the paper. The main tools here are the hook formula and the Murnaghan-Nakayama
 rule. Section~\ref{Sec:MainLemFuchs} gives the proof of Theorem~B (ii). In preparation for this
 argument, we derive a number of results concerning the statistics of symmetric groups, mostly dealing with
 the distribution of cycles in various subsets of $S_n$; cf. Subsection~\ref{Subsec:Statistics}. Again, this
 group of results is also used in other sections. The theory developed up to this point would already
 enable us to estimate the subgroup growth of Fuchsian groups with $s=0$; that is, none of the generators
 $y_1,\ldots,y_s$ in (\ref{Eq:GammaIntro}) are present. However, in order to deal with Fuchsian groups in
 full generality, we also need some insight into the growth behaviour of multiplicities of root number functions for
 symmetric groups, measured against the degree of the corresponding irreducible character;
 in particular, we have to establish Theorem~B (iii). Section~\ref{Sec:Multiplicities} is
 devoted to the proof of this and related results.\\[3mm]
 Proof and discussion of Theorem~A are the principal themes of Section~\ref{Sec:Growth}. Following the argument
 establishing Theorem~A, we demonstrate that condition (\ref{Eq:ThmACond}) is in fact necessary. More specifically,
 we show the following. \\[3mm]
 {\bf Theorem~D.} {\it Let $\Gamma$ be as in {\em (\ref{Eq:GammaIntro})} with $r=t=0$ and $\alpha(\Gamma)<0$ (that is,
 $\Gamma$ is the one-relator group associated with the defining relation $y_1^{e_1}y_2^{e_2}\cdots y_s^{e_s}=1$).
 Then, as $n$ tends to infinity, we have}
\[
s_n(\Gamma) \sim K (n!)^{\mu(\Gamma)-\alpha(\Gamma)/2}
\exp\Bigg(\sum_{j=1}^s
\underset{\nu_j<e_j}{\sum_{v_j\mid e_j}} \frac{n^{\nu_j/e_j}}{\nu_j}\hspace{.8mm}+\hspace{.8mm}
\frac{\alpha(\Gamma)-2\mu(\Gamma)+2}{4} \log n\Bigg).
\]
This is the contents of Theorem~\ref{Thm:OneRelAsymp}, where also the constant $K$ is given explicitly. According to
Theorem~D, the subgroup growth of these one-relator groups is faster than might be expected in view of Theorem~A.
 In Subsection~\ref{Subsec:CoeffComp}, we discuss the explicit computation of the coefficients $a_\nu(\Gamma)$ in general,
and, as an example, compute the first $22$ of these coefficients for the triangle group $\Gamma(2,3,7)$, only $10$ of
which turn out to be non-vanishing. As a further application of Theorem~B (iii), we also determine the subgroup growth of
Demu\v skin groups (Poincar\'e groups of dimension $2$). These are known to be one-relator groups essentially
associated with a defining relation of the form
\[
R=x_1^{p^h}[x_1, x_2][x_3, x_4]\cdots [x_{m-1}, x_m],\quad h\in\NN\cup\{\infty\},\; p\geq3\mbox{ prime}.
\]
In Subsection~\ref{subsec:Demu}, we prove the following slightly more general result.\\[3mm]
{\bf Theorem~E.} {\it For integers $q\geq 1$ and $d\geq 2,$ let
\[
\Gamma_{q, d} = \Big\langle x_1, y_1, \ldots, x_d, y_d\,\big|\,x_1^{q-1}[x_1, y_1]\cdots[x_d, y_d]=1\Big\rangle.
\]
Then there exist explicitly computable constants $\gamma_\nu(\Gamma_{q, d}),$ such that
\[
s_n(\Gamma_{q, d}) \approx \delta n (n!)^{2d-2}\left\{1+\sum_{\nu=1}^\infty
\gamma_{\nu}(\Gamma_{q, d})\hspace{.3mm}n^{-\nu}\right\},\quad n\rightarrow\infty,
\]
where}
\[
\delta=\begin{cases}1, &q\mbox{ even}\\2, &q\mbox{ odd.}\end{cases}
\]

Introduce an equivalence relation $\sim$ on the class of finitely generated groups via
\[
\Gamma \sim \Delta:\Leftrightarrow s_n(\Gamma)=(1+o(1))s_n(\Delta),\quad(n\rightarrow\infty).
\]
In \cite[Theorem 3]{MInvent} a characterization in terms of
structural invariants is given for the equivalence relation $\sim$
on the class of groups $\Gamma$ of the form
\[
\Gamma=G_1*G_2*\dots*G_s*F_r
\]
with $r, s\geq 0$ and finite groups $G_\sigma$, and it is shown
that  each $\sim$-class of free products decomposes into finitely many isomorphism
classes. Our final section is concerned with the analogous problems for
Fuchsian groups.\\[3mm]
{\bf Theorem~F.} {\it The multiset $\{a_1, a_2\ldots,a_r\}$ together with the numbers $\mu(\Gamma)$ and $\delta$
form a complete system  of invariants for the equivalence relation $\sim$ on the class $\mathcal{F}$ of all Fuchsian groups
$\Gamma$ satisfying $\alpha(\Gamma)>0$.}\\[3mm]
Theorem~F allows us to construct an infinite sequence of pairwise non-isomorphic Fuchsian groups, all of
which are $\sim$-equivalent to the same Fuchsian group $\Gamma$; in particular there cannot be a
finiteness result for the relation $\sim$ on $\mathcal{F}$. The situation changes, if we take into
account the full precision of (\ref{Eq:AsympIntro}) in Theorem~A. More specifically, consider three refinements
of the equivalence relation $\sim$ on $\F$:
 (i) the relation $\approx$ of strong equivalence defined via
\[
\Gamma\approx\Delta :\Leftrightarrow
s_n(\Gamma)=s_n(\Delta)(1+\mathcal{O}(n^{-A}))\mbox{ for every }A>0,
\]
(ii) isomorphy, and (iii) equality of the system of parameters
\[
(r,s,t;a_1,a_2\ldots,a_r,e_1,e_2,\ldots,e_s)
\]
in the Fuchsian presentation (\ref{Eq:GammaIntro}), denoted $\Gamma=\Delta$. Clearly,
\[
\Gamma = \Delta \Rightarrow \Gamma\cong\Delta \Rightarrow \Gamma\approx\Delta \Rightarrow \Gamma\sim\Delta.
\]
It can be shown that all these implications are in fact strict. For these relations, we have the
following surprising result.\\[3mm]
{\bf Theorem~G.} {\it Each $\approx$-equivalence class of $\F$ decomposes into finitely many classes
with respect to $=;$ that is, each group $\Gamma\in\F$ has only finitely many presentations of the
form {\em (\ref{Eq:GammaIntro}),} and is $\approx$-equivalent
to at most finitely many non-isomorphic $\F$-groups.}\\[3mm]
{\bf Some notation.} Permutations are denoted by $\pi$, $\sigma$, or $\tau$. For $\pi\in S_n$ and $1\leq i\leq n$,
let $s_i(\pi)$ be the number of $i$-cycles of $\pi$. The {\it support} $\mathrm{supp}(\pi)$ of $\pi$ is the set
of points moved by $\pi$. For number partitions we mostly follow the conventions of \cite[Chap.~I, Sec.~1]{Macdonald}.
Specifically, a partition $\lambda=(\lambda_1,\lambda_2,\ldots)$ is a weakly decreasing sequence of
non-negative integers $\lambda_j$, such that $\lambda_j=0$ for $j$ sufficiently large. The {\it weight} $|\lambda|$
of $\lambda$ is $|\lambda|=\sum_j \lambda_j$, and the {\it norm} $\|\lambda\|$ of $\lambda$ is  the largest $j$ such
that $\lambda_j\neq0$. As usual, we write $\lambda\vdash n$ for $|\lambda|=n$, and say that $\lambda$ is a partition
of $n$. For partitions $\lambda,\mu$ we write $\mu\subseteq\lambda$, if $\mu_j\leq \lambda_j$ for all $j$ (that is,
the Ferrers diagram of $\mu$ is contained in the Ferrers diagram of $\lambda$). For a partition $\lambda$, we denote
by $\lambda'$ the conjugate partition: $\lambda_i' = \max\{j: \lambda_j\geq i\}$ (that is, the Ferrers diagram
of $\lambda'$ is obtained from that of $\lambda$ by reflection through the main diagonal). By $\lambda\setminus\lambda_1$
we mean the partition $\lambda\setminus\lambda_1 = (\lambda_2,\lambda_3,\ldots)$ (that is, the Ferrers
diagram of $\lambda\setminus\lambda_1$ is obtained from that of $\lambda$ by deleting the first row). Whenever
convenient, we shall denote the integer $|\lambda|-\lambda_1$ by $\Delta$. For a partition
$\lambda\vdash n$, we denote by $\chi_\lambda$ the irreducible character of $S_n$ associated with $\lambda$.
For a finite group $G$, let $\mathrm{Irr}(G)$ be the set of irreducible characters of $G$. The usual scalar product
on the space $\C^G$ is denoted by $\langle\cdot,\cdot\rangle_G$, or simply by $\langle\cdot,\cdot\rangle$ if $G$ is
a symmetric group.\\[3mm]
We use what we believe to be standard number theoretic notation. Specifically, the partition function is
denoted $p(n)$, $\tau(n)$ and $\sigma(n)$ are the number of divisors and the sum of divisors of $n$, respectively,
$S(n,m)$ is the number of (set theoretic) partitions of an $n$-set into $m$ non-empty parts (a Stirling number
of the second kind). For integers $m$ and $n$, we denote their greatest common divisor and least common multiple
by $(m,n)$ respectively $[m,n]$. For arithmetic functions $f,g:\NN\rightarrow\R$ we use the
Vinogradov symbol $f(n)\ll g(n)$ to mean $f(n)=\mathcal{O}(g(n))$. If $f(n)\ll g(n)\ll f(n)$, we write $f(n)\asymp g(n)$.
Asymptotic equivalence is denoted by $\sim$: we write $f(n)\sim g(n)$ to mean $f(n)=g(n)(1+o(1))$. We use $\approx$
to denote asymptotic expansions in the sense of Poincar\'e; for instance we write
\[
f(n) \approx g(n)\bigg\{1\hspace{.6mm}+\hspace{.6mm}\sum_{\nu=1}^\infty a_\nu\hspace{.3mm}n^{-\nu/q}\bigg\},
\quad(n\rightarrow\infty).
\]
Finally, we use some notation from probability theory: ${\bf 1}_X$ denotes the characteristic function for a
subset $X\subseteq\Omega$
 of the sample space $\Omega$, and ${\bf E}\xi$ is the expected value of the random variable $\xi$.
\section{Character Estimates and Random Walks on symmetric groups}
\label{Sec:Random}
\subsection{Roichman's Conjectures} Let $\mathbf{E}=E_1, E_2, E_3, \ldots$ be a Markov
chain on a metric space $X$. The random walk $(x_k)_{k\geq 0}$ on
$X$ determined by $\mathbf{E}$ is by definition the collection of
all infinite pathes on $X$ with probability distribution induced
by $\mathbf{E}$. If $X$ is finite, then $\mathbf{E}$ can be
determined by its transition matrix $P$. In what follows, we will
be interested in the case when $X$ is a finite symmetric group
given with the discrete metric.\footnote{Cf. \cite {Diac} for more
details on random walks on finite groups and their applications.}
Let $X=S_n$, and let $1\neq\mathbf{c}\subseteq S_n$ be a
non-trivial conjugacy class. The random walk $w_{\mathbf{c}}$ {\em
generated by }$\mathbf{c}$ has initial state $x_0=1$ and the
transition matrix
$P_{\mathbf{c}}=(p^{\mathbf{c}}_{\pi\sigma})_{\pi, \sigma\in S_n}$
where
\[
p_{\pi\sigma}^{\mathbf{c}} := \begin{cases} \frac{1}{|{\bf c}|}, &
\pi\sigma^{-1}\in\mathbf{c}\\ 0, & \mbox{otherwise}.
\end{cases}
\]
The distribution in the $k$-th step of $w_{\mathbf{c}}$ is given
by the convolution formula
\[
P(x_k=\pi)=\frac{1}{|\mathbf{c}|}\sum_{\sigma\in S_n}
P(x_{k-1}=\pi\sigma^{-1}).
\]
More generally, for two functions $f, g: S_n\rightarrow\C$, the
convolution $f\ast g:S_n\rightarrow\C$ is given by
\[
(f\ast g)(\pi)=\sum_{\sigma\in S_n}f(\sigma)g(\pi\sigma^{-1});
\]
in particular, $P(x_k=\pi)$ is the $k$-fold convolution of the
function $\frac{1}{|\mathbf{c}|}{\bf 1}_{\mathbf{c}}$. In the
sequel we shall always take $k$ even, to avoid parity problems.
Given a norm $\|\cdot\|$ on the complex algebra $\C^{S_n}$ and
$\varepsilon>0$, we say that the random walk $w_{\mathbf{c}}$ has
reached $\varepsilon$-equidistribution with respect to $\|\cdot\|$ in
step $k$, if
\[
\left\|P(x_k=\pi)-\frac{2}{n!}{\bf 1}_{A_n}\right\|<\varepsilon\cdot
\left\|\frac{2}{n!}\mathbf{1}_{A_n}\right\|.
\]
We define the {\em statistical mixing time} $t_s(\mathbf{c})$ of
$w_{\mathbf{c}}$ as the least even integer $k$ for which
$w_{\mathbf{c}}$ has reached $\frac{1}{n}$-equidistribution with
respect to the $\ell^2$-norm. A first lower bound for
$t_s(\mathbf{c})$ is given by the combinatorial mixing time
$t_c(\mathbf{c})$ of $\mathbf{c}$, that is, the least even integer
$k$, such that any $k$ elements of $\mathbf{c}$ have no common
fixed point with probability at least $1-\frac{1}{n}$. In
\cite{Roich}, Roichman conjectured that for every non-trivial
conjugacy class $\mathbf{c}\subseteq A_n$,
\begin{equation}
\label{Eq:RoichConj1}
t_s(\mathbf{c})\ll t_c(\mathbf{c}).
\end{equation}
His main result \cite[Theorem 6.1]{Roich} establishes this
conjecture for classes $\mathbf{c}$ with $c n$ fixed points.
Roichman suggests an approach to the general conjecture
(\ref{Eq:RoichConj1}) via a certain estimate for character values in
symmetric groups. More precisely, he conjectures that, for every
$\varepsilon>0$, $n$ sufficiently large, each conjugacy class
$\mathbf{c}\subseteq S_n$, and every partition $\lambda\vdash n$
\begin{equation}
\label{Eq:RoichConj2}
|\chi_\lambda(\mathbf{c})|\leq
\chi_\lambda(1)\left(\max\left( \frac{\lambda_1}{n},
\frac{\|\lambda\|}{n},
\frac{1}{e}\right)\right)^{(1-\varepsilon)n\log\frac{n}{n-\mathrm{supp}(\mathbf{c})+1}},
\end{equation}
which would imply (\ref{Eq:RoichConj1}). Unfortunately, as it
stands, estimate (\ref{Eq:RoichConj2}) is false. This can be seen,
for instance, as follows. For $\mathbf{c}$ fixed-point free,
$\lambda$ such that $\lambda_1, \|\lambda\| \leq\frac{n}{e}$, and
$\varepsilon=\frac{1}{2}$, (\ref{Eq:RoichConj2}) simplifies to
\begin{equation}
\label{Eq:RoichSimpl} |\chi_\lambda(\mathbf{c})| \leq
\chi_\lambda(1) e^{-(n\log n)/2}.
\end{equation}
The right-hand side of (\ref{Eq:RoichSimpl}) is bounded above by
$\sqrt{n!}\,n^{-n/2}<1$; that is, for $\mathbf{c}$ and $\lambda$
as above, and $n$ sufficiently large, it would follow that
$\chi_\lambda(\mathbf{c})=0$. Since the irreducible characters
$\{\chi_\lambda\}_{\lambda\vdash n}$ form a basis for the space of
class functions on $S_n$, this would imply that, for $n$
sufficiently large, the characters
\[
\big\{\chi_\lambda:\lambda\vdash n, \max(\lambda_1,
\|\lambda\|)>n/e\big\}
\]
would generate the space of class functions on the set of fixed-point free conjugacy classes of $S_n$. Comparing the size of the
former set with the dimension of the latter space, we find that,
for large $n$,
\begin{equation}
\label{Eq:Roichpartition}
2\sum_{0\leq\nu\leq n-n/e} p(\nu) \geq
p(n)-p(n-1),
\end{equation}
where $p(n)$ is the number of partitions of $n$. The right-hand
side can be estimated via the first term of Rademacher's series
expansion for $p(n)$ (see for example \cite[Theorem 5.1]{Andrews})
to give
\[
p(n)-p(n-1) \sim \frac{\pi e^{\pi\sqrt{\frac{2n}{3}}}}
{12\sqrt{2}n^{3/2}},\quad n\rightarrow\infty.
\]
Bounding the left-hand sum in (\ref{Eq:Roichpartition}) by means
of the estimate\footnote{Cf., for instance, \cite[Satz 7.6]
{Kraetzel}.} $p(n)<\frac{\pi}{\sqrt{6n}}e^{\pi\sqrt{2n/3}}$ we
obtain
\[
2\sum_{0\leq\nu\leq n-n/e} p(\nu) \leq 2np(n-\lfloor n/e\rfloor)
\leq 2\pi\sqrt{n/6} e^{\pi\sqrt{2(1-1/e)n/3}}.
\]
From these two estimates it is clear that inequality
(\ref{Eq:Roichpartition}) is violated for large $n$.\\[3mm]
However, the basic idea behind Roichman's approach turns out to be
correct. As a substitute for (\ref{Eq:RoichConj2}), we prove the
following.
\begin{theorem}
\label{thm:Roichcorrect} For sufficiently large $n,$ a non-trivial
conjugacy class $\mathbf{c}\subseteq S_n,$ and a partition
$\lambda\vdash n,$ we have
\begin{equation}
\label{Eq:thm1:charest} |\chi_\lambda(\mathbf{c})| \leq
\big(\chi_\lambda(1)\big)^{1-\frac{1-1/(\log n)}{6
t_c(\mathbf{c})}}
\end{equation}
and, for $1\leq s_1(\mathbf{c})\leq n-2,$
\begin{equation}
\label{Eq:thm1:tcest} \left|t_c(\mathbf{c}) - \frac{2\log
n}{\log(n/s_1(\mathbf{c}))}\right|\leq 3,
\end{equation}
whereas $t_c(\mathbf{c})=2$ for $s_1(\mathbf{c})=0$.
\end{theorem}
This result in turn allows us to establish Roichman's conjecture
(\ref{Eq:RoichConj1}) for the mixing time of random walks on
symmetric groups.
\begin{theorem}
\label{thm:walks} For $n\geq 4000$ and each non-trivial
conjugacy class $\mathbf{c}\subseteq S_n,$ we have
\begin{equation}
\label{Eq:thm2} t_c(\mathbf{c})\leq t_s(\mathbf{c}) \leq 10
t_c(\mathbf{c}).
\end{equation}
\end{theorem}
The constants in Theorems \ref{thm:Roichcorrect} and
\ref{thm:walks} are most certainly not optimal, but we have not
attempted to tighten our numerical estimates.\\[3mm]
This section is organized as follows. In Subsection
\ref{Subsec:CharAndWalk} we describe the general connection
between random walks on finite groups and character estimates, and
we explain, how Theorem \ref{thm:walks} can be deduced from
Theorem \ref{thm:Roichcorrect}. The next subsection establishes
certain elementary estimates for values and degrees of irreducible
characters of symmetric groups, which will be used throughout the
paper. Finally, Subsections \ref{Subsec:ProofThmRoich} and
\ref{Subsec:ProofLemRoich} contain the proof of Theorem 1.
\subsection{Character Theory and Random walks}
\label{Subsec:CharAndWalk} Here, we describe the connection between
character theory and probability measures on finite groups. For a
more detailed presentation, see \cite{Diac} and \cite{Terras}. Let
$G$ be a finite group. For a class function $f:G\rightarrow\C$ and
an irreducible character $\chi$ of $G$, define the Fourier
coefficient $\alpha_\chi(f)$ by means of the equation
\[
f(g)=\sum_{\chi}\alpha_\chi(f)\chi(g), \quad g\in G.
\]
Since irreducible characters form a basis of the space of class
functions, the Fourier coefficients exist and are uniquely defined
by this equation. We first state without proof some facts about
Fourier coefficients.
\begin{lemma}
\label{Lem:FourierGen}
\begin{enumerate}
\item $\alpha_\chi(f)=\langle f, \chi\rangle,$\\[-3mm]
\item $\alpha_\chi(f\ast g) = \frac{\alpha_\chi(f)\alpha_\chi(g)}
{\chi(1)},$\\[-3mm]
\item $\sum_{\chi}|\alpha_\chi(f)|^2=\frac{1}{|G|}\sum_{g\in G}|f(g)|^2.$
\end{enumerate}
\end{lemma}
Using Lemma~\ref{Lem:FourierGen}, we obtain a bound for the
statistical mixing time of $w_{\mathbf{c}}$.
\begin{lemma}
\label{Lem:DistributionBound} Let $\mathbf{c}\subseteq G$ be a
conjugacy class, and let $k$ be an integer. Then we have
\[
|G|\left\|\left(\frac{1}{|\mathbf{c}|}
\mathbf{1}_{\mathbf{c}}\right)^{\ast k} -
\frac{1}{|G|}\mathbf{1}\right\|_2^2 = \sum_{\chi\neq\chi_0}
\frac{|\chi(\mathbf{c})|^{2k}}{(\chi(1))^{2k-2}},
\]
where $\chi_0$ is the trivial character.
\end{lemma}
\begin{proof}
Noting that $\alpha_\chi(\mathbf{1})=\delta_{\chi,\chi_0}$, this follows from Lemma~\ref{Lem:FourierGen}.
\end{proof}
In contrast to \cite{Diac}, we will only consider the
$\ell^2$-norm. The passage from the $\ell^2$-norm to
$\ell^p$-norms with $p\leq 2$ is immediate by H\"older's
inequality; however, in the context of random walks, the passage
from $\ell^2$ to $\ell^\infty$ is possible as well.
\begin{lemma}
\label{Lem:l2tolinfty} Let $G$ be a finite group, and let
$f:G\rightarrow [0, \infty)$ be a probability measure. Then we
have
\[
\left\|f\ast f - \frac{1}{|G|}\mathbf{1}\right\|_\infty \leq
\left\|f-\frac{1}{|G|}\mathbf{1}\right\|_2^2.
\]
In particular, if $w_\mathbf{c}$ reaches
$\varepsilon$-equidistribution with respect to the $\ell^2$-norm
after $k$ steps, it reaches $\varepsilon^2$-equidistribution with
respect to the $\ell^\infty$-norm after $2k$ steps.
\end{lemma}
\begin{proof}
We have
\begin{eqnarray*}
(f\ast f)(g) & = & \sum_{h\in G} f(h)f(gh^{-1})\\
 & = & \sum_{h\in G}\left(f(h)-\frac{1}{|G|}\right)
\left(f(gh^{-1})-\frac{1}{|G|}\right) + \frac{2}{|G|}\sum_{h\in G} f(h)  - \frac{1}{|G|}\\
 & = & \frac{1}{|G|} + \sum_{h\in G}\left(f(h)-\frac{1}{|G|}\right)
\left(f(gh^{-1})-\frac{1}{|G|}\right),
\end{eqnarray*}
since $\sum f(h)=1$. Applying the Cauchy-Schwarz-inequality to the
last sum, our claim follows.
\end{proof}
We can now explain how to deduce Theorem~\ref{thm:walks} from
Theorem~\ref{thm:Roichcorrect}. Let $\mathbf{c}\subseteq S_n$ be a
non-trivial conjugacy class. Arguing as in the proof of
Lemma~\ref{Lem:DistributionBound} we know that $w_\mathbf{c}$
reaches $\frac{1}{n}$-equidistribution with respect to the
$\ell^2$-norm after $k$ steps, if and only if
\[
\underset{\chi(1)\neq 1}{\sum_{\chi}}
\frac{|\chi(\mathbf{c})|^{2k}}{(\chi(1))^{2k-2}} \leq \frac{1}{n}.
\]
By Theorem \ref{thm:Roichcorrect}, the left-hand side can be
bounded above by
\[
\underset{\chi(1)\neq 1}{\sum_{\chi}}
\big(\chi(1)\big)^{2-\frac{2k(1-1/(\log n))}{6t_c(\mathbf{c})}}.
\]
For $k\geq 10t_c(\mathbf{c})$, this in turn is less than
\[
\underset{\chi(1)\neq 1}{\sum_{\chi}} \big(\chi(1)\big)^{-5/4},
\]
and from \cite[Theorem 1]{MuPuSurf} we deduce that the latter quantity
is $\mathcal{O}(n^{-5/4})$. Hence, for $n$ sufficiently large, we
obtain the bound $t_s(\mathbf{c})\leq 10 t_c(\mathbf{c})$.
We postpone the argument leading to the lower bound $n\geq 4000$
to the end of the proof of Theorem 1 in Subsection~\ref{Subsec:ProofThmRoich}.
\subsection{Estimates for Character Values}
\label{Subsec:CharElemEst} Our main tools in this subsection are
the hook formula for the dimension of $\chi_\lambda$ and the
Murnaghan-Nakayama rule.\\[2mm]
{\bf The hook formula.} {\it We have
\begin{equation}
\chi_\lambda(1) = \frac{n!}{\prod_{(i,j)\in\lambda} h_{i,j}},
\end{equation}
where $h_{i,j}$ is the hook length of the box $(i,j)$.}\\[2mm]
The Murnaghan-Nakayama rule describes a procedure to recursively
compute arbitrary character values.\\[2mm]
{\bf Murnaghan-Nakayama rule.} {\it Let $\pi = \sigma \gamma$ be
the disjoint product of $\sigma\in S_{n-k}$ and a $k$-cycle
$\gamma$. Then we have
\begin{equation}
\label{Eq:Nakayama} \chi_\lambda(\pi) = \sum_\tau (-1)^{l(\tau)}
\chi_{\lambda\setminus\tau}(\sigma).
\end{equation}
 Here, the sum extends over all rim hooks $\tau$ of size $k$ in
 $\lambda,$ and $l(\tau)$ is the leg length of $\tau$.}\\[2mm]
Let $\lambda$ be a partition of $n$. By $\sq(\lambda)$ we mean the
side length of the largest square contained in the Ferrers
diagram of $\lambda$; that is, the largest $j$ such that
$\lambda_j\geq j$. Note that for $\lambda\vdash n$ we have
\begin{equation}\label{Eq:SqBound}
(\sq(\lambda)-1)\sq(\lambda)\leq n-\lambda_1,
\end{equation}
which we will apply mostly in the simpler version
$\sq(\lambda)\leq\sqrt{n-\lambda_1}+1$. The quantity $\sq(\lambda)$
leads to a useful inequality for $\chi_\lambda(1)$.
\begin{lemma}
\label{Lem:squaredegree} Let $\lambda$ be a partition of $n$, and
let $s=\mathrm{sq}(\lambda)$. Then
\[
\chi_\lambda(1)\geq{n\choose s^2} \left(\frac{s}{n}\right)^{s^2}
(s^2)!\hspace{0.2mm}.
\]
\end{lemma}
\begin{proof}
Each of the $n-s^2$ points of $\lambda$ outside the maximal square
lies in precisely $s$ hooks $h_{ij}$ with $i, j\leq s$, while the
point $(i, j)$ with $i, j\leq s$ lies in exactly $i+j-1$ such
hooks. Hence,
\[
\sum_{i, j\leq s}h_{ij} = s(n-s^2) + \sum_{i, j\leq s}(i+j-1) =
sn.
\]
By the arithmetic-geometric inequality, this gives
\[
\prod_{i, j\leq s} h_{ij} \leq \left(\frac{1}{s^2}\sum_{i, j\leq
s} h_{ij}\right)^{s^2} = \left(\frac{n}{s}\right)^{s^2}.
\]
Bounding the product of the hook lengths corresponding to points
outside the maximal square by $(n-s^2)!$, our claim follows from
the hook formula.
\end{proof}
Our next result gives an upper bound for the modulus of character
values $\chi_\lambda(\mathbf{c})$ in terms of
$\mathrm{sq}(\lambda)$ and the number of cycles of $\mathbf{c}$.
\begin{lemma}
\label{Lem:squarevalue} Let $\mathbf{c}\subseteq S_n$ be a
conjugacy class with $c$ cycles, and let $\lambda\vdash n$ be a
partition. Then we have
\[
|\chi_\lambda(\mathbf{c})| \leq
(2\hspace{.3mm}\mathrm{sq}(\lambda))^c.
\]
\end{lemma}
\begin{proof}
If $\mu\subseteq\lambda$ is a partition, then certainly
$\mathrm{sq}(\mu)\leq\mathrm{sq}(\lambda)$, hence, arguing by
induction on $c$ and applying the Murnaghan-Nakayama rule, it
suffices to show that for any given $k$, a partition $\lambda$ has
at most $2\hspace{.3mm}\mathrm{sq}(\lambda)$ removable rim hooks
of length $k$. Let $r$ be such a rim hook. The right-uppermost box of
$r$ is either to the right of the maximal square contained in the
Ferrers diagram of $\lambda$, or the left-lowest box of $r$ is below
the maximal square of
$\lambda$, or both. Since the right-uppermost box of a rim hook is
always at the end of a row, while the left-lowest box is at the end of a
column, our claim follows.
\end{proof}
If $\lambda_1$ is of similar size as $n$, Lemma
\ref{Lem:squarevalue} is of little use. In this case we will apply
the following.
\begin{lemma}
\label{Lem:LargeL1CharEst} Let $\lambda\vdash n$ be a partition,
and let $\mathbf{c}\subseteq S_n$ be a conjugacy class with $c$
cycles and $f$ fixed points. Then we have the bounds
\begin{equation}\label{Eq:L1CharEst1}
|\chi_\lambda(\mathbf{c}))| \leq
\chi_{\lambda\setminus\lambda_1}(1)\underset{a+2b\leq
n-\lambda_1}{\sum_{a,b\geq0}}{f\choose a}{{c-f}\choose b}.
\end{equation}
and
\begin{equation}\label{Eq:L1CharEst2}
|\chi_\lambda(\mathbf{c})| \leq n\max_{\nu\leq n-\lambda_1}
(2\sqrt{n-\lambda_1})^\nu\binom{c}{\nu},
\end{equation}
which improves on {\em (\ref{Eq:L1CharEst1})} if $c$ is
considerably smaller then $n-\lambda_1$.
\end{lemma}
\begin{proof}
Neglecting the sign in the Murnaghan-Nakayama rule, we see that
the modulus of $\chi_\lambda(\mathbf{c})$ is bounded above by the
number of possible ways to completely deconstruct $\lambda$ by
removing rim hooks of sizes given by the cycle structure of
$\mathbf{c}$. To prove the first estimate, we classify these
deconstructions by means of the number $a$ of fixed points of
$\mathbf{c}$ contained within $\lambda\setminus\lambda_1$, and the
corresponding number $b$ of cycles of lengths $\geq 2$. Given $a$
and $b$, there are ${f\choose a}$ ways to choose the fixed points
of $\mathbf{c}$ to be removed from $\lambda\setminus\lambda_1$,
and ${{c-f}\choose b}$ ways to choose the corresponding set of
cycles. Once these sets are given, their are at most
$\chi_{\lambda\setminus\lambda_1}(1)$ ways to remove these fixed
points and cycles.\\[3mm]
For the second bound we argue in a similar manner, this time
bounding the number of deconstructions of
$\lambda\setminus\lambda_1$ as in the proof of
Lemma~\ref{Lem:squarevalue}.
\end{proof}
The next lemma will be useful in computing values of characters of
small degrees.
\begin{lemma}\label{Lem:CharPolynom}
Let $\lambda\vdash n$ be a partition,
$\mu=\lambda\setminus\lambda_1,$ and let $\pi\in S_n$ be a
permutation. Then
\[
\chi_\lambda(\pi) =
\sum_{\underset{\tilde{\mu}_1=1}{\tilde{\mu}\subseteq\mu}}(-1)^{|\tilde{\mu}|}
\sum_{\mathbf{c}\subseteq S_{|\mu|-|\tilde{\mu}|}} \chi_{\mu,
\tilde{\mu}}(\mathbf{c})\prod_{i\leq |\mu|}\binom{s_i(\pi)}{c_i},
\]
where $\mathbf{c}$ runs over all conjugacy classes of
$S_{|\mu|-|\tilde{\mu}|},$ $\chi_{\mu, \tilde{\mu}}(\mathbf{c})$
denotes the number of ways to obtain $\tilde{\mu}$ from $\mu$ by
removing rim hooks according to the cycle structure of
$\mathbf{c},$ counted with the sign prescribed by the
Murnaghan-Nakayama rule, and $c_i$ is the number of $i$-cycles of an element of
$\mathbf{c}$.
\end{lemma}
\begin{proof}
We may assume that $\lambda_1>2|\mu|$, and that $\pi$ contains a
cycle of length $>|\mu|$. For, if we replace $\lambda_1$ by
$\lambda_1+2n$, and add a cycle of length $2n$ to $\pi$, both
conditions are satisfied; on the other hand, the only way to
remove a cycle of length $2n$ from the new partition is within the
first row, hence $\chi_\lambda(\pi)$ is not affected by these
changes. Now we use the Murnaghan-Nakayama rule to remove all
cycles of length $\leq |\mu|$. If we are left with a partition not
of the form $\tilde{\lambda}=(\tilde{\lambda}_1, 1, \ldots, 1)$,
the remaining partition cannot be removed by deleting rim hooks of
length $>|\mu|$, so
$\tilde{\mu}=\tilde{\lambda}\setminus\tilde{\lambda}_1$ can be
assumed to be of the form $(1, \ldots, 1)$. Since
$\tilde{\lambda}$ can be removed in precisely one way by deleting
rim hooks of lengths $>|\mu|$, and all but the last one are contained
in the first row, the value of
$\chi_{\tilde{\lambda}}(\tilde{\pi})$ is $(-1)^{|\tilde{\mu}|}$,
where $\tilde{\pi}$ is the element obtained from $\pi$ by removing
all cycles of lengths $\leq |\mu|$. The rim hooks which are
removed and do not contain any box from the first row define a
conjugacy class $\mathbf{c}$ within $S_{|\mu|-|\tilde{\mu}|}$;
counting all possible placements of the cycles of $\mathbf{c}$ in
$\mu\setminus\tilde{\mu}$ with the sign prescribed by the
Murnaghan-Nakayama rule yields $\chi_{\mu,
\tilde{\mu}}(\mathbf{c})$. Finally, for fixed $\mathbf{c}$, the
set of $i$-cycles to be placed in $\mu\setminus\tilde{\mu}$ among
all $i$-cycles of $\pi$ can be chosen in $\binom{s_i(\pi)}{c_i}$
ways, and our claim follows.
\end{proof}
Finally, we shall also need the following lower bounds for
character degrees.
\begin{lemma}
\label{Lem:chi1lower} Let $\lambda\vdash n$ be a partition. Then
\vspace{-3mm}
\begin{enumerate}
\item $\displaystyle\chi_\lambda(1)\geq 2^{n/8},\quad \|\lambda\|\leq\lambda_1\leq
3n/4;$\\[-2mm]
\item $\displaystyle\chi_\lambda(1)\geq\binom{\lambda_1}{n-\lambda_1}
\chi_{\lambda\setminus\lambda_1}(1),\quad \lambda_1\geq n/2$.
\end{enumerate}
\end{lemma}
\begin{proof} (i) We distinguish the cases $\lambda_1\leq n/4$,
$\lambda_1\geq n/4$ and $\lambda_2\leq n/8$, and $\lambda_1\geq
n/4, \lambda_2\geq n/8$. In the first case it was shown in
\cite[Formula (23)]{MuPuSurf} that $\chi_\lambda(1)\geq
(3/2)^{n/4}\geq 2^{n/8}$. In the second case, for any given
$\lfloor n/8\rfloor$-tuple $(t_1, \ldots, t_{\lfloor n/8\rfloor})$
of 0's and 1's, we can start to deconstruct $\lambda$ by choosing
in the $i$-th step a box from the first row, if $t_i=1$, and from
$\lambda\setminus\lambda_1$, if $t_i=0$. Hence, there are at least
$2^{n/8}$ ways of deconstruction. In the final case, note that
$\chi_\lambda(1)\geq\chi_\mu(1)$ for any partition $\mu$ contained
in $\lambda$; choosing $\mu=(\lfloor n/4\rfloor, \lfloor
n/8\rfloor)$, our claim follows by applying to $\mu$ the argument
used in the second case.\\[2mm]
(ii) This follows as in the proof of \cite[Formula (21)]{MuPuSurf},
observing that the hook product $H[(\lambda_2, \ldots,
\lambda_k)]$ equals
$\frac{(n-\lambda_1)!}{\chi_{\lambda\setminus\lambda_1}(1)}$.
\end{proof}
\subsection{Proof of Theorem 1}
\label{Subsec:ProofThmRoich} The proof of Theorem
\ref{thm:Roichcorrect} makes use of the following two auxiliary
results, whose proofs will be given in the next subsection.
\begin{lemma}
\label{Lem:FewFixed} Let $\mathbf{c}\subseteq S_n$ be a
non-trivial conjugacy class, and let $\pi$ be the element visited
by the random walk $w_\mathbf{c}$ after $3t_c(\mathbf{c})$ steps.
Then, for each $k\geq 1,$ the probability that $\pi$ has more than
$k$ fixed points is bounded above by
\[
\max\left(\frac{2^k}{(k-1)!}, \frac{2^{n/2}}{(\lfloor
n/2\rfloor-1)!}\right).
\]
\end{lemma}
\begin{lemma}
\label{Lem:FewCycles} Let $\mathbf{c}_1, \mathbf{c}_2\subseteq
S_n$ be non-trivial conjugacy classes with $f_1$ respectively
$f_2$ fixed points. For $i=1, 2,$ let $x_i\in\mathbf{c}_i$ be
chosen at random. Then the probability that $x_1$ and $x_2$ have
$l$ common fixed points, is at most
\[
\binom{n}{l}\left(\frac{f_1f_2}{n^2}\right)^l.
\]
Moreover, the probability for $x_1x_2$ to have $k$ cycles on
$\mathrm{supp}(x_1)\cup\mathrm{supp}(x_2)$ is bounded above by
\[
(\log n)^{k-1}/(k-1)!\hspace{.2mm}.
\]
\end{lemma}
The proof of Theorem \ref{thm:Roichcorrect} now proceeds as
follows. Define $g_1:S_n\rightarrow [0, \infty)$ to be the density
of the random walk $w_\mathbf{c}$ after $3t_c(\mathbf{c})$ steps,
and let $g_2$ be the corresponding density after
$6t_c(\mathbf{c})$ steps. Using the fact that $g_1$ is a
class function, we decompose $g_1$ as
\[
g_1=\sum_{\mathbf{c}'}\alpha_{\mathbf{c}'}
\mathbf{1}_{\mathbf{c}'},
\]
and compute $g_2$ in the form
\[
g_2(\pi)=\sum_{\sigma\in S_n} g_1(\sigma) g_1(\pi\sigma^{-1}) =
\sum_{\mathbf{c}_1, \mathbf{c}_2}
\alpha_{\mathbf{c}_1}\alpha_{\mathbf{c}_2} \big|\big\{(c_1,
c_2)\in \mathbf{c}_1\times\mathbf{c}_2: c_1c_2=\pi \big\}\big|.
\]
From Lemmas \ref{Lem:squarevalue}, \ref{Lem:FewFixed}, and
\ref{Lem:FewCycles} we deduce that
\begin{eqnarray*}
|\langle g_2, \chi_\lambda\rangle| & \leq & \frac{1}{n!}
\sum_{\mathbf{c}_1, \mathbf{c}_2} \alpha_{\mathbf{c}_1}
\alpha_{\mathbf{c}_2} \underset{c_2\in\mathbf{c}_2}
{\sum_{c_1\in\mathbf{c}_1}}
(2\hspace{0.3mm}\mathrm{sq}(\lambda))^{\mbox{\scriptsize$\#$
cycles of $c_1c_2$}}\\
 & = & \frac{1}{n!}\sum_{k, l}
(2\hspace{.3mm}\mathrm{sq}(\lambda))^{k+l} \sum_{\mathbf{c}_1,
\mathbf{c}_2}\alpha_{\mathbf{c}_1}\alpha_{\mathbf{c}_2}|\mathbf{c}_1|
|\mathbf{c}_2| P_{k, l}(\mathbf{c}_1, \mathbf{c}_2)\\
 & \leq & \frac{1}{n!} \sum_{k, l}(2\hspace{.3mm}\mathrm{sq}(\lambda))^{k+l}
\sum_{f_1, f_2\geq l}\underset{s_1(\mathbf{c}_2)=f_2}
{\underset{s_1(\mathbf{c}_1)=f_1}{\sum_{\mathbf{c}_1,
\mathbf{c}_2}}} \alpha_{\mathbf{c}_1}\alpha_{\mathbf{c}_2}
|\mathbf{c}_1| |\mathbf{c}_2| \binom{n}{l}
\left(\frac{f_1f_2}{n^2}\right)^l
\frac{(\log n)^{k-1}}{(k-1)!}\\
 & = & \frac{1}{n!}\sum_{k, l}\sum_{f_1, f_2\geq l}(2\hspace{.3mm}\mathrm{sq}(\lambda))^{k+l}
\binom{n}{l}\left(\frac{f_1f_2}{n^2}\right)^l \frac{(\log
n)^{k-1}}{(k-1)!}\\
&&\hspace{4cm}\times\hspace{1mm}
\Bigg(\underset{s_1(\mathbf{c}_1)=f_1}{\sum_{\mathbf{c}_1}}
\alpha_{\mathbf{c}_1}|\mathbf{c}_1|\Bigg)
\Bigg(\underset{s_1(\mathbf{c}_2)=f_2}{\sum_{\mathbf{c}_2}}
\alpha_{\mathbf{c}_2}|\mathbf{c}_2|\Bigg)\\
 & \leq & \frac{1}{n!}\sum_{k, l}\sum_{f_1, f_2\geq l}
(2\hspace{.3mm}\mathrm{sq}(\lambda))^{k+l} \binom{n}{l}
\left(\frac{f_1f_2}{n^2}\right)^l \frac{(\log n)^{k-1}}{(k-1)!}\\
&&\hspace{2.5cm}\times\hspace{1mm} \max\Big(\frac{f_1
2^{f_1}}{f_1!}, \frac{2^{n/2}}{(\lfloor n/2\rfloor-1)!}\Big)
\max\Big(\frac{f_2 2^{f_2}}{f_2!}, \frac{2^{n/2}}{(\lfloor
n/2\rfloor-1)!}\Big).
\end{eqnarray*}
Here, $P_{k, l}(\mathbf{c}_1, \mathbf{c}_2)$ denotes the
probability that, for $x_1$ and $x_2$ chosen at random from
$\mathbf{c}_1$ respectively $\mathbf{c}_2$, $x_1$ and $x_2$
have $l$ common fixed points, and the product $x_1x_2$ has
precisely $k$ cycles on the remaining $n-l$ points. If $f_1$ is in
the interval $[l, n/2]$, increasing $f_1$ by 1 changes the value
of a summand by a factor
\[
2\left(\frac{f_1+1}{f_1}\right)^{l+1}\frac{1}{f_1} \leq
\frac{2e}{f_1},
\]
while each summand increases with $f_1$ in the range $n/2\leq
f_1\leq n$. The same is true for $f_2$; hence, by symmetry, we
obtain
\begin{eqnarray*}
n!\hspace{.3mm}|\langle g_2, \chi_\lambda\rangle| & \leq & 80 n^2 \sum_{k,
l}(2\hspace{.3mm}\mathrm{sq}(\lambda))^{k+l} \binom{n}{l}
\left(\frac{l^2}{n^2}\right)^l \frac{(\log n)^{k-1}}{(k-1)!}
\left(\frac{2^l}{(l-1)!}\right)^2\\
 && + 80n^2 \sum_{k,
l}(2\hspace{.3mm}\mathrm{sq}(\lambda))^{k+l} \binom{n}{l}
\frac{(\log n)^{k-1}}{(k-1)!} \left(\frac{2^{n/2}}{(\lfloor
n/2\rfloor-1)!}\right)^2.
\end{eqnarray*}
For $n\rightarrow\infty$, the second sum tends to zero; thus,
applying Stirling's formula,
\begin{eqnarray*}
n!\hspace{.3mm}|\langle g_2, \chi_\lambda\rangle| & \leq & 1 + 80 n^2 \sum_{k,
l} \frac{(2\hspace{.3mm}\mathrm{sq}(\lambda) \log n)^k}{(k-1)!}
\cdot \frac{(8\hspace{.3mm}\mathrm{sq}(\lambda)n
l^2)^l}{l!((l-1)!)^2n^{2l}}\\
 & \leq & 1 + 80 n^5 \sum_{k, l}
\frac{(2\hspace{.3mm}\mathrm{sq}(\lambda) \log n)^k} {k!} \cdot
\left(\frac{8e^3\hspace{.3mm}\mathrm{sq}(\lambda)} {ln}\right)^l.
\end{eqnarray*}
Since $\mathrm{sq}(\lambda)\leq\sqrt{n}$, the second factor tends
to zero, as $n$ tends to infinity, while the summation over $k$
yields an exponential series. Hence, we obtain the bound
\[
n!\hspace{.3mm}|\langle g_2, \chi_\lambda\rangle| \leq
n^{2\hspace{.3mm}\mathrm{sq}(\lambda)+7}.
\]
On the other hand, from Lemma \ref{Lem:squaredegree}, we deduce
the bound $\chi_\lambda(1)\geq
(\mathrm{sq}(\lambda)/e)^{\mathrm{sq}(\lambda)^2}$, and for
$\chi_\lambda(1)>e^{3(\log n)^4}$, we deduce that
\[
n!\hspace{.3mm}|\langle g_2, \chi_\lambda\rangle| \leq \chi_\lambda(1)^{1/(\log
n)}.
\]
By Lemma~\ref{Lem:FourierGen} (ii), we have
\[
n!\hspace{.3mm}|\langle g_2, \chi_\lambda\rangle| =
n!\frac{\big|\big\langle\frac{1}{|c|}\mathbf{1}_\mathbf{c},
\chi_\lambda\big\rangle\big|^{6t_c(\mathbf{c})}
(n!)^{6t_c(\mathbf{c})-1}} {\chi_\lambda(1)^{6t_c(\mathbf{c})-1}}
= \frac{|\chi_\lambda(\mathbf{c})|^{6t_c(\mathbf{c})}}
{\chi_\lambda(1)^{6t_c(\mathbf{c})-1}},
\]
and together with our estimate for $|\langle g_2,
\chi_\lambda\rangle|$ we obtain
\[
|\chi_\lambda(\mathbf{c})|\leq (\chi_\lambda(1))^{1-
\frac{1-1/(\log n)}{6t_c(\mathbf{c})}}.
\]
Before establishing the upper bound for
$|\chi_\lambda(\mathbf{c})|$ for characters associated to
partitions $\lambda$ satisfying $n-\lambda_1\leq 3\log^4 n$, we
prove the estimate for $t_c(\mathbf{c})$. Let $\xi_k$ be the
random variable which for permutations $\pi_1, \ldots, \pi_k$
chosen independently at random from $\mathbf{c}$ counts the
number of points fixed by all of them. The probability that $1$ is
fixed by all of the permutations equals
$\left(\frac{s_1(\mathbf{c})}{n}\right)^k$, whereas the
probability that both 1 and 2 are fixed by all the permutations
equals
\[
\left(\frac{s_1(\mathbf{c}) (s_1(\mathbf{c})-1)}{n(n-1)}\right)^k.
\]
Hence, $E\xi_k=n\left(\frac{s_1(\mathbf{c})}{n}\right)^k$, and
\begin{equation*}
E(\xi_k^2) =  n\left(\frac{s_1(\mathbf{c})} {n}\right)^k +
n(n-1)\left(\frac{s_1(\mathbf{c})(s_1(\mathbf{c})-1)}
{n(n-1)}\right)^k
  \leq  E\xi_k + (E\xi_k)^2.
\end{equation*}
Using the fact that $\xi_k$ takes only integral values, we deduce
the inequality
\[
1-E\xi_k\leq P(\xi_k=0)\leq
1-E\xi_k+\frac{1}{2}\big(E\xi_k\big)^2;
\]
hence, we obtain for $t_c(\mathbf{c})$ the bounds
\[
\min\left\{k\in2\NN: \left(\frac{s_1(\mathbf{c})} {n}\right)^k \leq
\frac{1}{n(n-1)}\right\} \leq t_c(\mathbf{c}) \leq
\min\left\{k\in2\NN: \left(\frac{s_1(\mathbf{c})} {n}\right)^k \leq
\frac{1}{n^2}\right\}.
\]
Since $s_1(\mathbf{c})\leq n-2$ for every non-trivial class, the
solutions of the equations
\[
\left(\frac{s_1(\mathbf{c})} {n}\right)^k  =
\frac{1}{n(n-1)}\quad\mbox{and}\quad\left(
\frac{s_1(\mathbf{c})}{n}\right)^k = \frac{1}{n^2}
\]
differ by less than 1. Solving for $k$ gives our claim.

Now we bound $\chi_\lambda(\mathbf{c})$ using
Lemma~\ref{Lem:LargeL1CharEst}. We have
\[
|\chi_\lambda(\mathbf{c})|\leq\chi_{\lambda\setminus\lambda_1}(1)
\sum_{a+2b\leq n-\lambda_1}\binom{s_1(\mathbf{c})}{a}
\binom{\lfloor n/2\rfloor}{b}.
\]
Assume first that $s_1(\mathbf{c})\leq n^{2/3}$. Then, using
Lemma~\ref{Lem:chi1lower} (ii),
\[
|\chi_\lambda(\mathbf{c})|\leq (n-\lambda_1)!^{1/2} \sum_{a+2b\leq
n-\lambda_1} n^{2a/3+b} \leq 2(n-\lambda_1)!^{3/2}\binom{n}
{\lfloor2(n-\lambda_1)/3\rfloor}\leq\chi_\lambda(1)^{2/3+\varepsilon},
\]
which is sufficiently small, since $t_c(\mathbf{c})\geq 2$. If, on
the other hand, $s_1(\mathbf{c})>n^{2/3}$, replacing $b$ by $b-1$
and $a$ by $a+2$ changes the value of a summand by
\[
\frac{(s_1(\mathbf{c})-a)(s_1(\mathbf{c})-a-1)b}{(a+1)(a+2)(n/2-b+1)}
> n^{1/4},
\]
and, if $a+2b<n-\lambda_1$, replacing $a$ by $a+1$ changes the
value of a summand by $\frac{s_1(\mathbf{c})-a}{a+1}>n^{1/2}$,
hence, for $n$ sufficiently large, the whole sum over $a$ and $b$
is at most twice its greatest term. Again using
Lemma~\ref{Lem:chi1lower}, we deduce from this that
\[
|\chi_\lambda(\mathbf{c})|\leq
2\chi_{\lambda\setminus\lambda_1}(1)\binom{s_1(\mathbf{c})}{n-\lambda_1}
\leq 2\chi_\lambda(1) \left(\frac{s_1(\mathbf{c})}
{n}\right)^{n-\lambda_1} < 2\chi_\lambda(1)^{1-\frac{\log
(n/s_1(\mathbf{c}))}{\log n}},
\]
which is again sufficiently small by the lower bound for
$t_c(\mathbf{c})$.\\[3mm]
We now sketch the computations needed to show that Theorem~\ref{thm:walks}
holds in fact for all $n\geq 4000$. In the form given above, the proof
only applies to $n\geq e^{40}$. First, in the deduction of Theorem~\ref{thm:walks}
from Theorem~\ref{thm:Roichcorrect}, we can handle the characters $\chi_{(n-1, 1)}$,
$\chi_{(2, 1, 1, \ldots, 1)}$ seperately, noting that for $n\geq 4000$ we have
\[
{\sum_{\lambda}}^* \chi_\lambda(1)^{-0.7}<\frac{1}{n}
\]
where the summation is extended over all partitions $\lambda\vdash n$ apart
from $(n)$, $(n-1, 1)$, $(2, 1, 1, \ldots, 1)$ and $(1, 1, \ldots, 1)$. Following
through the proof of Theorem~\ref{thm:Roichcorrect}, we find that the contribution
of all characters $\chi_\lambda$ with $\lambda_1\leq 3n/4$ of $n-\lambda_1\leq n^{1/3}$
is sufficiently small. To close this gap, we estimate $|\langle g_2, \chi_\lambda\rangle|$
as above, but use Lemma~\ref{Lem:LargeL1CharEst} instead of Lemma~\ref{Lem:squarevalue},
that is, the bound for $|\chi_\lambda(\mathbf{c})|$ now reads
\[
\min\left(\chi_\lambda(1), \max_{\nu\leq n-\lambda_1}(2\sqrt{n-\lambda_1})^\nu\binom{k+l}{\nu}\right)
\quad\mbox{instead of }(2\sq(\lambda))^{k+l};
\]
as a result, we obtain a bound of sufficient quality for all characters $\chi_\lambda$ satisfying
$\chi_\lambda(1)>2e^{\frac{10}{19}\log^2 n+\frac{20}{19}\log n}$; and we conclude the proof
by noting that for $n\geq 4000$, the degree of a character $\chi_\lambda$ with $\|\lambda\|\leq\lambda_1<n-n^{1/3}$
is larger than this bound.\\[3mm]
Note that in the intermediate range, we established a bound somewhat weaker than Theorem~\ref{thm:Roichcorrect};
in particular, we do not claim that Theorem~\ref{thm:Roichcorrect} holds for all $n\geq 4000$.
However, it may well be true that both theorems are in fact true for all integers $n$ without any exception.
\subsection{Proof of Lemmas~\ref{Lem:FewFixed}
and~\ref{Lem:FewCycles}} \label{Subsec:ProofLemRoich} To complete
the proof of Theorem~\ref{thm:Roichcorrect}, it remains to
establish Lemmas \ref{Lem:FewFixed} and \ref{Lem:FewCycles}.
\begin{proofof}{Lemma~{\em \ref{Lem:FewFixed}}} Let $a$, $b$ be integers, and let $\pi_1,
\ldots, \pi_{3t_c(\mathbf{c})}\in\mathbf{c}$ be elements chosen at
random. Denote by $P(a, b)$ the probability that the points $1,
\ldots, a$ are fixed by all the $\pi_i$, and that, for each
$\beta$ with $a+1\leq\beta\leq a+b$, there is some $i$ such that
$\pi_i$ moves $\beta$, while the product $\pi_1\cdots
\pi_{3t_c(\mathbf{c})}$ fixes $\beta$. We have
\[
P(a, 0) = P(\pi_1\mbox{ fixes }1, \ldots, a)^{3t_c(\mathbf{c})}
\leq P(\pi_1\mbox{ fixes }1)^{3at_c(\mathbf{c})} =
\Bigg(\left(\frac{s_1(\mathbf{c})}
{n}\right)^{t_c(\mathbf{c})}\Bigg)^{3a}.
\]
By the definition of $t_c(\mathbf{c})$, we have
$\left(\frac{s_1(\mathbf{c})}
{n}\right)^{t_c(\mathbf{c})}<\frac{1}{n}$; hence, $P(a, 0)$ is
bounded by $n^{-3a}$. Next, consider $P(0, b)$. The product
$\pi_1\cdots \pi_{3t_c(\mathbf{c})}$ fixes $\beta$ if and only if
$$\big(\pi_1\cdots \pi_{3t_c(\mathbf{c})-1}\big)(\beta) =
\pi_{3t_c(\mathbf{c})}^{-1}(\beta).$$
Let
$h\in\mbox{Sym}\big([n]\setminus\{(\pi_1\cdots
\pi_{3t_c(\mathbf{c})-1})(\gamma): \gamma\leq b,
\gamma\neq\beta\}\big)$ be chosen at random. Then replacing
$\pi_{3t_c(\mathbf{c})}$ by $\pi_{3t_c(\mathbf{c})}^h$ does not
alter $(\pi_1\cdots \pi_{3t_c(\mathbf{c})-1})(\gamma)$ for
$\gamma\leq b, \gamma\neq\beta$, while $\big(\pi_1\cdots
\pi_{3t_c(\mathbf{c})-1}\big)(\beta) =
\big(\pi_{3t_c(\mathbf{c})}^h\big)^{-1}(\beta)$ holds with
probability $\frac{1}{n-b-1}$ or 0. Since the $\pi_i$ are chosen
from a conjugacy class, conjugating with a random element from
some subgroup does not affect the equidistribution of the $\pi_i$,
hence, we obtain $P(0, b)\leq (n-b+1)^b$. Finally, a permutation
$\pi\in S_n$ that fixes the points $1, \ldots, a$ can be viewed as
a permutation on $n-a$ elements, thus $P(a, b)\leq n^{-3a}
(n-a-b+1)^{-b}$. If $a+b\leq n/2$, we deduce $P(a,
b)\leq\frac{2^b}{n^{3a+b}}$, and, summing over all pairs $a, b$
with $a+b\geq k$ we obtain our claim. Note that $P$ is decreasing
in both $a$ and $b$; hence, if $a+b>n/2$, we may replace the pair
$(a, b)$ by some pair $(a', b')$ satisfying $a'+b'=\lfloor
n/2\rfloor$ and use our estimate for the latter pair. Since the
probability that there exist $k$ points which are fixed by the
product $\pi_1\cdots \pi_{3t_c(\mathbf{c})}$ is at most
$\binom{n}{k}$ times the probability that $\pi_1\cdots
\pi_{3t_c(\mathbf{c})}$ fixes the points $1, \ldots, k$, our claim
follows.
\end{proofof}
\begin{proofof}{Lemma~{\em \ref{Lem:FewCycles}}}
The probability that both $x_1$ and $x_2$ fix 1 equals
$\frac{f_1f_2}{n^2}$, and the conditional probability
\[
P\big(\mbox{$x_1$ and $x_2$ fix 1$\,|\,\exists a_1, \ldots a_k:\,a_i\neq1,\,x_1$ and
$x_2$ fix $a_i, \forall i\leq k$}\big)
\]
is smaller, thus, the first claim of the lemma follows. Now let
$x_1, x_2$ be chosen at random from $\mathbf{c}_1$ and $\mathbf{c}_2$,
respectively. Assume that not both of them fix 1, without loss we
assume that $x_1(1)\neq 1$. Then 1 lies in a cycle of length $i$
of $x_1x_2$ if and only if
\begin{equation}\label{Eq:icyclecondition}
\big(x_2(x_1x_2)^{i-1}\big)(1)=x_1^{-1}(1)
\end{equation}
with $i$ chosen minimal among all positive integers with this
property. Choose an element
$h\in\mbox{Sym}\big([n]\setminus\{x_2(x_1x_2)^j(1): 0\leq j\leq
i-2\}\big)$ at random, and replace $x_1$ by $x_1^h$. Then
(\ref{Eq:icyclecondition}) becomes true with probability
$\frac{1}{n-i+1}$. Increasing $i$ until (\ref{Eq:icyclecondition})
happens to hold, we obtain one cycle of $x_1x_2$. Next, choose
some point outside this cycle, and repeat the procedure, where $h$
is to be chosen in such a way that $h$ fixes all points in all
cycles already determined as well as the points already
constructed in the current cycle. In this way, we obtain the
number $c$ of cycles of $x_1x_2$ as the value returned by the
following\\[3mm]
{\bf Stochastic Algorithm}
\vspace{-3mm}
\begin{enumerate}
\item Set $m:= n, i:= 0$ and $c:= 0$.\\[-3mm]
\item Increase $i$ by 1.\\[-3mm]
\item With probability $1-\frac{1}{m-i+1}$, go to (ii); otherwise
continue with (iv).\\[-3mm]
\item Set $m:= m-i$, $c:=c+1$ and $i:=0$.\\[-3mm]
\item If $m=0$, stop and return $c$; otherwise go to (ii).
\end{enumerate}
Let $P$ be the probability that $x_1x_2$ has $k$ cycles, and that
their lengths are $c_1, \ldots c_k$. Then, in step (iii) of the
algorithm, the second possibility was chosen $k$ times, and the
probabilities were $\frac{1}{n-c_1+1}, \frac{1}{n-c_1-c_2+1},
\ldots, \frac{1}{n-c_1-\ldots-c_k+1}$, respectively. Hence, $P$ is
bounded above by the product of these probabilities; and, writing
$i_j:=c_{k-j+1}+\dots+c_k+1$, we obtain
\begin{eqnarray*}
P(x_1x_2\mbox{ has $k$ cycles}) & \leq &
\sum_{1=i_1<i_2<\dots<i_k\leq n} \prod_{j=1}^k\frac{1}{i_j}\\[1mm]
 & \leq & \frac{1}{(k-1)!}\sum_{2\leq i_2, \ldots, i_k\leq n}
\prod_{j=1}^k\frac{1}{i_j}\\[1mm]
 & = & \frac{1}{(k-1)!}\left(\sum_{2\leq
i\leq n} \frac{1}{i}\right)^{k-1}\\[1mm]
 & \leq & \frac{(\log n)^{k-1}}{(k-1)!},
\end{eqnarray*}
which proves our claim.
\end{proofof}
\section{Character estimates for elements of prescribed order}
\label{Sec:MainLemFuchs}
As it stands, Theorem~\ref{thm:Roichcorrect} is not strong enough
to obtain an asymptotic estimate for the subgroup growth of
Fuchsian groups. When combined with the methods of
Subsection~\ref{subsec:MainLemFuchs}, it could be used to
determine the asymptotics of $s_n(\Gamma)$ for certain Fuchsian
groups $\Gamma$, namely those given by a presentation of the form
\[
\Gamma=\Big\langle x_1, \ldots, x_r\,\big|\,
x_1^{a_1}=x_2^{a_2}=\cdots=x_r^{a_r}=x_1 x_2\cdots x_r=1\Big\rangle
\]
where $a_i\geq 2$ and $r\geq 73$. Unfortunately, this would
exclude all better known examples in this class. From the point of view of  an application
to the subgroup growth of Fuchsian groups, the weakness of
Theorem~\ref{thm:Roichcorrect} is caused by its generality.
Combining instead an estimate by Fomin and Lulov \cite{FoLu} for
character values $\chi_\lambda(\pi)$ where all cycle lengths of
$\pi$ are equal, with combinatorial arguments plus the estimates
of Subsection~\ref{Subsec:CharElemEst}, we shall derive the
following sharper estimate.
\begin{proposition}
\label{prop:MainLemFuchs} Let $q\geq 2$ be an integer, $\varepsilon>0,$
and let $n$ be sufficiently large. Then, for every partition
$\lambda\vdash n,$ we have
\begin{equation}\label{Eq:MainLemFuchs}
\underset{\pi^q=1}{\sum_{\pi\in S_n}} |\chi_\lambda(\pi)| \leq
\big(\chi_\lambda(1)\big)^{\frac{1}{q}+\varepsilon} |\Hom(C_q, S_n)|.
\end{equation}
\end{proposition}
We begin with some results concerning the statistical distribution
of permutations in $S_n$ and in wreath products, which will be
used in the proof of Proposition~\ref{prop:MainLemFuchs} as well
as in the next section.
\subsection{Statistics of the symmetric group}
\label{Subsec:Statistics}
For integers $n\geq 1$ and $q\geq 2$, define $N(n, q)$ to be the
number of elements $\pi\in S_n$ with $\pi^q=1$. Furthermore, for
integers $s_t$ with $t|q$ and $t<q$, define $N(n, q, s_1, \ldots,
s_T)$, to be the number of elements $\pi\in S_n$ with $\pi^q=1$
and $s_t(\pi)=s_t$ for all $t|q,t<q$.
\begin{lemma}\label{Lem:statUppbound}
Let $q\geq 2$ be an integer, $n$ sufficiently large, and let $s_1,
\ldots, s_T$ be given in such a way that $n\equiv 1\cdot s_1 +
\cdots + T\cdot s_T \pmod{q}$. Then we have the estimate
\begin{equation}
\frac{N(n, q, s_1, \ldots, s_T)}{N(n, q)} \leq q^{\sigma(q)}
\prod\limits_{\underset{s_t>2en^{t/q}}{t}}\left(\frac{en^{t/q}}{ts_t}\right)^{s_t},
\end{equation}
where $\sigma(q)$ denotes the sum of divisors of $q$.
\end{lemma}
\begin{proof}
Put $\mathcal{S}:=\underset{t\neq q}{\sum\limits_{t|q}} ts_t$. We
have
\[
N(n, q, s_1, \ldots, s_T) = \frac{n!} {((n-\mathcal{S})/q)!
s_1!\cdots s_T! 1^{s_1}\cdots T^{s_T}\cdot q^{(n-\mathcal{S})/q}}.
\]
Let
\[
\tilde{s_t}:= \begin{cases}
s_t, & s_t\leq 2en^{t/q},\\[1mm]
s_t\bmod q/t, & s_t>2en^{t/q}
\end{cases}\qquad (t|q,\,t<q),
\]
where $s_t\bmod{q/t}$ takes values in $\{0, 1, \ldots, q/t-1\}$,
and put $\tilde{\mathcal{S}}:=\underset{t\neq
q}{\sum\limits_{t|q}} t\tilde{s}_t$. Then we get
\begin{eqnarray*}
\frac{N(n, q, s_1, \ldots, s_T)}{N(n, q)} & \leq &
\frac{N(n, q, s_1, \ldots, s_T)}{N(n, q, \tilde{s}_1, \ldots, \tilde{s}_T)}\\[1mm]
 & = & \frac{((n-\tilde{\mathcal{S}})/q)!
\tilde{s}_1!\cdots \tilde{s}_T! 1^{\tilde{s}_1} \cdots
T^{\tilde{s}_T}\cdot
q^{(n-\tilde{\mathcal{S}})/q}}{((n-\mathcal{S})/q)! s_1!\cdots
s_T!
1^{s_1} \cdots T^{s_T}\cdot q^{(n-\mathcal{S})/q}}\\[1mm]
 & \leq & \prod_{\underset{s_t>2en^{t/q}}{t}}\frac{n^{ts_t/q} (q/t)!}{s_t! t^{s_t-q/t}}\\[1mm]
 & \leq & \prod_{\underset{s_t>2en^{t/q}}{t}}q^{q/t} \left(\frac{en^{t/q}}{ts_{t}}\right)^{s_{t}}\\[1mm]
 & \leq & q^{\sigma(q)}\prod_{\underset{s_t>2en^{t/q}}{t}}\left(\frac{en^{t/q}}{ts_{t}}\right)^{s_{t}},
\end{eqnarray*}
as claimed.
\end{proof}
For a conjugacy class $\mathbf{c}$ of $S_n$, denote by $C_{S_n}(\mathbf{c})$
the centralizer of $\mathbf{c}$ in $S_n$. Then $C_{S_n}(\mathbf{c})$ is isomorphic
to a direct product of the form
\[
C_{S_n}(\mathbf{c}) \cong \prod\limits_{d|q} C_d\wr S_{s_d(\sigma)},
\]
where $\sigma$ is some element of ${\bf c}$.
\begin{lemma}\label{Lem:CycleBound}
Let $q$ and $l$ be integers, and let $\mathbf{c}\in S_n$ be a conjugacy
class with $\mathbf{c}^q=1$.\vspace{-3mm}
\begin{enumerate}
\item The number of $\pi\in S_n$ with at least $k$ cycles is bounded above by
\[
n!\hspace{.4mm}(\log n)^{k-1}/((k-1)!).
\]
\item The number of $\pi\in C_{S_n}(\mathbf{c})$ with at least $k$ cycles is bounded above by
\[
|C_{S_n}(\mathbf{c})|\left(\frac{3q\log n}{k}\right)^{k/q},
\]
 provided that
$k\geq (\log n)^3$ and $n\geq n_0(q)$.
\end{enumerate}
\end{lemma}
\begin{proof}
(i) Let $\pi\in S_n$ be chosen at random. Then 1 is a fixed point
of $\pi$ with probability $\frac{1}{n}$. If it is not a fixed
point, then it lies in a 2-cycle with probability $\frac{1}{n-1}$,
and, more generally, the conditional probability for 1 to lie in a
cycle of length $c$, provided that it lies in a cycle of length at
least $c$ is $\frac{1}{n-c+1}$. Arguing now as in the proof of
Lemma~\ref{Lem:FewCycles} establishes our claim.\\[3mm]
(ii) Consider a single direct factor $G:=C_d\wr S_{s_d(\sigma)}$ of $C_{S_n}(\mathbf{c})$,
and let $\phi:G\rightarrow S_{s_d(\sigma)}$ be the canonical projection.
Let $c$ be a cycle in $G$. Then, the projection $\overline{c}$ of $c$ in $S_{s_d(\sigma)}$
is a cycle, too, and there are at most $d$ cycles $c_1, \ldots, c_d\in G$ which have
the same image in $S_{s_d(\sigma)}$. We deduce that the probability that a permutation $\pi$,
chosen at random in $G$, has $k$ cycles is at most the probability that a permutation
chosen at random in $S_{s_d(\sigma)}$ has $\lceil k/d\rceil$ cycles. Together with part (i) of
this lemma, we obtain
\begin{eqnarray*}
\frac{1}{|C_{S_n}(\mathbf{c})|}\hspace{.4mm}\Big|\Big\{\pi\in C_{S_n}(\mathbf{c}):\hspace{1mm}
|\mbox{Orbits}(\pi)|\geq k\Big\}\Big| & \leq &
\sum_{\sum_{d|q}d\kappa_d=k}\prod_{d|q} \min\left(1,
\frac{(\log n)^{\kappa_d-1}}{(\kappa_d-1)!}\right)\\[1mm]
 & \leq & k^{\tau(q)} \max\limits_{\sum_{d|q}d\kappa_d=k}\prod_{d|q}
\min\left(1, \frac{(\log n)^{\kappa_d-1}}{(\kappa_d-1)!}\right)\\[1mm]
 & \leq & k^{\tau(q)} \max\limits_{\sum_{d|q}d\kappa_d=k}\underset{\kappa_d\geq 3\log n}{\prod_{d|q}}
\frac{(\log n)^{\kappa_d-1}}{(\kappa_d-1)!}.
\end{eqnarray*}
If we replace $\kappa_d$ by $\kappa_d+\frac{q}{d}$, and $\kappa_q$ by
$\kappa_q-1$, a single summand is changed by a factor
\[
\frac{(\log n)^{\frac{q}{d}-1}(\kappa_q-1)}{\kappa_d(\kappa_d+1)\cdots(\kappa_d+q/d-1)},
\]
which is less than 1, provided that $\kappa_d>(\log n)\sqrt{\kappa_q}$. Hence, the maximum
is attained for some tuple $(\kappa_1, \ldots, \kappa_q)$ satisfying $\kappa_q\geq k/q-q\sqrt{k}\log n$.
From this we deduce
\begin{eqnarray*}
\frac{1}{|C_{S_n}(\mathbf{c})|}\hspace{.4mm}\Big|\Big\{\pi\in C_{S_n}(\mathbf{c}):
|\mbox{Orbits}(\pi)|\geq k\Big\}\Big| & \leq &
k^{\tau(q)}\frac{(\log n)^{k/q- q\sqrt{k}(\log n)-1}} {\lfloor
k/q- q\sqrt{k}(\log n)-1\rfloor!}\\[1.5mm]
 & \leq & \left(\frac{3q\log n}{k}\right)^{k/q},
\end{eqnarray*}
provided that $k\geq (\log n)^3$ and $n\geq n_0(q)$.
\end{proof}
\begin{lemma}
\label{Lem:Poisson}
Let $\pi\in S_n$ be chosen at random, and let $d, d_1, d_2$ be positive integers.\vspace{-3mm}
\begin{enumerate}
\item As $n\rightarrow\infty,$ the distribution of $s_d(\pi)$ converges to
a Poisson distribution with mean $\frac{1}{d},$ and we have
\[
\frac{1}{n!}\sum_{\pi\in S_n} \big(s_d(\pi)\big)^q \rightarrow \sum_{\nu=1}^q d^{-\nu} S(q, \nu),\quad n\rightarrow\infty,
\]
where the $S(q, \nu)$ are Stirling numbers of the second kind.\\[-3mm]
\item As $n\rightarrow\infty,$ the random variables $s_{d_1}(\pi)$ and $s_{d_2}(\pi)$
are asymptotically independent.
\end{enumerate}
\end{lemma}
\begin{proof}
(i) Let $P(d, k)$ be the probability that for $\pi$ chosen at random from $S_n$,
$\pi$ contains the $d$-cycles $(12\ldots d), (d+1\ldots 2d), \ldots, ((k-1)d+1\ldots kd)$.
Then, as $P(d, k)=\frac{(n-kd)!}{n!}$, we have for $k\leq n/d$,
\[
\frac{1}{n!}\sum_{\pi\in S_n}\binom{s_d(\pi)}{k} = \frac{n!}{d^k k! (n-kd)!}P(d, k) = \frac{1}{d^k k!}.
\]
On the other hand, for a random variable $\xi$ which has Poisson distribution with mean
$1/d$, we have
\[
\mathbf{E}\binom{\xi}{k} = \sum_{\nu=0}^\infty \binom{\nu}{k}\frac{e^{-1/d}}{d^\nu \nu!}
= \frac{e^{-1/d}}{d^k k!} \sum_{\nu=k}^\infty \frac{1}{(\nu-k)!d^{\nu-k}} = \frac{1}{d^k k!}.
\]
We conclude that the first $\lfloor n/d\rfloor$ moments of $\xi$ and $s_d(\pi)$ coincide,
hence, by the method of moments,\footnote{Cf., for instance, \cite{Elliott}.} the distributions are identical,
proving the first assertion.
Let $\xi$ be a random variable with mean $1/d$ and Poisson distribution. Then $\mathbf{E}\xi^q$ is
the expected number of $q$-multisets in $[\xi]$. A set $S$ of $\nu$ elements gives rise to $S(q, \nu)$
different multisets $M$, such that $M=S$ as sets, and the computation of $\mathbf{E}\binom{\xi}{\nu}$
shows, that the expected number of $\nu$-sets is $d^{-\nu}$, whence the second assertion.\\[3mm]
(ii) For $i\leq t, i\neq d$, fix integers $e_i$. We compute the conditional expectation
\[
\mathbf{E}\left(\binom{s_d(\pi)}{k}\,\Big|\, s_i(\pi)=e_i, i\leq t, \,i\neq d\right) =
\frac{(n-\mathcal{E})!}{d^k k! (n-\mathcal{E}-kd)!}\hspace{.4mm}P(d, k),
\]
where $\mathcal{E}=\sum\limits_{q\neq i\leq t} ie_i$. As $n$ tends to infinity, the right
hand side converges to $\frac{1}{d^k k!}$, proving our claim.
\end{proof}
Our next group of results describes the distribution of cycles in permutations
of prescribed order. The proof makes use of the following purely analytic result.
\begin{lemma}
\label{Lem:Hayman}
Let $P(z)=\sum_{k=1}^q a_k z^k$ be a real polynomial, and let $Q(z)=\sum_{k=1}^q |a_k| z^k$.
Assume that $P(z)\neq\pm Q(\pm z)$, and that $\{k:a_k\neq 0\}$ has greatest common divisor $1$.
Define the sequences $(b_n), (b_n^+)$ by means of the equations
\[
e^{P(z)} = \sum_{n=0}^\infty b_n\,z^n/(n!),\quad
e^{Q(z)} = \sum_{n=0}^\infty b_n^+\,z^n/(n!).
\]
Then we have $|b_n|<b_n^+e^{-cn^{1/q}}$ for some $c>0$ and sufficiently large $n$.
\end{lemma}
\begin{proof}
We first claim that there is some constant $c$, such that for all real
numbers $r$ sufficiently large, and all complex numbers $z$ with $|z|=r$,
we have $\Re P(z)\leq Q(r)-cr$. For otherwise we would have
$\Re a_k z^k \geq |a_k| r^k- cr$, that is, $\big|\frac{z}{|z|}-\zeta\big|<\varepsilon$
for some $(2k)$-th root of unity $\zeta$, and $\varepsilon$ arbitrarily small.
Since by assumption the set $\{k:a_k\neq 0\}$ has greatest common divisor 1, we deduce
that $|\arg z|<\varepsilon$ or $|\arg z-\pi|<\varepsilon$. However, the assumptions $P(z)\neq\pm Q(\pm z)$
 and $a_k\neq0$ for at least one odd $k$ imply, that in these cases $\Re a_k z^k$ was negative for at least one value of $k$.
From this we obtain that $|e^{P(z)}|\leq e^{Q(r)-cr}$ for some $c>0$ and all $z$ with
$|z|=r$. Let $r_n$ be the solution of the equation $rQ'(r)=n$. Then we deduce from
Cauchy's bound that $\beta_n\leq\frac{e^{Q(r_n)-cr_n}}{r_n^n}$, while from \cite[Corollary II]{Hayman}
we obtain the lower bound $\beta_n^+\geq \frac{e^{Q(r_n)}}{r_n^n n^c}$ with some
absolute constant $c$. From these bounds and the asymptotics $r_n\sim \sqrt[q]{n/(qa_q)}$
the lemma follows.
\end{proof}

\begin{lemma}
\label{Lem:CycleDistAsymp}
Let $q\geq 2$ and $e_t\geq 0\, (t|q,\, t<q)$ be integers.
\begin{enumerate}
\item There exist constants $\alpha_{e_1, \ldots, e_T}^{(d)},$ such that,
for all $n\geq 1,$
\begin{equation}
\label{Eq:CycledistAsymp}
\sum_{\pi^q=1}\underset{t<q}{\prod_{t|q}} \big(s_t(\pi)\big)^{e_t} =
\sum_{\nu=0}^D \alpha_{e_1, \ldots, e_T}^{(\nu)}\hspace{.4mm}\frac{n!}{(n-\nu)!}\hspace{.4mm}|\Hom(C_q, S_{n-\nu})|,
\end{equation}
where $D=\sum_t te_t$. The coefficients $\alpha_{e_1, \ldots, e_T}^{(d)}$
are recursively determined by means of the equations
\[
\alpha_{e_1, \ldots, e_q}^{(d)} = \begin{cases}\displaystyle\frac{1}{d}\underset{t<q}{\sum_{t|q}}
\sum_{\nu=1}^{e_t}\binom{e_t}{\nu}
\alpha_{e_1, \ldots, e_t-\nu, \ldots, e_T}^{(d-t)},& d\geq 1\\[2mm] 1, & d=0.
\end{cases}
\]
\item We have
\[
\sum_{\pi^q=1}\underset{t<q}{\prod_{t|q}} \big(s_t(\pi)\big)^{e_t} = \Big(1+\mathcal{O}(n^{-1/q})\Big)
\alpha_{e_1, \ldots, e_T}^{D} n^{D/q}.
\]
\item If $q$ is even, then there exists a constant $c>0$ such that, for $n$ sufficiently large,
\begin{equation}
\label{Eq:AlterCycledist}
\sum_{\pi^q=1}\underset{t<q}{\prod_{t|q}}(-1)^{(t-1)s_t(\pi)} \big(s_t(\pi)\big)^{e_t}  <
e^{-cn^{1/q}}|\Hom(C_q, S_{n-\nu})|.
\end{equation}
\end{enumerate}
\end{lemma}
\begin{proof} (i) It will be convenient to allow $t$ to run over all divisors of $q$, setting $e_q:=0$.
Abbreviate the left-hand side of (\ref{Eq:CycledistAsymp}) as
$\mathcal{S}_{e_1, \ldots, e_q}(n)$, and let $\pi$ be a permutation in $S_n$ with
$\pi^q=1$. Then $n$ occurs in some cycle of $\pi$ of length $t$ for some $t|q$.
Let $\pi'\in S_{n-t}$ be the permutation resulting from $\pi$ by deleting the cycle
containing $n$. We have $s_d(\pi')=s_d(\pi)$ for $d\neq t$, and $s_t(\pi')=s_t(\pi)-1$,
that is,
\[
{\prod_{d|q}} \big(s_d(\pi)\big)^{e_d} = \big(s_t(\pi')+1\big)^{e_t}\prod_{d|q} \big(s_d(\pi')\big)^{e_d}.
\]
The $t$-cycle containing $n$ can be chosen in $\frac{(n-1)!}{(n-t)!}$ ways; hence, we obtain for $n\geq 1$
the recursion formula
\begin{equation}
\label{Eq:CycleDisRec}
\mathcal{S}_{e_1, \ldots, e_q}(n) = \sum_{t|q} \frac{(n-1)!}{(n-t)!} \sum_{\nu=0}^{e_t}
\binom{e_t}{\nu} \mathcal{S}_{e_1, \ldots, e_t-\nu, \ldots, e_q}(n-t),
\end{equation}
where
\[
S_{e_1, \ldots, e_q}(0) = \begin{cases} 1, & (e_1, \ldots, e_q)=(0, \ldots, 0)\\ 0, & \mbox{otherwise},
\end{cases}
\]
and $S_{e_1, \ldots, e_q}(n)=0$ if one of the $e_t$ or $n$ is negative.
Introducing the exponential generating functions
\[
\mathfrak{S}_{e_1, \ldots, e_q}(z) = \sum_{n=0}^\infty \mathcal{S}_{e_1, \ldots, e_q}(n)\, z^n/(n!),
\]
multiplying (\ref{Eq:CycleDisRec}) by $\frac{z^{n-1}}{(n-1)!}$, and summing over $n\geq 1$,
this recurrence relation translates into the differential equation
\[
\mathfrak{S}'_{e_1, \ldots, e_q}(z) - \bigg(\sum_{t|q} z^{t-1}\bigg)\mathfrak{S}_{e_1, \ldots, e_q}(z)
= \sum_{t|q}z^{t-1}\sum_{\nu=1}^{e_t}\binom{e_t}{\nu} \mathfrak{S}_{e_1, \ldots, e_t-\nu, \ldots, e_q}(z).
\]
Integrating the latter equation, we find that
\begin{multline}
\label{Eq:CycleDiffEq}
\mathfrak{S}_{e_1, \ldots, e_q}(z) = \exp\bigg(\sum_{t|q}\frac{z^t}{t}\bigg) \sum_{t|q}\sum_{\nu=1}^{e_t}
\binom{e_t}{\nu}\int_0^z \bigg(\zeta^{t-1}\exp\bigg(-\sum_{t|q}\frac{\zeta^t}{t}\bigg)
\mathfrak{S}_{e_1, \ldots, e_t-\nu, \ldots, e_q}(\zeta)\Big)\;d\zeta\\
+ \exp\bigg(\sum_{t|q}\frac{z^t}{t}\bigg),
\end{multline}
where the value of the integration constant has been determined by comparing
the coefficients of $z$. We claim that there exist polynomials $P_{e_1, \ldots, e_q}(z)$, such that
\begin{equation}
\label{Eq:CycleDiffSol}
\mathfrak{S}_{e_1, \ldots, e_q}(z) = P_{e_1, \ldots, e_q}(z) \exp\bigg(\sum_{t|q}\frac{z^t}{t}\bigg).
\end{equation}
The proof is by induction on $\mathbf{e}=\sum_{t|q} e_t$. If $e_t=0$ for all $t$,
then $\mathcal{S}_{0, \ldots, 0}(n)=|\Hom(C_q, S_n)|$, that is\footnote{Cf.,
for instance, \cite[Prop.~1]{DM}.},
$\mathfrak{S}_{0, \ldots, 0}(z)=\exp\big(\sum_{t|q}\frac{z^t}{t}\big)$, and
(\ref{Eq:CycleDiffSol}) holds with $P_{0, \ldots, 0}(z)=1$.
Suppose now that our claim holds for all tuples $(e_1', \ldots, e_q')$ with
$\sum_{t|q} e_t'<\mathbf{e}$, and some $\mathbf{e}\geq 1$, and let $(e_1, \ldots, e_q)$
be a tuple with $\sum_{t|q} e_t=\mathbf{e}$. Inserting (\ref{Eq:CycleDiffSol})
into the right-hand side of (\ref{Eq:CycleDiffEq}), we find that
\[
\mathfrak{S}_{e_1, \ldots, e_q}(z) = \exp\bigg(\sum_{t|q}\frac{z^t}{t}\bigg)
\sum_{t|q}\sum_{\nu=1}^{e_t}\binom{e_t}{\nu}\int_0^z \zeta^{t-1}
P_{e_1, \ldots, e_t-\nu, \ldots, e_q}(\zeta)\;d\zeta
+ \exp\bigg(\sum_{t|q}\frac{z^t}{t}\bigg),
\]
that is, (\ref{Eq:CycleDiffSol}) holds for $(e_1, \ldots, e_q)$ with
\[
P_{e_1, \ldots, e_q}(z) = 1+\sum_{t|q}\sum_{\nu=1}^{e_t}\binom{e_t}{\nu}\int_0^z \zeta^{t-1}
P_{e_1, \ldots, e_t-\nu, \ldots, e_q}(\zeta)\;d\zeta.
\]
Comparing coefficients in (\ref{Eq:CycleDiffSol}) and in the recurrence relation
determining the polynomials $P_{e_1, \ldots, e_q}(z)$, the assertions of the lemma follow.\\[3mm]
(ii) By \cite[Eq. (22)]{MInvent}, we have
\[
\frac{|\Hom(C_q, S_n)\|}{|\Hom(C_q, S_{n-k})|} = \big(1+\mathcal{O}(n^{-1/q})\big)n^{k(1-1/q)}.
\]
Together with part (i), our claim follows.\\[3mm]
(iii) Denote the left-hand side of (\ref{Eq:AlterCycledist}) by $\mathcal{S}_{e_1, \ldots, e_q}^*(n)$.
In this case, we obtain the recurrence relation
\[
\mathcal{S}_{e_1, \ldots, e_q}^*(n) =
\sum_{t|q} (-1)^{t-1}\frac{(n-1)!}{(n-t)!} \sum_{\nu=0}^{e_t}
\binom{e_t}{\nu} \mathcal{S}_{e_1, \ldots, e_t-\nu, \ldots, e_q}^*(n-t).
\]
Arguing as in part (i), the corresponding exponential generating function $\mathfrak{S}_{e_1, \ldots, e_q}^*(z)$
is found to satisfy
\[
\mathfrak{S}_{e_1, \ldots, e_q}^*(z) = P_{e_1, \ldots, e_q}^*(z)
\exp\bigg(\sum_{t|q}(-1)^{t-1}\frac{z^t}{t}\bigg)
\]
with certain polynomials $P_{e_1, \ldots, e_q}^*$.
Our claim follows from this and Lemma~\ref{Lem:Hayman}.
\end{proof}
\begin{lemma}
\label{Lem:CycleAsympExp} Define the polynomials $P_{e_1, \ldots,
e_q}$ as in the proof of Lemma~{\em \ref{Lem:CycleDistAsymp}}.
Then we have
\begin{eqnarray*}
P_{\vec{e}_t}(z) & = & 1\hspace{.6mm}+\hspace{.6mm}\frac{z^t}{t},\\
P_{\vec{e}_{t_1}+\vec{e}_{t_2}}(z) & = & 1\hspace{.6mm}+\hspace{.6mm}\frac{z^{t_1}}{t_1}\hspace{.6mm}+\hspace{.6mm}
\frac{z^{t_2}}{t_2}\hspace{.6mm}+\hspace{.6mm}\frac{t_1+t_2}{(t_1t_2)^2}\hspace{.4mm}z^{t_1+t_2},\\
P_{k\cdot\vec{e}_t}(z) & = & \sum_{\nu=0}^k S(k+1, \nu+1)\hspace{.4mm}
\frac{z^{\nu t}}{t^\nu},
\end{eqnarray*}
where $\vec{e}_t$ denotes the tuple with $e_t=1$ and $e_d=0$ for
$d\neq t$.
\end{lemma}
\begin{proof}
The first two equations follow directly from the definition. The
third equation is established by induction on $k$. For $k=1$ the
claim is already proven. Assuming the result for
$P_{\kappa\cdot\vec{e}_t}(z)$ with $\kappa<k$ and some $k\geq2$,
we find that
\begin{eqnarray*}
P_{k\cdot\vec{e}_t}(z) & = & 1\hspace{.6mm}+\hspace{.6mm}
\sum_{\nu=1}^k\binom{k}{\nu}\int_0^z
\zeta^{t-1}P_{(k-\nu)\cdot\vec{e}_t}(\zeta)\;d\zeta\\
 & = & 1\hspace{.6mm}+\hspace{.6mm}\sum_{\nu=1}^k\binom{k}{\nu} \sum_{\mu=0}^{k-\nu}
S(k-\nu+1, \mu+1)\hspace{.4mm}\frac{z^{(\mu+1)t}}{(\mu+1)t^{\mu+1}}\\
 & = & 1\hspace{.6mm}+\hspace{.6mm}\sum_{\mu=0}^{k-1} \frac{z^{(\mu+1)t}}{t^{\mu+1}}\cdot
\frac{1}{\mu+1}\sum_{\nu=1}^{k-\mu} \binom{k}{\nu}S(k-\nu+1, \mu+1).
\end{eqnarray*}
Hence, our claim would follow from the identity
\[
(\mu+1) S(k+1, \mu+2) = \sum_{\nu=1}^{k-\mu}\binom{k}{\nu}S(k-\nu+1, \mu+1),
\]
which can be seen to hold as follows: the left-hand side counts the number
of partitions of a set with $k+1$ points, with one special point distinguished,
into $\mu+2$ non-empty parts with one part distinguished, which does
not contain the special point. On the right-hand side, we first determine the size $\nu$
of the distinguished part, then select the points for this part avoiding the distinguished point,
and finally partition the remaining $k+1-\nu$-set into $\mu+1$ non-empty parts.
The number of possibilities in the latter case clearly matches the combinatorial description
of the left-hand side, and the result follows.
\end{proof}
\subsection{Proof of Proposition~\ref{prop:MainLemFuchs}}
\label{subsec:MainLemFuchs}
\begin{lemma}\label{Lem:FoLu}
Let $\pi\in S_n$ be a permutation consisting only of cycles of
length $q$. Then, for every irreducible character $\chi,$ we have
$|\chi(\pi)|\leq\sqrt{n}\,\big(\chi(1)\big)^{1/q}$.
\end{lemma}
This follows from \cite[Theorem 1.1]{FoLu} together with the
Murnaghan-Nakayama rule.
\begin{lemma}\label{Lem:FoLuWithFixed}
Let $q\geq 2$ be an integer, let $\pi\in S_n$ be a permutation
with $s$ cycles of lengths different from $q,$ and let $\lambda$
be a partition of $n$. Then
\begin{equation}\label{Eq:Hoeldereq}
|\chi_\lambda(\pi)|\leq\sqrt{n}\hspace{.4mm}\big(2\,\sq(\lambda)\big)^{s(q-1)/q}
\big(\chi_\lambda(1)\big)^{1/q}.
\end{equation}
\end{lemma}
\begin{proof}
For every partition $\mu$ denote by $N(\mu, \lambda)$ the number
of ways to obtain $\mu$ from $\lambda$ by stripping off $s$ rim
hooks of lengths according to the cycle structure of $\pi$, ending in an
element $\pi_0$ containing only cycles of length $q$. By the
Murnaghan-Nakayama rule, we have
\begin{equation}\label{Eq:DimEstH}
\chi_\lambda(1)\geq\sum_\mu N(\mu, \lambda)\chi_\mu(1).
\end{equation}
Neglecting the sign in the Murnaghan-Nakayama rule, and applying
Lemma \ref{Lem:FoLu}, H\"older's inequality, and the estimate
(\ref{Eq:DimEstH}), we get
\begin{eqnarray*}
|\chi_\lambda(\pi)| & \leq & \sum_\mu N(\mu, \lambda)
\hspace{.4mm}|\chi_\mu(\pi_0)|\\
 & \leq & \sqrt{n}\sum_\mu N(\mu, \lambda) \big(\chi_\mu(1)\big)^{1/q}\\
& \leq & \sqrt{n}\left(\sum_\mu N(\mu, \lambda)\right)^{(q-1)/q}
\left(\sum_\mu \chi_\mu(1)\hspace{.4mm}N(\mu, \lambda)
\right)^{1/q}\\
& \leq & \sqrt{n}\left(\sum_\mu N(\mu, \lambda)\right)^{(q-1)/q}
\big(\chi_\lambda(1)\big)^{1/q}.
\end{eqnarray*}
Arguing as in the proof of Lemma~\ref{Lem:squarevalue}, we see
that the $s$ cycles of length different from $q$ can be chosen in
at most $(2\,\mathrm{sq}(\lambda))^f$ ways, that is,
\[
\sum_\mu N(\mu, \lambda) \leq \big(2\,\mathrm{sq}(\lambda)\big)^s,
\]
and our claim follows.
\end{proof}
\begin{lemma}\label{Lem:MainLemCharEst}
Let $q\geq 2$ be an integer, let $\pi\in S_n$ be a permutation
with cycle lengths $\leq q,$ and let $\lambda\vdash n$ be a
partition of $n$ satisfying $\lambda_1\geq\|\lambda\|$. Set
$\Delta=n-\lambda_1,$ and let $\varepsilon>0$ be given. Then there
exists a constant $C=C(q, \varepsilon),$ such that, for
$\Delta\in[C, n/C]$ and $n$ sufficiently large,
\[
|\chi_\lambda(\pi)|\leq\chi_\lambda(1)^{\frac{1}{q}+\varepsilon}
\prod\limits_{t<q} E(n, \Delta, t, s_t(\pi)),
\]
where
\[
E(n, \Delta, t, s) = \begin{cases}\displaystyle
\Big(\frac{\Delta^{(q+2t-1)/2q}}{n^{t/q}}\Big)^s, &\mbox{\em if
$ts\leq
\Delta$}\\[3mm]
\displaystyle\Big(\frac{s}{\Delta^{\frac{2qt-2t-q+1}{2qt}}n^{t/q}}\Big)^\Delta,
& \mbox{\em if $ts> \Delta$.}
\end{cases}
\]
\end{lemma}
\begin{proof}
The proof is by induction on the minimum $m(\pi)$ of the cycle
lengths in $\pi$ from $q$ down to 1. For $m(\pi)=q$, our claim
follows immediately from Lemma \ref{Lem:FoLu}. Suppose that our
assumption holds for all $\pi$ with $m(\pi)\geq t+1$ and some $t\in[q-1]$.
Write $\mu=\lambda\setminus\lambda_1$. For
a partition $\nu\subseteq\lambda$ denote by $N(\nu, \lambda)$ the
number of ways to obtain $\nu$ from $\lambda$ by removing
$s_t(\pi)$ rim hooks of length $t$, let $\pi\in S_n$ be a
permutation with $m(\pi)=t$, and let $\pi_0\in S_{n-ts_t(\pi)}$ be
a permutation whose cycle structure is the same as that of $\pi$
with all $t$-cycles removed. Then we have
\[
|\chi_\lambda(\pi)| \leq \sum_\nu N(\nu, \lambda)\hspace{.4mm}
|\chi_\lambda(\pi_0)|.
\]
We claim that for a partition $\nu$ such that $N(\nu, \lambda)\neq 0$,
\[
|\chi_\nu(\pi_0)| \leq \big(\chi_\lambda(1)\big)^\varepsilon
\big(\chi_\nu(1)\big)^{\frac{1}{q}+\varepsilon} \prod_{t<\tau<q} E(n, \Delta,
\tau, s_\tau(\pi)).
\]
Indeed, if $\Delta'=n-t s_t(\pi)-\nu_1\geq\varepsilon \Delta$, then
this assumption holds (even without the factor
$(\chi_\lambda(1))^\varepsilon$) by the inductive hypothesis;
otherwise, using Lemma~\ref{Lem:chi1lower} and the fact that
$\Delta\leq n/C$, we get
\[
|\chi_\nu(\pi_0)| \leq \chi_\nu(1)\leq n^{\Delta'}\leq
n^{\varepsilon \Delta} \leq (\chi_\lambda(1))^{2\varepsilon}.
\]
Setting
\[
E:=\prod\limits_{t<\tau<q} E(n, \Delta, \tau, s_\tau(\pi)),
\]
this gives
\begin{eqnarray}
\label{Eq:asplit}
|\chi_\lambda(\pi)| & \leq & E \big(\chi_\lambda(1)\big)^\varepsilon \sum_\nu
N(\nu, \lambda) \big(\chi_\nu(1)\big)^{1/q}\nonumber\\[1.5mm]
 & = & E \big(\chi_\lambda(1)\big)^\varepsilon \sum_{a\leq s_t(\pi)}
\sum_{\underset{\nu_1=\lambda_1-(s_t(\pi)-a)t}{\nu}} N(\nu, \lambda)
\big(\chi_\nu(1)\big)^{1/q}.
\end{eqnarray}
Put $\kappa:= \nu\setminus\nu_1$. Given $a$, we bound
$\chi_\nu(1)$ as follows: we choose a set $I$ of $|\kappa| =
\Delta-at$ integers in $[n-ts_t(\pi)]$, and then count
the number of ways of removing all boxes of $\nu$ in such a way
that, in the $i$-th step, a box outside the first row is removed
if and only if $i\in I$. We claim that the latter number is
bounded by $\chi_\kappa(1)\leq \chi_\mu(1) N(\kappa, \mu)^{-1}$.
Indeed, refining the removal of a rim hook into a sequence of
removals of single boxes, we see that there are at least
$N(\kappa, \mu)$ ways to obtain $\kappa$ from $\mu$ by removing
single boxes, while there are $\chi_\kappa(1)$ ways to remove
$\kappa$ completely by deleting boxes. Hence, we have
$\chi_\mu(1)\geq N(\kappa, \mu)\chi_\kappa(1)$, from which the
last claim follows. Since $I$ can be chosen in
${n-ts_t(\pi)}\choose{\Delta-at}$ different ways, we
obtain that
\begin{equation}
\label{Eq:ChiNuEstLoc} \chi_\nu(1)\leq
{{n-ts_t(\pi)}\choose{\Delta-at}}\hspace{.4mm}\chi_\mu(1)N(\kappa,
\mu)^{-1}.
\end{equation}
Next, if $\nu_1=\lambda_1-(s_t(\pi)-a)t$, $N(\nu, \lambda)$ is the number of
ways to remove $s_t(\pi)-a$ rim hooks from the first row of $\lambda$, and
$a$ rim hooks from $\mu$. The position in the sequence of
steps where a rim hook is removed from the first row can be chosen
in $\binom{s_t(\pi)}{a}$ ways, thus
\begin{equation}
\label{Eq:ChiNuEstLoc2}
N(\nu, \lambda)=\binom{s_t(\pi)}{a}
N(\kappa, \mu).
\end{equation}
Inserting (\ref{Eq:ChiNuEstLoc}) and (\ref{Eq:ChiNuEstLoc2}) into
(\ref{Eq:asplit}), we get
\begin{equation*}
|\chi_\lambda(\pi)| \leq E \big(\chi_\lambda(1)\big)^\varepsilon \sum_{a\leq
s_t(\pi)} \sum_{\underset{\kappa\vdash
\Delta-at}{\kappa\subseteq\mu}} \binom{s_t(\pi)}
{a}\binom{n-ts_t(\pi)}{\Delta-at}^{1/q} N(\kappa,
\mu)^{1-1/q} \big(\chi_\mu(1)\big)^{1/q}.
\end{equation*}
Finally,
\[\sum_{\underset{\kappa\vdash \Delta-at}
{\kappa\subseteq\mu}}N(\kappa, \mu)\] is the number of ways to
remove $a$ rim hooks from $\mu$, which can be done in
at most $(2\sq(\mu))^a\leq
(2\sqrt{\Delta})^a$ ways, as we saw in the proof of
Lemma 2. Applying H\"older's inequality, and observing the fact
that the number of partitions of $\Delta$ is
$e^{c\sqrt{\Delta}}\leq (\chi_\lambda(1))^\varepsilon$, we obtain
\begin{equation*}\label{penultimate}
|\chi_\lambda(\pi)| \leq E \big(\chi_\lambda(1)\big)^\varepsilon
\!\!\sum_{a\leq s_t(\pi)}
\binom{s_t(\pi)}{a}\binom{n-ts_t(\pi)}{\Delta-at}^{1/q}
(2\sqrt{\Delta})^{(1-1/q)a}\big(\chi_\mu(1)\big)^{1/q}.
\end{equation*}
Since by Lemma~\ref{Lem:chi1lower} (ii) we have for $C>3/\varepsilon$
\[
\big(\chi_\lambda(1)\big)^{\frac{1}{q}+\varepsilon}\geq
\left(\frac{n^\Delta}{\Delta!}\right)^{1/q}\big(\chi_\mu(1)\big)^{1/q},
\]
we obtain
\begin{equation}\label{norepest}
\frac{|\chi_\lambda(\pi)|}{(\chi_\lambda(1))^{\frac{1}{q}+\varepsilon}}
\leq E \left(\frac{\Delta!}{n^\Delta}\right)^{1/q} \sum_{a\leq s_t(\pi)}
\binom{s_t(\pi)}{a} \binom{n-ts_t(\pi)}{\Delta-at}^{1/q}
(2\sqrt{\Delta})^{(1-1/q)a}.
\end{equation}
Since by assumption $\Delta\geq C$, the number $s_t(\pi)\leq n$ of summands is
of order at most $(\chi_\lambda(1))^\varepsilon$; in particular, we can
estimate the sum over $a$ by its largest term. We now distinguish two
cases, according to whether $\Delta\geq ts_t(\pi)$ or $\Delta<ts_t(\pi)$.\\[3mm]
Case 1: $\Delta\geq ts_t(\pi)$. We first note that terms of the order of
magnitude $e^{c\Delta}$ can be neglected on the right-hand side of
(\ref{norepest}), since they are bounded above by
$(\chi_\lambda(1))^\varepsilon$; in particular, $\binom{s_t(\pi)}{a}\leq 2^{\Delta/t}$ and
$2^a\leq 2^{\Delta/t}$ are absorbed into the term
$(\chi_\lambda(1))^\varepsilon$. We now split the summation over $a$ into
the ranges $a\leq\varepsilon s_t(\pi)$, $\varepsilon s_t(\pi)<a<s_t(\pi)-\varepsilon\Delta$ and
$a\geq s_t(\pi)-\varepsilon\Delta$. In the last case, we have
\begin{eqnarray*}
\left(\frac{\Delta!}{n^\Delta}\right)^{1/q} \binom{n-ts_t(\pi)}{\Delta-at}^{1/q}
\Delta^{\frac{q-1}{2q}a}
 & \leq & \left(\frac{\Delta!\hspace{.4mm}n^{\Delta-ts_t(\pi)+\varepsilon \Delta t}}
{n^\Delta(\Delta-ts_t(\pi)+\varepsilon \Delta t)!}\right)^{1/q} \Delta^{\frac{q-1}{2q}s_t(\pi)}\\[1.5mm]
 & \leq & n^{\varepsilon \Delta t/q}
\left(\frac{\Delta^{\frac{1}{2}+\frac{t}{q}-\frac{1}{2q}}}{n^{t/q}}\right)^{s_t(\pi)}\\[1.5mm]
 & \leq & \big(\chi_\lambda(1)\big)^\varepsilon
\left(\frac{\Delta^{\frac{1}{2}+\frac{t}{q}-\frac{1}{2q}}}{n^{t/q}}\right)^{s_t(\pi)}.
\end{eqnarray*}
Hence, every term in this range is of the desired magnitude, and
therefore this part of the sum is sufficiently small.\\[2mm]
Next, we turn our attention to terms with $a\leq \varepsilon s_t(\pi)$.  In this case, we obtain
\begin{eqnarray*}
\left(\frac{\Delta!}{n^\Delta}\right)^{1/q} \binom{n-ts_t(\pi)}{\Delta-at}^{1/q}
\Delta^{\frac{q-1}{2q}a}
 & \leq &\left(\frac{\Delta!}{n^\Delta}\right)^{1/q}\binom{n}{\Delta}^{1/q}
\Delta^{\varepsilon\frac{q-1}{2q}s_t(\pi)}\\[1.5mm]
 & \leq & \big(\chi_\lambda(1)\big)^\varepsilon.
\end{eqnarray*}
Finally, consider the range $\varepsilon s_t(\pi) < a < s_t(\pi)-\varepsilon\Delta$. It
suffices to consider the case where the terms in this range are not
dominated by terms in the other ranges, that is, we may assume that the
maximal term lies within this range. If we increase $a$ by 1, a single
summand in (\ref{norepest}) is changed by a factor
\[
F(a) := \frac{\big((n-ts_t(\pi)-\Delta+at)\cdots
(n-ts_t(\pi)-\Delta+at-t+1)\big)^{1/q}}
{\big((\Delta-at+1)\cdots(\Delta-at+t)\big)^{1/q} \Delta^{1/2-1/2q}}.
\]
$F(a)$ is monotonically increasing, hence there is a unique positive
solution $x_0$ of the equation $F(x)=1$, and the value $a_{\max}$ for which
the corresponding summand is maximal, differs from $x_0$ by
at most 1. Using the bounds for $a$, we find that
\[
F(a)\asymp\frac{n^{t/q}}{\Delta^{\frac{q+2t-1}{2q}}},\quad a\in(\varepsilon s_t(\pi),
s_t(\pi)-\varepsilon\Delta);
\]
hence, we can neglect this range, unless the expression on the
right-hand side, which does not depend on $a$, is $\asymp 1$;
that is, $n\asymp \Delta^{(q+2t-1)/2t}$. In the latter case, we have
\begin{eqnarray*}
\left(\frac{\Delta!}{n^\Delta}\right)^{1/q} \binom{n-ts_t(\pi)}{\Delta-at}^{1/q}
\Delta^{\frac{q-1}{2q}a}
 & \leq & \frac{\Delta^{at/q+\frac{q-1}{2q}a}}{n^{at/q}}\\[1mm]
 & \leq & \left(\frac{\Delta^{q+2t-1}}{n^{2t}}\right)^{\frac{a}{2q}}\\[1mm]
 & \leq & e^{c\Delta}\\[1mm]
 & \leq & \big(\chi_\lambda(1)\big)^\varepsilon.
\end{eqnarray*}
Case 2: $ts_t(\pi)> \Delta$. Note that we have $a\leq \Delta/t$. We begin
with $a$ in the range $a\geq (1-\varepsilon)\Delta/t$. Then
we have
\begin{eqnarray*}
\left(\frac{\Delta!}{n^\Delta}\right)^{1/q}{s_t(\pi)\choose a}
\binom{n-ts_t(\pi)}{\Delta-at}^{1/q} \Delta^{\frac{q-1}{2q}a}
&\leq&\left(\frac{\Delta!}{n^\Delta}\right)^{1/q} {s_t(\pi)\choose \Delta}
\binom{n-ts_t(\pi)}{\varepsilon \Delta}^{1/q} \Delta^{\frac{q-1}{2qt}\Delta}\\[1.5mm]
&\leq&\big(\chi_\lambda(1)\big)^\varepsilon
\left(\frac{s_t(\pi)}{\Delta^{\frac{2qt-2t-q+1}{2qt}}n^{1/q}}\right)^{\Delta},
\end{eqnarray*}
which is the desired result.\\[3mm]
If $a\leq \varepsilon \Delta/t$, then we have
\begin{eqnarray*}
\left(\frac{\Delta!}{n^\Delta}\right)^{1/q} {s_t(\pi)\choose a}
{{n-ts_t(\pi)}\choose{\Delta-at}}^{1/q} \Delta^{\frac{q-1}{2q}a}
&\leq&  \left(\frac{\Delta!}{n^\Delta}\right)^{1/q} \binom{s_t(\pi)}{\varepsilon \Delta/t}
\binom{n}{\Delta}^{1/q} \Delta^{\varepsilon\frac{q-1}{2q}\Delta}\\[1.5mm]
&\leq& \big(\chi_\lambda(1)\big)^\varepsilon
\left(\frac{s_t(\pi)}{\Delta}\right)^{\varepsilon \Delta},
\end{eqnarray*}
which is less than $(\chi_\lambda(1))^\varepsilon\hspace{.4mm}E(n, \Delta, t, s_t(\pi))$.\\[3mm]
Finally, if $a\in(\varepsilon\Delta/t, (1-\varepsilon) \Delta/t)$,
increasing $a$ by 1 changes a single summand by a factor
\begin{multline*}
F(a):=  \frac{(s_t(\pi)-a)\big((n-\Delta-t s_t(\pi)+at)\cdots(n-\Delta-t s_t(\pi)+at-t+1)\big)^{1/q}}
{(a+1)\Delta^{1/2-1/2q} \big((\Delta-at+1)\cdots(\Delta-at+t)\big)^{1/q}}\\[2mm]
 \asymp  \frac{n^{t/q}\Delta^{1/2-t/q+1/2q}}{s_t(\pi)};
\end{multline*}
and, as in the first case, the summands corresponding to these values of $a$ can be neglected,
unless the last expression is $\asymp 1$. If this is the case, we compute
a single summand to be
\begin{eqnarray*}
\left(\frac{\Delta!}{n^\Delta}\right)^{1/q} {s_t(\pi)\choose a}
{{n-ts_t(\pi)}\choose{\Delta-at}}^{1/q} \Delta^{\frac{q-1}{2q}a}
&\leq&\big(\chi_\lambda(1)\big)^\varepsilon\left(\frac{\Delta!}{n^\Delta}\right)^{1/q}
\left(\frac{s_t(\pi)}{\Delta^{1/2+1/2q}}\right)^a
\left(\frac{n^{1/q}}{\Delta^{1/q}}\right)^{\Delta-at}\\[1.5mm]
&\leq&\big(\chi_\lambda(1)\big)^\varepsilon\left(\frac{s_t(\pi)}
{\Delta^{1/2-t/q+1/2q}n^{t/q}}\right)^a\\[1.5mm]
&\leq&\big(\chi_\lambda(1)\big)^\varepsilon e^{c\Delta}.
\end{eqnarray*}
\end{proof}
We are now in a position to prove Proposition~\ref{prop:MainLemFuchs}.
Let $\lambda$ be a partition of $n$, $q\geq 2$ an integer, and $\pi\in S_n$
a permutation such that $\pi^q=1$. Then
$s_t(\pi)=0$, unless $t|q$. We prove (\ref{Eq:MainLemFuchs}), using
estimates in different ranges for cycle numbers of $\pi$ and $\Delta=n-\lambda_1$.\\[3mm]
If $\Delta=0$, the assertion is trivial, and if $1\leq \Delta\leq n^{\varepsilon}$, and
$s_t(\pi)\geq 2en^{t/q}$ for some $t$, then we
use the trivial estimate $\chi_\lambda(\pi)\leq\chi_\lambda(1)\leq n^\Delta$,
together with Lemma \ref{Lem:statUppbound} to see
that the contribution of such terms to the sum in question is
$<e^{-cn^{1/q}}$. If on the other hand $s_t(\pi)<2en^{t/q}$ for all $t$, we estimate
$|\chi_\lambda(\pi)|$ as follows. As in the proof of Lemma~\ref{Lem:LargeL1CharEst},
we choose $a_t$ cycles of length $t$ without
boxes from the first row, and obtain
\[
|\chi_\lambda(\pi)| \leq \Delta! \sum_{\sum_t ta_t\leq\Delta}
\prod\limits_t{2en^{t/q}\choose a_t}\leq \Delta!\hspace{.4mm}\Delta^{\tau(q)} (2en)^{\Delta/q} \leq (\Delta!)^2
\big(\chi_\lambda(1)\big)^{1/q},
\]
and this is of the desired order of magnitude, since
\[
\big(\chi_\lambda(1)\big)^{3\varepsilon}\geq \left(\frac{n}{\Delta}\right)^{3\varepsilon\Delta}\geq
n^{2\varepsilon\Delta} \geq \Delta^{2\Delta}.
\]
Next, we consider the case where $\Delta$ is in the range
$[n^\varepsilon,n/C]$ for some sufficiently large constant $C$.
Assume first,
that there is some $t$ such that $s_t(\pi)<2en^{t/q}$. Then we have
$E(n, \Delta, t, s_t(\pi))\leq (\chi_\lambda(1))^\varepsilon$. For, either
$\Delta<ts_t(\pi)$, which implies
\[
s_t(\pi) < 2en^{t/q} \leq \Delta^{1/2q} n^{t/q} \leq \Delta^{\frac{q-2t+1}{2q}} n^{t/q};
\]
or $ts_t(\pi)\leq \Delta\leq ts_t(\pi)^{1+1/2q}$, in which case
\[
\Delta^{\frac{q-2t+1}{2q}} \leq \Delta^{1-1/(2q)} \leq t\big(s_t(\pi)\big)^{1-1/(4q^2)} \leq n^{t/q};
\]
or, finally, $\Delta> ts_t(\pi)^{1+1/2q}$, which implies
\[
E(n, \Delta, t, s_t(\pi))\leq n^{s_t(\pi)}\leq 2^\Delta \leq \big(\chi_\lambda(1)\big)^\varepsilon.
\]
Hence, disregarding the factor corresponding to such a
value of $t$ does not change the estimate in Lemma \ref{Lem:MainLemCharEst}
significantly. From this observation and Lemmas~\ref{Lem:statUppbound} and
\ref{Lem:MainLemCharEst} we obtain for $\pi$ of order dividing $q$
\begin{equation}\label{Eq:MedCharEst}
\frac{N(q, n, s_1(\pi), \ldots s_T(\pi))}{N(q, n)}\hspace{.4mm}|\chi(\pi)| \leq
\big(\chi_\lambda(1)\big)^{\frac{1}{q}+\varepsilon} \underset{s_t(\pi)\geq 2en^{t/q}}
{\underset{t<q}{\prod\limits_{t|q}}}
\left(\frac{en^{t/q}}{ts_t(\pi)}\right)^{s_t(\pi)} E(n, \Delta, t, s_t(\pi)).
\end{equation}
If $t$ is such that $s_t(\pi)>\Delta$, then
\[
\left(\frac{en^{t/q}}{s_t(\pi)}\right)^{s_t(\pi)} E(n, \Delta, t, s_t(\pi)) =
\left(\frac{en^{t/q}}{s_t(\pi)}\right)^{s_t(\pi)-\Delta} \Delta^{-\frac{2qt-2t-q+1}{2qt}\Delta} < 1;
\]
if $\Delta/t\leq s_t(\pi)\leq\Delta$, we find that
\begin{eqnarray*}
\left(\frac{en^{t/q}}{s_t(\pi)}\right)^{s_t(\pi)} E(n, \Delta, t, s_t(\pi))
 & = & \left(\frac{en^{t/q}}{s_t(\pi)}\right)^{s_t(\pi)-\Delta} \Delta^{-\frac{2qt-2t-q+1}{2qt}\Delta}\\[1.5mm]
 & \leq & e^{c\Delta} \left(\frac{\Delta}{n^{t/q}}\right)^{(1-\frac{1}{t})\Delta}\Delta^{-\frac{2qt-2t-q+1}{2qt}\Delta}\\[1.5mm]
 & \leq & e^{c\Delta} \left(\frac{\Delta^{-\frac{1}{2t}+\frac{1}{q}-\frac{1}{2qt}}}{n^{(t-1)/q}}\right)^\Delta\\[1.5mm]
 & \leq & \big(\chi_\lambda(1)\big)^\varepsilon,
\end{eqnarray*}
since $t$ is bounded by $q/2$.
Finally, if $ts_t(\pi)\leq \Delta$, we have
\[
\left(\frac{en^{t/q}}{s_t(\pi)}\right)^{s_t(\pi)} E(n, \Delta, t, s_t(\pi)) =
\left(\frac{e\Delta^{\frac{q+2t-1}{2q}}}{s_t(\pi)}\right)^{s_t(\pi)} \leq
\left(\frac{e\Delta^{1-1/2q}}{s_t(\pi)}\right)^{s_t(\pi)},
\]
and either $s_t(\pi)<\Delta^{1-1/4q}$, in which case this factor is
$\leq(\chi_\lambda(1))^\varepsilon$, or $s_t(\pi)\geq \Delta^{1-1/4q}$, in which
case it is less than 1. Hence, in any case the right-hand side of
(\ref{Eq:MedCharEst}) is bounded by $(\chi_\lambda(1))^\varepsilon$. Summing over all
possible values for the $s_t(\pi)$ gives an additional factor $\leq n^{\tau(q)}$, which is
absorbed into $(\chi_\lambda(1))^\varepsilon$ as well. Hence, for these characters
(\ref{Eq:MainLemFuchs}) holds.\\[3mm]
Finally, we have to consider partitions with $\Delta>n/C$. By Lemma~\ref{Lem:chi1lower}, this
implies $\chi_\lambda(1)\geq e^{cn}$. By Lemma \ref{Lem:FoLuWithFixed}, we have
\[
|\chi_\lambda(\pi)|\leq (2\sqrt{n})^{\frac{q-1}{q}s}\big(\chi_\lambda(1)\big)^{\frac{1}{q}+\varepsilon},
\]
where $s=\sum_{t<q} s_t(\pi)$. Together with Lemma~\ref{Lem:statUppbound}, we deduce
\begin{eqnarray*}
\frac{|\chi_\lambda(\pi)|N(n, q, s_1(\pi), \ldots, s_T(\pi))}
{(\chi_\lambda(1))^{1/q+\varepsilon}N(n, q)} & \leq & (4n)^{\frac{q-1}{2q}s}
\underset{t<q}{\prod_{t|q}}\left(\frac{en^{t/q}}{s_t(\pi)}\right)^{s_t(\pi)}\\[1.5mm]
 & = & \underset{t<q}{\prod_{t|q}}\left(\frac{4en^{\frac{q+2t-1}{2q}}}{s_t(\pi)}\right)^{s_t(\pi)}.
\end{eqnarray*}
In the last product, all factors with $s_t(\pi)>4en^{\frac{q+2t-1}{2q}}$ are bounded by 1,
while for $t$ such that $s_t(\pi)\leq 4en^{\frac{q+2t-1}{2q}}$ the corresponding factor is
at most $n^{n^{1-\frac{1}{2q}}}<(\chi_\lambda(1))^\varepsilon$. Hence, also in this case,
the left-hand side of (\ref{Eq:MainLemFuchs}) has the desired order of magnitude, and the proof of
Proposition~\ref{prop:MainLemFuchs} is complete.
\section{The multiplicities of root number functions}
\label{Sec:Multiplicities}
In order to be able to estimate the subgroup growth of Fuchsian groups,
we also need to establish certain properties of the multiplicities of root number functions
for symmetric groups. These are summarized in our next two results.
For a positive integer $q$, define the $q$-th root number function
$r_q:S_n\rightarrow\mathbb{N}_0$ via
\[
r_q(\pi):=
\big|\big\{\sigma\in S_n:\sigma^q=\pi\big\}\big|,
\]
and, for each
irreducible character $\chi$ of $S_n$, let
\[
m_\chi^{(q)}:=\langle r_q, \chi\rangle
\]
be the multiplicity of $\chi$ in $r_q$. It is known that the functions $r_q$ are proper
characters, that is, the $m_\chi^{(q)}$ are non-negative integers; cf. \cite{Scharf}.
\begin{proposition}
\label{Prop:RootNumberEst}
Let $q\geq 2$ be an integer, $\varepsilon>0,$ $n\geq n_0(\varepsilon),$ and let
$\lambda\vdash n$ be a partition with corresponding character $\chi_\lambda$.
\begin{enumerate}\vspace{-3mm}
\item We have $m_{\chi_\lambda}^{(q)} < (\chi_\lambda(1))^{1-2/q+\varepsilon}$.\\[-3mm]
\item Given a partition $\mu\vdash \Delta,$ there exists some constant $C_\mu^q,$ depending only
on $\mu$ and $q,$ such that, for $\lambda\setminus\lambda_1=\mu$
and $n$ sufficiently large, we have $m_{\chi_\lambda}^{(q)}=C_\mu^q$. In particular,
we have
\begin{eqnarray*}
C_{(1)}^q & = & \tau(q)-1,\\
C_{(2)}^q & = & \frac{1}{2}\big(\sigma(q)+\tau(q)^2-3\tau(q)+\tau_\mathrm{odd}(q)\big),\\
C_{(1,1)}^q & = & \frac{1}{2}\big(\sigma(q)+\tau(q)^2-3\tau(q)-\tau_\mathrm{odd}(q)\big)+1,
\end{eqnarray*}
where $\tau_\mathrm{odd}(q)$ is the number of odd divisors of $q$.
\item For a partition $\mu\vdash\Delta$, $q$ odd, and sufficiently large $n$, we have
$m_{\chi_{\lambda'}}^{(q)}=0$.
\end{enumerate}
\end{proposition}
\begin{proposition}
\label{Prop:RootNumberAsymp}
Let $\varepsilon>0$ be given, let $q\geq 2$ be an integer, $\Delta\geq\Delta_0(q, \varepsilon)$
and $n\geq n_0(q, \Delta, \varepsilon)$. Then, for partitions
$\mu\vdash \Delta$ and $\lambda\vdash n$ with
$\lambda\setminus\lambda_1=\mu,$ we have that
\[
\left|\frac{m_{\chi_\lambda}^{q} \Delta!}{\chi_\mu(1) Q_\Delta(q)} -
1\right| < \varepsilon,
\]
where $Q_\Delta(q) = \sum_{i=1}^\Delta S(\Delta, \Delta-i)\hspace{.4mm}q^i$.
\end{proposition}
This section is devoted to the proofs of these results.
\subsection{Proof of Proposition~\ref{Prop:RootNumberEst}}
\label{subsec:prop2}
We begin by translating the problem
of bounding $m_\chi^{(q)}$ from an algebraic into a combinatorial question.
\begin{lemma}
For $q\geq2$ and $\chi\in\mbox{\em Irr}(S_n),$ we have
\[
m_\chi^{(q)}\leq \sum_{{\bf c}^q=1}\left\langle {\bf 1},
\big|\chi\!\!\downarrow_{C_{S_n}({\bf c})}\big|\right\rangle_{C_{S_n}({\bf c})},
\]
where the summation extends over all conjugacy classes ${\bf c}$ in $S_n$, whose
orders divide $q$, and $C_{S_n}({\bf c})$ denotes the
centralizer of some element $\pi\in {\bf c}$.
\end{lemma}
\begin{proof}
It is shown in \cite{Scharf} that, for every class ${\bf c}$, there exists a linear
character $\phi_{\bf c}$ of $C_{S_n}({\bf c})$ such that $r_q=\sum\limits_{{\bf c}^q=1}
\phi_{\bf c}\!\!\uparrow^{S_n}$. By Frobenius reciprocity this implies
\begin{eqnarray*}
m_\chi^{(q)} & = & \sum_{{\bf c}^q=1}
\left\langle\phi_{\bf c}\!\!\uparrow^{S_n}, \chi\right\rangle_{S_n}\\[1.5mm]
 & = & \sum_{{\bf c}^q=1}\left\langle \phi_{\bf c},
 \chi\!\!\downarrow_{C_{S_n}({\bf c})}\right\rangle_{C_{S_n}({\bf c})}\\[1.5mm]
 & \leq & \sum_{{\bf c}^q=1}\left\langle {\bf 1},
\big|\chi\!\!\downarrow_{C_{S_n}({\bf c})}\big|\right\rangle_{C_{S_n}({\bf c})}.
\end{eqnarray*}
\end{proof}
\begin{proofof}{Proposition~\em\ref{Prop:RootNumberEst}}
(i) Let $\varepsilon>0$ be given.
We first consider characters $\chi_\lambda$ with $\chi_\lambda(1)>n^{n^{1-\varepsilon/4}}$,
starting from the formula
\[
m_{\chi_\lambda}^{(q)}\leq\sum_{\mathbf{c}^q=1}\frac{1}{|C_{S_n}({\bf c})|}\sum_{\pi\in
C_{S_n}({\bf c})}|\chi_\lambda(\pi)|.
\]
The number of summands of the outer sum is $\leq n^{\tau(q)}$ and therefore
negligible. In the inner sum, we bound $\chi_\lambda(\pi)$ by either
using Lemma \ref{Lem:squarevalue} or via the trivial bound
$\chi_\lambda(1)$, depending on the number of cycles of $\pi$.
Estimating the number of elements $\pi\in C_{S_n}(\mathbf{c})$ with $k$ cycles using
Lemma~\ref{Lem:CycleBound} (ii) for $k\geq (\log n)^3$, and trivially otherwise, we obtain
\[
m_{\chi_\lambda}^{(q)} \leq \big(\chi_\lambda(1)\big)^\varepsilon (2\sqrt{n})^{(\log n)^3}
+ \big(\chi_\lambda(1)\big)^\varepsilon \sum_{k\geq (\log n)^3}
\min\big((2\sqrt{n})^k, \chi_\lambda(1)\big)\frac{(3q\log
n)^{k/q}}{\lfloor k/q\rfloor!}.
\]
The first summand is negligible. The greatest term of the sum is
coming from one of $k=\lfloor\frac{2\log\chi_\lambda(1)}{\log
n}\rfloor$ and $k=\lceil\frac{2\log\chi_\lambda(1)}{\log
n}\rceil$, and these terms differ by a factor $n$ at most; hence,
using Stirling's formula, we obtain
\[
m_{\chi_\lambda}^{(q)} \leq \big(\chi_\lambda(1)\big)^{1+\varepsilon}\,
e^{-\frac{2}{q}\log\chi_\lambda(1)\frac{\log\log\chi_\lambda(1)}{\log
n}} \leq \big(\chi_\lambda(1)\big)^{1-\frac{2}{q}+2\varepsilon},
\]
and our claim is proven in this case. In particular, setting
$\Delta=n-\lambda_1$, our first claim holds true for all characters
$\chi_\lambda$ belonging to a
partition $\lambda$ with $\Delta>n^{1-\varepsilon/2}$.\\[3mm]
Next, we consider the range $(\log n)^4\leq \Delta\leq
n^{1-\varepsilon/2}$. In this range we estimate $\chi_\lambda(\pi)$ by means of
Lemma~\ref{Lem:LargeL1CharEst} and Equation (\ref{Eq:L1CharEst2}), and we
bound the number of centralizer elements with $k$
cycles again via Lemma~\ref{Lem:CycleBound} (ii), or trivially. In this way, we obtain
\begin{eqnarray*}
m_{\chi_\lambda}^{(q)} & \leq & \big(\chi_\lambda(1)\big)^\varepsilon +
\big(\chi_\lambda(1)\big)^\varepsilon \sum_{k\geq (\log n)^3}
\min\bigg(\max_{\nu\leq \Delta}(2\sqrt{\Delta})^\nu\binom{k}{\nu},
\chi_\lambda(1)\bigg)\frac{(3q\log n)^{k/q}}{\lfloor k/q\rfloor!}\\[1mm]
 & \leq & \big(\chi_\lambda(1)\big)^{2\varepsilon} \sum_{k\geq (\log n)^3}
\min\Big((4\sqrt{\Delta})^k,  \chi_\lambda(1)\Big)\frac{(3q\log
n)^{k/q}}{\lfloor k/q\rfloor!}\\[1mm]
 & \leq & \big(\chi_\lambda(1)\big)^{1+3\varepsilon}\,
e^{-\frac{2}{q}\log\chi_\lambda(1)\frac{\log\log\chi_\lambda(1)}
{\log \Delta}}.
\end{eqnarray*}
On the other hand, from Lemma~\ref{Lem:chi1lower} (ii) we see that
$\frac{\log\log\chi_\lambda(1)}{\log \Delta}\geq 1$ for all
$\Delta<n/5$,
hence, our claim holds in this case as well.\\[3mm]
Finally, for $\Delta<(\log n)^4$, we see from Lemma~\ref{Lem:CycleBound} (ii)
that we may neglect
the contribution of permutations with at least $(\log n)^3$
cycles. For the remaining ones we have $\chi_\lambda(\pi)\leq (\log
n)^{3\Delta}$ by Lemma~\ref{Lem:LargeL1CharEst}, Equation (\ref{Eq:L1CharEst1}),
which is less then $(\chi_\lambda(1))^\varepsilon$.
This completes the proof of Proposition~\ref{Prop:RootNumberEst} (i).\\[3mm]
(ii) We describe the computation of  $m_{\chi_\lambda}^{(q)}$ for bounded $\Delta$. We have
\[
m_{\chi_\lambda}^{(q)} = \langle \chi_\lambda, r_q\rangle = \frac{1}{n!}\sum_{\pi\in S_n} \chi_\lambda(\pi^q).
\]
Let $\mu\vdash \Delta$ be a partition. By Lemma~\ref{Lem:CharPolynom}, there exists a polynomial
$P_\mu(x_1, \ldots, x_\Delta)$, such that, for every $n$ and $\lambda\vdash n$ with $\lambda\setminus\lambda_1=\mu$,
we have $\chi_\lambda=P_\mu(s_1, \ldots, s_\Delta)$. Moreover, if $x_d$ has weight $d$, then $P_\mu$ has
weight at most $|\mu|$. We have
\begin{equation}
\label{Eq:PowerCycle}
s_d(\pi^q) = \underset{\kappa/(\kappa, q) = d}{\sum_\kappa}  (\kappa, q) s_\kappa(\pi),
\end{equation}
since the $q$-th power of a cycle of length $\kappa$ consists of $(\kappa, q)$ cycles, each of length
$\frac{\kappa}{(\kappa, q)}$. Thus, we have to compute
\[
\frac{1}{n!}\sum_{\pi\in S_n} Q_{\mu, q}\big(s_1(\pi), s_2(\pi), \ldots, s_{\Delta q}(\pi)\big)
\]
for a certain polynomial $Q_{\mu, q}$ of weight at most $\Delta q$.

By Lemma~\ref{Lem:Poisson}, the last expression converges to some constant.
Hence, $m_{\chi_\lambda}^{(q)}$ converges to some real number, but as $m_{\chi_\lambda}^{(q)}$
is integral, it has to become constant for $n$ sufficiently large.

We now consider the cases $\mu=(1), (2)$ and $(1, 1)$. If $\lambda=(n-1, 1)$, we have
\[
m_{\chi_\lambda}^{(q)} = \frac{1}{n!}\sum_{\pi\in S_n} \chi_\lambda(\pi^q) = \frac{1}{n!}\sum_{\pi\in S_n}
s_1(\pi^q) - 1 = \frac{1}{n!}\sum_{\pi\in S_n}\sum_{d|q} ds_d(\pi) -1.
\]
By Lemma~\ref{Lem:Poisson}, the expected number of $d$-cycles is $\frac{1}{d}$; hence,
$m_{\chi_\lambda}^{(q)}=\tau(q)-1$. Next, for $\lambda=(n-2, 2)$, Lemma~\ref{Lem:CharPolynom} yields
\[
\chi_\lambda = \frac{s_1^2}{2} - \frac{3s_1}{2} + s_2.
\]
Inserting (\ref{Eq:PowerCycle}), we deduce
\begin{eqnarray*}
m_{\chi_\lambda}^{(q)} & = & \frac{1}{n!}\sum_{\pi\in S_n} \Bigg(\frac{1}{2}\Big(\sum_{d|q} ds_d(\pi)\Big)^2
-\frac{3}{2}\sum_{d|q} ds_d(\pi) + \underset{(2d, q)=d}{\sum_{d|q}} ds_d(\pi)\Bigg)\\[1.5mm]
 & = & \frac{1}{2}\sum_{d_1, d_2|q} \frac{1}{n!}\sum_{\pi\in S_n} \big(d_1d_2s_{d_1}(\pi)s_{d_2}(\pi)\big)
 -\frac{3}{2}\sum_{d|q}\frac{1}{n!}\sum_{\pi\in S_n} ds_s(\pi)\\[1mm]
 &&\hspace{6.1cm}+ \underset{(2d, q)=d}{\sum_{d|q}} \frac{1}{n!}\sum_{\pi\in S_n} ds_d(\pi).
\end{eqnarray*}
If $d_1\neq d_2$ then, by Lemma~\ref{Lem:Poisson}, $s_{d_1}(\pi)$ and $s_{d_2}(\pi)$ are asymptotically independent for
$\pi\in S_n$ chosen at random, and have mean values $\frac{1}{d_1}$ and $\frac{1}{d_2}$, respectively, while
$(s_d(\pi))^2$ has mean value $\frac{1}{d} + \frac{1}{d^2}$. Hence, the first sum is
asymptotically equal to $\sigma(q)+\tau(q)^2$. The second sum equals
$\tau(q)$, whereas the last sum yields $|\{d|q: (2d, q)=d\}|$, which equals $\tau_\mathrm{odd}(q)$.
We deduce that, as $n\rightarrow\infty$,
\[
m_{\chi_\lambda}^{(q)}\rightarrow \frac{1}{2}\big(\sigma(q)+\tau(q)^2-3\tau(q)+\tau_\mathrm{odd}(q)\big),
\]
which implies our claim, since $m_{\chi_\lambda}^{(q)}$ is always integral. A similar computation leads to
the value of $m_{\chi_\lambda}^{(q)}$ for $\lambda=(n-2, 1, 1)$.\\[3mm]
(iii) For $q$ odd, we have
\[
m_{\chi_{\lambda'}}^{(q)} = \sum_{\pi\in S_n}\chi_{\lambda'}(\pi^q)=
\sum_{\pi\in S_n}\chi_\lambda(\pi^q)\hspace{.4mm}\epsilon(\pi^q)=
\sum_{\pi\in S_n}\chi_\lambda(\pi^q)\hspace{.4mm}\epsilon(\pi),
\]
where $\epsilon$ is the sign character. As in the proof of part (ii), we can write
$m_{\chi_{\lambda'}}^q$ as a linear combination of sums of the form
\[
\frac{1}{n!}\sum_{\pi\in S_n} \epsilon(\pi)\prod_{j=1}^{lq}\big(s_j(\pi)\big)^{e_j}.
\]
Observe first that the contribution coming from permutations $\pi$ with $\sum_{j=1}^{ql} js_j(\pi)>n/2$
is $o(1)$ by Lemma~\ref{Lem:Poisson}. The sum over permutations $\pi$ with $\sum_{j=1}^{ql} js_j(\pi)\leq n/2$
vanishes, since it can be expressed as a linear combination of sums of the form
\[
\underset{s_i(\pi)\geq s_i\,(i\leq lq)}{\sum_{\pi\in S_n}}\epsilon(\pi) =
\pm\sum_{\pi\in S_{n-\mathcal{S}}}\epsilon(\pi) = 0,
\]
where $\mathcal{S}=\sum_{j=1}^{lq} s_j$. Hence, for $q$
odd and $n$ large, $m_{\chi_{\lambda'}}^{(q)}<1$ and therefore vanishes.
\end{proofof}
\subsection{Proof of Proposition~\ref{Prop:RootNumberAsymp}}
\label{subsec:Prop3}
We compute the scalar product $\langle r_q,\chi_\lambda\rangle$,
and evaluate $\chi_\lambda$ using Lemma~\ref{Lem:CharPolynom}, to obtain
\begin{equation}\label{Eq:RootNumberCharPol}
m_{\chi_\lambda}^{(q)} = \frac{\chi_\mu(1)}{n!}\sum_{\pi\in
S_n}{{s_1(\pi^q)}\choose{\Delta}} +
\mathcal{O}\left(\frac{\Delta\chi_\mu(1)}{n!}\sum_{i\leq \Delta-1}\sum_{\pi\in
S_n}{{s_1(\pi^q)+\ldots+s_\Delta(\pi^q)}\choose{i}}\right).
\end{equation}
By Lemma~\ref{Lem:Poisson}, we have as $n\rightarrow\infty$
\[
\frac{1}{n!}\big|\big\{\pi\in S_n: s_l(\pi) = a\big\}\big| \sim
\frac{e^{-1/l}}{l^a a!},
\]
and the
events ``$s_l(\pi)=a$'' and ``$s_{l'}(\pi)=b$'' are asymptotically
independent.
From this together with Equation (\ref{Eq:PowerCycle}) we deduce
\begin{eqnarray*}
\frac{1}{n!} \sum_{\pi\in S_n} \big(s_1(\pi^q)\big)^\Delta & = & \frac{1}{n!}
\sum_{\pi\in S_n} \left(\sum_{t|q} ts_t(\pi)\right)^\Delta\\
 & \sim & \sum_{\underset{\sum_t d_t=\Delta}{d_t, t|q}} \prod_{t|q}
\left(\frac{1}{n!}\sum_{\pi\in S_n} t s_t(\pi)\right)^{d_t}\\
 & \sim & \sum_{\underset{\sum_t d_t=\Delta}{d_t, t|q}} \prod_{t|q}
\sum_{a=1}^\infty \frac{e^{-1/t}(ta)^{d_t}}{t^a a!}\\
 & = & \sum_{\underset{\sum_t d_t=\Delta}{d_t, t|q}} \prod_{t|q}
t^{d_t} \sum_{\nu=1}^{d_t} S(d_t,\nu)\hspace{.4mm}t^{-\nu}.
\end{eqnarray*}
The last quantity can be written as
\begin{equation}
\label{Eq:StirlingSum}
\sum_{\underset{\sum_t d_t=\Delta}{d_t, t|q}} \prod_{t|q} Q_{d_t}(t),
\end{equation}
where
\[
Q_n(t) = \sum_{\nu=0}^{n-1} S(n, n-\nu)\hspace{.4mm}t^\nu.
\]
For fixed
$n$, the sequence $S(n, m)$ is unimodal in $m$; define $m_0$
to be the least $m$ with $S(n, m_0)=\max_m S(n, m)$. Kanold
\cite{Kanold} has shown that
\[
m_0\sim \frac{n}{\log n};
\]
moreover, it can be deduced from his proof that
\begin{equation}
\label{Eq:KanoldEst}
\sum_{m=(1-\varepsilon)m_0}^{(1+\varepsilon) m_0} S(n, m) \sim
\sum_{m=1}^n S(n, m),\quad n\rightarrow\infty.
\end{equation}
This estimate implies $Q_{d_t}(t) \leq
(t/q)^{d_t(1-\varepsilon)} Q_{d_t}(q)$ for all $t|q$ and $d_t\geq
d_0(\varepsilon)$.\\[3mm]
Next, we establish the inequality $Q_n(x)Q_{n'}(x)\leq
Q_{n+n'}(x)$ for all real positive $x$ and $n, m\geq 1$. More precisely, we
show that each single coefficient of $Q_n(x)Q_{n'}(x)$ is less than or equal to
the corresponding coefficient of $Q_{n+n'}(x)$, which implies our
claim. Computing the $m$-th coefficient explicitly, we have to show that
\[
\sum_{i+j=m} S(n, n'-i) S(n', n'-j) \leq S(n+n', n+n'-m),
\]
which is true since the right-hand side is the number of
partitions of a set with $n+n'$ elements into $n+n'-m$ parts,
while the left-hand side is the number of these partitions that
respect some fixed partition of the large set into two sets with
$n$ and $n'$ elements. Together with (\ref{Eq:KanoldEst}), we deduce that, for $d$ sufficiently large, and
a fixed tuple $d_1, \ldots, d_q$ with
$d_1+\dots+d_q=\Delta$,
\[
\prod_{t|q} Q_{d_t}(t) \leq \prod_{t|q} (t/q)^{2d_t/3} Q_{d_t}(q)
\leq (2/3)^{\Delta-d_q} \prod_{t|q} Q_{d_t}(q) \leq (2/3)^{\Delta-d_q}
Q_\Delta(q).
\]
We now split the sum (\ref{Eq:StirlingSum}) into three ranges, according to
whether $d_q=\Delta$, $\Delta-C\tau(q)\leq d_q\leq \Delta-1$, or $d_q<\Delta-C\tau(q)$,
and we want to show that the sum over the latter two ranges is negligible compared to the first term.
From the last estimate we find that, for $\varepsilon>0$ and
$\Delta>\Delta_0(q, \varepsilon)$, there is some $C=C(\varepsilon)$, such
that\footnote{Complex integration shows that we may take $C(\varepsilon)=10+7\log\varepsilon^{-1}$.}
\[
\sum_{\underset{d_q\leq \Delta-C\tau(q)}{\underset{\sum_t d_t=\Delta}{d_t,
t|q}}} \prod_{t|q} Q_{d_t}(t) \leq Q_\Delta(q)\sum_{\nu\geq C\tau(q)}
(2/3)^\nu {{\nu+\tau(q)-1}\choose{\tau(q)-1}}\leq\varepsilon Q_\Delta(q).
\]
In the range $\Delta-C\tau(q) \leq d_q\leq \Delta-1$, the number of summands is bounded,
and their sum can be estimated by
$CQ_{\Delta-1}(q)$, and we obtain
\[
\sum_{\underset{d_q\neq \Delta}{\underset{\sum_t d_t=\Delta}{d_t, t|q}}}
\prod_{t|q} Q_{d_t}(t) \leq \varepsilon Q_\Delta(q) + C Q_{\Delta-1}(q).
\]
Now we use the inequalities (cf. \cite[Satz 1]{Kanold})
\[
\frac{(m+1)^n}{m!}\left(1-\frac{m}{(1-1/m)^n}\right) \leq S(n, m)
\leq \frac{(m+1)^n}{m!},
\]
to see that $S(n, n-\mu)\hspace{.4mm}q^\mu$
in monotonically decreasing with $\mu$ for $n$ sufficiently large
and $\mu>n-\frac{n}{2\log n}$; in particular, we have
\[
\sum_{\mu>n-n/(2\log n)} S(n, n-\mu)\hspace{.4mm}q^\mu < \frac{1}{n} Q_n(q).
\]
Since $S(\Delta+1, \Delta-\mu)\geq (\Delta-\mu) S(\Delta, \Delta-\mu)$, the latter inequality implies
\begin{eqnarray*}
Q_\Delta(q) & \geq & \sum_{\mu\leq \Delta-\Delta/(2\log \Delta)} S(\Delta, \Delta-\mu) \hspace{.4mm}q^\mu\\
 & \geq & \frac{\Delta}{2\log \Delta} \sum_{\mu\leq \Delta-\Delta/(2\log \Delta)} S(\Delta-1, \Delta-\mu)
\hspace{.4mm}q^\mu\\
& \geq & \frac{\Delta}{2\log \Delta}\left(1-\frac{1}{\Delta-1}\right)Q_{\Delta-1}(q),
\end{eqnarray*}
and we deduce that
\[
Q_\Delta(q) \leq \frac{1}{n!}\sum_{\pi\in S_n}\big(s_1(\pi^q)\big)^\Delta \leq
(1+\varepsilon) Q_\Delta(q),
\]
provided that $\Delta\geq\Delta_0(q, \varepsilon)$ and $n\geq n_0(\Delta, q, \varepsilon)$.
Estimating the error term in (\ref{Eq:RootNumberCharPol})
in a way similar to our treatment of the main term, and using the fact that
\[
{{s_1(\pi^q)}\choose \Delta} = \frac{1}{\Delta!}\big(s_1(\pi^q)\big)^\Delta +
\mathcal{O}\big(s_1(\pi^q)^{\Delta-1}\big),
\]
Proposition~\ref{Prop:RootNumberAsymp} follows.
\section{Subgroup growth of Fuchsian groups}
\label{Sec:Growth}
\subsection{The generic case}
\label{Subsec:FuchsGeneric}
Let $r, s, t\geq 0$ be integers,
$a_1, \ldots, a_r\geq 2$ in $\mathbb{N}\cup\{\infty\}$, and let $e_1,
\ldots, e_s\geq 2$ be integral. Define the group $\Gamma = \Gamma(t;
a_1, \ldots, a_r; e_1, \ldots, e_s)$ associated with these data
by\footnote{$\Gamma$ is a Fuchsian group in the most general sense
met in the literature; cf., for instance, \cite[Prop.~5.3]{LynSch}
or \cite[Sec.~II.7]{Magnus}.}
\begin{multline}
\label{Eq:GammaDef}
\Gamma = \Big\langle x_1, \ldots, x_r, y_1, \ldots, y_s, u_1, v_1, \ldots, u_t, v_t\,\big|\,\\
x_1^{a_1}=x_2^{a_2}=\cdots=x_r^{a_r} =  x_1 x_2 \cdots x_r y_1^{e_1}y_2^{e_2}\cdots y_s^{e_s} [u_1, v_1] [u_2,v_2]\cdots
[u_t, v_t] = 1\Big\rangle.
\end{multline}
Define
\begin{eqnarray*}
\mu(\Gamma) & = & \sum_{i=1}^r\left(1-\frac{1}{a_i}\right)+s+2(t-1),\\
\alpha(\Gamma) & = & \mu(\Gamma)-\sum_{j=1}^s\left(1-\frac{2}{e_j}\right),\\
m_\Gamma & = & [a_1, \ldots, a_r].
\end{eqnarray*}
The main result of this section provides an asymptotic expansion for the
number of index $n$ subgroups of these groups.
\begin{theorem}\label{Thm:FuchsAsymp}
Let $\Gamma$ be given as in {\em (\ref{Eq:GammaDef}),} and suppose that $\alpha(\Gamma)>0$. Then
the number $s_n(\Gamma)$ of index $n$ subgroups in $\Gamma$ satisfies an asymptotic expansion
\[
s_n(\Gamma) \approx \delta L_\Gamma (n!)^{\mu(\Gamma)} \Phi_\Gamma(n)
\left\{1+\sum_{\nu\geq 1} a_\nu(\Gamma) n^{-\nu/m_\Gamma}\right\},
\qquad (n\rightarrow\infty).
\]
Here,
\begin{eqnarray*}
\delta  &=& \begin{cases} 2, &\forall i: a_i \text{\em\ finite and odd}, \forall j: e_j\text{\em\ even}\\[1mm]
1, & \text{\em otherwise,}
\end{cases}\\
L_\Gamma & = & (2\pi)^{-1/2-\sum_i(1-1/a_i)} (a_1 a_2\cdots a_r)^{-1/2} \exp\Bigg(
-\sum_{\underset{2|a_i}{i}}\frac{1}{2a_i}\Bigg),\\[1mm]
\Phi_\Gamma(n) & = & n^{3/2-\sum_i(1-1/a_i)} \exp\Bigg(\sum_{i=1}^r\sum_{\underset{t<a_i}{t|a_i}}\frac{n^{t/a_i}}{t}\Bigg),
\end{eqnarray*}
and the $a_\nu(\Gamma)$ are explicitly computable constants depending only on $\Gamma$.
\end{theorem}
\begin{corollary}
Let $\Gamma$ be as in {\rm (\ref{Eq:GammaDef})}, and suppose that $\alpha(\Gamma)>0$.
Then we have
\[
\frac{s_{n+1}(\Gamma)}{s_n(\Gamma)} \sim (n+1)^{\mu(\Gamma)},
\quad (n\rightarrow\infty);
\]
in particular, $s_n(\Gamma)$ is strictly increasing for
sufficiently large $n$.
\end{corollary}
\begin{proofof}{Theorem~{\rm \ref{Thm:FuchsAsymp}}}
The proof proceeds in three steps: first we express $h_n(\Gamma):= |\Hom(\Gamma, S_n)|/n!$
in character theoretic terms; next, we use results from Sections
\ref{Sec:MainLemFuchs} and \ref{Sec:Multiplicities} to obtain an asymptotic
estimate for $h_n(\Gamma)$. The assertions of the theorem are then deduced by means of an asymptotic
method for divergent power series due to Bender \cite{Bender}.\\[3mm]
Let $w$ be a segment of
\[
R=x_1 x_2\cdots x_r y_1^{e_1} y_2^{e_2}\cdots y_s^{e_s} [u_1, v_1] [u_2,v_2]\cdots[u_t, v_t],
\]
where powers $y_j^{e_j}$ and commutators are treated as single letters. For such a segment $w$
and $\pi\in S_n$, define $N_w(\pi)$ to be the number of solutions of the equation $w=\pi$, subject to
the conditions $x_i^{a_i}=1$ for those $i$ for which $x_i$ occurs in $w$. Observing that $N_w(\pi)$
is a class function, define coefficients $\beta_\chi^{(w)}$ via the expansion
\begin{equation}
\label{Eq:BetaDef}
N_w(\pi) = \frac{\prod_i|\Hom(C_{a_i}, S_n)|}{n!}\sum_{\chi\in \text{Irr}(S_n)} \beta_\chi^{(w)}
\chi(\pi),
\end{equation}
where the product is extended over those $i\in[r]$, such that $x_i$
occurs in $w$. Moreover, for a positive integer $q$ and an
irreducible character $\chi$, we set
\[
\alpha_\chi^{(q)} := \frac{1}{|\Hom(C_q, S_n)|}\sum_{\pi^q=1}\chi(\pi).
\]
We claim that
\begin{equation}\label{Eq:betaclaim}
\beta_\chi^{(w)} = (n!)^{|\{j:y_j\in w\}|+2|\{k:u_k\in w\}|}\;
\frac{\prod\limits_{x_i\in w} \alpha_\chi^{(a_i)}
\prod\limits_{y_j\in w}m_\chi^{(e_j)}}{(\chi(1))^{l(w)-1}},
\end{equation}
where $m_\chi^{(q)}=\langle r_q, \chi\rangle$, and $l(w)$ is the
number of different letters occurring in $w$. We prove
(\ref{Eq:betaclaim}) by induction on $l(w)$. If $w$ has length 1,
(\ref{Eq:betaclaim}) holds by the definition of
$\alpha_\chi^{(a_1)}$, the definition of $m_\chi^{(e_1)}$, or
\cite[Lemma 1]{MuPuSurf}, according to which letter occurs in $w$.
Similarly, the induction step breaks into three cases with regard
to the new letter; we only have to treat the case where the new
letter is $x_i$ for some  $2\leq i\leq r$, the other cases follow directly from
Lemma~\ref{Lem:FourierGen} (ii). We have
\begin{eqnarray*}
N_{wx_i}(\pi) & = & \sum_{\sigma\in S_n} N_w(\sigma)\cdot |\Hom(C_{a_i}, S_n)|\sum_\chi\alpha_\chi^{(a_i)}
\chi(\sigma^{-1}\pi) \\
 & = & \frac{\prod_{\nu=1}^i |\Hom(C_{a_\nu}, S_n)|}{n!} \sum_\sigma\sum_\chi
\alpha_\chi^{(a_i)}\beta_\chi^{(w)} \chi(\sigma)\chi(\sigma^{-1}\pi)\\
 & = & \frac{\prod_{\nu=1}^i |\Hom(C_{a_\nu}, S_n)|}{n!} \sum_\chi
\alpha_\chi^{(a_i)}\beta_\chi^{(w)} \sum_{\mathbf{c}_1, \mathbf{c}_2}
\frac{|\mathbf{c}_1|\cdot|\mathbf{c}_2|}{n!}\cdot
\frac{\chi(\mathbf{c}_1)\chi(\mathbf{c}_2)\chi(\pi)}{\chi(1)}\\
 & = & \frac{\prod_{\nu=1}^i |\Hom(C_{a_\nu}, S_n)|}{n!} \sum_\chi
\frac{\alpha_\chi^{(a_i)}\beta_\chi^{(w)}}{\chi(1)}\chi(\pi).
\end{eqnarray*}
Here, we have used orthogonality, and the fact\footnote{Cf., for instance, \cite[Prop.
9.33]{CurRei} or \cite[Theorem 6.3.1]{Kerber}.} that the number of
solutions in a finite group $G$ of the equation $x_1\cdot x_2=g$ with $x_i\in \mathbf{c}_i$ is
given by
\[
\frac{|\mathbf{c}_1|\cdot |\mathbf{c}_2|}{n!}\sum_\chi\frac{\chi(\mathbf{c}_1)\chi(\mathbf{c}_2)\chi(g^{-1})}{\chi(1)}.
\]
From this, (\ref{Eq:betaclaim}) follows.\\[3mm]
Now take $w$ as the full relation $R$, and set $\pi=1$ in the definition of $\beta_\chi^{(w)}$ to
obtain
\begin{equation}\label{hncharacter}
h_n(\Gamma) = (n!)^{s+2t-2}\prod_{i=1}^r |\Hom(C_{a_i}, S_n)| \sum_{\lambda\vdash n}
\frac{\prod_{i=1}^r\alpha_{\chi_\lambda}^{(a_i)} \prod_{j=1}^s m_{\chi_\lambda}^{(e_j)}}
{(\chi_\lambda(1))^{r+s+2t-2}}.
\end{equation}
We concentrate on the sum over characters. Let $A>0$ be given, and split the sum into three
parts $\sum_1, \sum_2, \sum_3$, according to whether $\lambda_1\geq n-A$, $\|\lambda\|\geq n-A$, or
$\lambda_1, \|\lambda\|<n-A$. Note that the first two cases are mutually exclusive for $n>2A+1$.
For $\varepsilon>0$ and $n$ sufficiently large, we have by
Propositions~\ref{prop:MainLemFuchs} and \ref{Prop:RootNumberEst} (i)
\begin{eqnarray*}
\sum\nolimits_3 & < & \sum_{\underset{\|\lambda\|, \lambda_1<n-A}{\lambda\vdash n}}
\big(\chi_\lambda(1)\big)^{\sum_i 1/a_i + \sum_j (1-2/e_j) -r-s-2t+2+\varepsilon}\\
 & \leq & 2\sum_{\underset{\|\lambda\|\leq \lambda_1<n-A}{\lambda\vdash n}}\big(\chi_\lambda(1)\big)^{-\alpha(\Gamma)
 +\varepsilon}.
\end{eqnarray*}
If $\lambda_1>3n/4$, Lemma~\ref{Lem:chi1lower} (ii) gives $\chi_\lambda(1)>{{\lambda_1}\choose{n-\lambda_1}}$, thus
\[
\sum_{\underset{3n/4<\lambda_1<n-A}{\lambda\vdash n}}\big(\chi_\lambda(1)\big)^{-\alpha(\Gamma)+\varepsilon}
<\sum_{3n/4<\nu<n-A}{\nu\choose{n-\nu}}^{-\alpha(\Gamma)+\varepsilon} p(n-\nu) \ll
{{n-A}\choose A}^{-\alpha(\Gamma)+\varepsilon};
\]
whereas for $\lambda_1\leq 3n/4$, Lemma~\ref{Lem:chi1lower} (i) implies $\chi_\lambda(1)>2^{n/8}$, hence
\[
\sum_{\underset{\|\lambda\|\leq\lambda_1\leq 3n/4}{\lambda\vdash n}} \big(\chi_\lambda(1)\big)^{-\alpha(\Gamma)+\varepsilon}
< 2^{-n/8} p(n) < 2^{-n/9}.
\]
We conclude that $\sum_3\ll n^{-\alpha(\Gamma)A+\varepsilon}$.\\[3mm]
Next, we consider $\sum_2$. Suppose that $\delta=2$. Then every permutation of order $a_i$, every
$e_j$-th power of a permutation and every commutator is even, hence
$\chi_\lambda(\pi)=\chi_{\lambda'}(\pi)$ for $\pi$ with $\pi^{a_i}=1$, and $m_{\chi_\lambda}^{(e_j)} =
m_{\chi_{\lambda'}}^{(e_j)}$, and we obtain $\sum_2=\sum_1$ in this case. If, on the other hand,
$\delta=1$, then either there is some $i$ such that $a_i$ is even, or there is some $j$ such
that $e_j$ is odd. In the first case, write
\begin{eqnarray*}
\alpha_{\chi_\lambda}^{(a_i)} & = & \frac{1}{|\Hom(C_{a_i}, S_n)|}\sum_{\pi^{a_i}=1}\chi_\lambda(\pi)\\[1.5mm]
 & = & \frac{1}{|\Hom(C_{a_i}, S_n)|}\sum_{\pi^{a_i}=1} \sum_{\sum_{t\leq \Delta} te_t=\Delta}
 \gamma_{e_1, \ldots, e_\Delta}
\prod_{t|a_i} (-1)^{(t-1)s_t(\pi)} \big(s_t(\pi)\big)^{e_t},
\end{eqnarray*}
where $\Delta=n-\|\lambda\|$, and the coefficients $\gamma_{e_1, \ldots, e_\Delta}$
are determined via Lemma~\ref{Lem:CharPolynom}. Interchanging summations and
applying Lemma~\ref{Lem:CycleDistAsymp} (ii), we see that $\alpha_{\chi_\lambda}^{(a_i)}<e^{-cn^{1/a_i}}$,
and $\sum_2$ is negligible. In the second case, we have $m_{\chi_\lambda}^{(e_i)}=0$ for some
$i$ and $n$ sufficiently large by Proposition~\ref{Prop:RootNumberEst} (iii), hence $\sum_3$ vanishes.
By what we have shown so far,
\[
\sum\nolimits_1+\sum\nolimits_2+\sum\nolimits_3 = \delta\sum\nolimits_1\hspace{.8mm}+\hspace{.8mm}
\mathcal{O}(n^{-A\alpha(\Gamma)+\varepsilon}).
\]
To deal with
$\sum_1$, we fix a partition $\lambda\vdash n$ with $\lambda_1\geq
n-A$. Then $\prod_{j=1}^s m_{\chi_\lambda}^{(e_j)}$ is ultimately
constant, and $\chi_\lambda(1)$ is a polynomial in $n$ of degree
$n-\lambda_1$. Using Lemmas~\ref{Lem:CharPolynom} and \ref{Lem:CycleDistAsymp} (i),
we compute
\begin{eqnarray*}
\alpha_{\chi_\lambda}^{(a_i)} & = &\frac{1}{|\Hom(C_{a_i}, S_n)|}\sum_{\pi^{a_i}=1}\chi_\lambda(\pi)\\[1.5mm]
 & = & \frac{1}{|\Hom(C_{a_i}, S_n)|}\sum_{\pi^{a_i}=1} \sum_{\sum_{t|a_i} te_t\leq\Delta} \gamma_{e_1, \ldots, e_{a_i}}
\prod_{t|a_i} \big(s_t(\pi)\big)^{e_t}\\[1.5mm]
 & = & \sum_{\sum_{t|a_i} te_t\leq\Delta}\gamma_{e_1, \ldots, e_{a_i}}\sum_{k\leq\sum_{t|a_i}te_t}
 \alpha_{e_1, \ldots, e_{a_i}}^{(k)} \frac{n!}{(n-k)!}\frac{|\Hom(C_{a_i}, S_{n-k})|}{\Hom(C_{a_i}, S_n)|}.
\end{eqnarray*}
By \cite[Eq. (22)]{MInvent} we have, for every finite group $G$ and each fixed $k$,
an asymptotic expansion
\begin{equation}
\label{Eq:HomQuotAsymp}
\frac{|\Hom(G, S_n)|}{|\Hom(G, S_{n-k})|} \approx n^{k(1-1/m)}
\exp\bigg(\sum_{\nu=1}^\infty Q_\nu^{(k)} n^{-\nu/m}\bigg),\quad(n\rightarrow\infty),
\end{equation}
where the coefficients $Q_\nu^{(k)}$ are given in \cite{MInvent}
after Equation (22). Putting $G=C_{a_i}$, we find that
\[
\alpha_{\chi_\lambda}^{(a_i)} \approx
n^{\Delta/a_i}\bigg(\sum_{\nu=0}^\infty A_\nu^{\lambda, a_i} n^{-\nu/a_i}\bigg),\quad n\rightarrow\infty.
\]
Inserting these results into (\ref{Eq:betaclaim}), we obtain an
asymptotic formula
\[
\beta_{\chi_\lambda}^{(R)} = (n!)^{s+2t}n^{(1-l(r+s+2t)+\sum_i
1/a_i)\Delta}\bigg(\sum_{\nu=0}^\infty B_\nu^{\lambda, R}\hspace{.3mm}n^{-\nu/{m_\Gamma}}\bigg),\quad n\rightarrow\infty.
\]
For fixed $A$, there are only finitely many partitions
$\lambda\vdash n$ with $\Delta\leq A$, hence, we obtain an
asymptotic expansion for $\sum_1$ with leading term given by the
partition $\lambda=(n)$. In this case,
$\alpha_{\chi_{(n)}}^{(a_i)}=m_{\chi_{(n)}}^{(e_i)}=1$, thus,
$\beta_{\chi_{(n)}}^{(R)} = (n!)^{s+2t}$, and therefore, as $n\rightarrow\infty$,
\[
|\Hom(\Gamma, S_n)| = N_R(1) \approx \delta
(n!)^{s+2t-1}\prod_{i=1}^r |\Hom(C_{a_i}, S_n)|
\bigg(1\hspace{.7mm}+\hspace{.7mm}\sum_{\nu=1}^\infty C_{\nu}(\Gamma)\hspace{.3mm}n^{-\nu/m_\Gamma}\bigg).
\]
Using the asymptotic expansion \cite[Theorem 5]{MComb} for
$|\Hom(G, S_n)|$, the main term of the last expression is found to
be
\[
\delta L_\Gamma (n!)^{\mu(\Gamma)+1} n^{-1} \Phi_\Gamma(n),
\]
with $\mu(\Gamma), L_\Gamma, \Phi_\Gamma$ and $\delta$ as defined
above. Mimicking the proof of \cite[Prop.~1]{MInvent} and using
our assumption that $\mu(\Gamma)\geq \alpha(\Gamma)>0$, we find
that, for each fixed $K\geq 1$,
\[
\sum_{k=K}^{n-K} \frac{h_k(\Gamma)h_{n-k}(\Gamma)}{h_n(\Gamma)}
\ll n^{-K\mu(\Gamma)}.
\]
Combining the latter estimate with the Lemma from \cite[Sec.
3]{MInvent} (cf. also \cite{Bender}) and the transformation formula\footnote{Cf.
\cite[Prop.~1]{DM} or \cite[Prop.~1]{MJLMS}. A far
reaching generalization of this counting principle is found in
\cite{MRepres}.}
\begin{equation}
\label{Eq:transform}
s_n(\Gamma) = nh_n(\Gamma)-\sum_{k=1}^{n-1}h_{n-k}(\Gamma)s_k(\Gamma),
\end{equation}
now gives
\begin{equation}
\label{Eq:DkTransform}
\frac{s_n(\Gamma)}{nh_n(\Gamma)} \approx 1+\sum_{k=1}^\infty
d_k(\Gamma) \frac{h_{n-k}(\Gamma)}{h_n(\Gamma)}, \quad
(n\rightarrow\infty),
\end{equation}
where $d_k(\Gamma)$ is the coefficient of $z^k$ in the formal
power series $\left(\sum_{n\geq 0}h_n(\Gamma) z^n\right)^{-1}$.
Expanding $\frac{h_{n-k}(\Gamma)}{h_n(\Gamma)}$ by means of the
asymptotic formula for $h_n(\Gamma)$ and the Taylor-series of
$\Phi_\Gamma(n)$, the theorem follows.
\end{proofof}
The condition $\alpha(\Gamma)>0$ in Theorem~\ref{Thm:FuchsAsymp}
is essentially necessary. It can be violated in one of two
ways: either $\mu(\Gamma)\leq 0$, or $\mu(\Gamma)>0$, but there
are sufficiently many large values among the $e_j$ to keep
$\alpha(\Gamma)$ small. Here, we deal with the first possibility;
cf. Subsection~\ref{Subsec:OneRel} for the second case.\\[3mm]
The groups $\Gamma$ with a presentation of
the form (\ref{Eq:GammaDef}) and $\mu(\Gamma)\leq 0$ naturally fall
into three classes, according to whether $\mu(\Gamma)=0$, and either $s=0$, or $e_j=2$ for all $j\leq s$;
or $\mu(\Gamma)=0$ and $e_j>2$ for some $j$; or $\mu(\Gamma)<0$. We will show that in each of these cases the assertion
of Theorem~\ref{Thm:FuchsAsymp} fails to hold.\vspace{-3mm}
\begin{enumerate}
\item $\mu(\Gamma)=0$, and either $s=0$ or $e_j=2$ for all $j\leq s$. Then $\Gamma$ is virtually abelian of rank
2, hence, $s_n(\Gamma)\ll n^c$; cf. \cite[Chapter III, Prop. 7.10]{LynSch}. This would certainly contradict the assertion
of Theorem~\ref{Thm:FuchsAsymp}, provided that $r\neq 0$. If $r=0$, we have $\Gamma=\langle x, y\,|\,[x, y]=1\rangle$
or $\Gamma=\langle x, y\,|\,x^2y^2=1\rangle$, and $s_n(\Gamma)=\sigma(n)$ in both cases, whereas Theorem~\ref{Thm:FuchsAsymp} would
predict $s_n(\Gamma)\sim n^{3/2}$.\\[-3mm]
\item $\mu(\Gamma)=0$, and $e_1>2$, say. Then either $r=2$, $a_1=a_2=2$, or $r=0, s=2, t=0$. In the first case,
$\Gamma$ maps homomorphically onto $C_2\ast C_{e_1}$, whereas in the second case, $\Gamma$ maps homomorphically
onto $C_{e_1}\ast C_{e_2}$; that is, in both cases $\Gamma$ maps onto a free product with negative Euler characteristic,
and therefore $s_n(\Gamma)\gg s_n(C_2\ast C_3)\gg (n!)^{1/6}$, while Theorem~\ref{Thm:FuchsAsymp} would predict that
$\Gamma$ is of subexponential growth.\\[-3mm]
\item $\mu(\Gamma)<0$. In this case, $\Gamma$ is finite, and $s_n(\Gamma)$ is ultimately 0, which contradicts the assertion
of Theorem~\ref{Thm:FuchsAsymp} as well.\vspace{-3mm}
\end{enumerate}
\subsection{Computation of the coefficients $a_\nu(\Gamma)$}
\label{Subsec:CoeffComp} In this section we describe how the
coefficients $a_\nu(\Gamma)$ can be computed explicitly. We begin
by giving explicit values to some of the quantities shown to exist
in the previous sections.\\[3mm]
Combining Lemmas~\ref{Lem:CharPolynom},
\ref{Lem:CycleDistAsymp}, and \ref{Lem:CycleAsympExp}, we evaluate
the constants $\alpha_{\chi_\lambda}^{(q)}$ for partitions
$\lambda\vdash n$ with $\Delta\leq2$.
\begin{lemma}
\label{Lem:MainLemFuchsExpl} Let $q\geq 2$ be an integer. Then we
have
\begin{align*}
\alpha_{\chi_{(n-1,1)}}^{(q)} &= |\Hom(C_q, S_n)|^{-1}\sum_{\pi^q=1}\chi_{(n-1, 1)}(\pi) & =&\hspace{2mm}
n\frac{|\Hom(C_q, S_{n-1})|}{|\Hom(C_q, S_n)|},
\intertext{and for $q$ even}
\alpha^{(q)}_{\chi_{(n-2,2)}} &= |\Hom(C_q, S_n)|^{-1}\sum_{\pi^q=1}\chi_{(n-2, 2)}(\pi) & =&\hspace{2mm}
\frac{n^2}{2}\frac{|\Hom(C_q, S_{n-2})|}{|\Hom(C_q, S_n)|}\\[1.5mm]
\alpha^{(q)}_{\chi_{(n-2,1,1)}} &= |\Hom(C_q, S_n)|^{-1}
\sum_{\pi^q=1}\chi_{(n-2, 1, 1)}(\pi) & =&\hspace{2mm} 1 +
\frac{n^2}{2}\frac{|\Hom(C_q, S_{n-2})|}{|\Hom(C_q, S_n)|}\\[1.5mm]
\intertext{whereas for $q$ odd}
\alpha^{(q)}_{\chi_{(n-2,2)}} &= |\Hom(C_q, S_n)|^{-1}\sum_{\pi^q=1}\chi_{(n-2,
2)}(\pi)
& =&\hspace{2mm} -1 + \frac{n(n-1)}{2}\frac{|\Hom(C_q, S_{n-2})|}{|\Hom(C_q, S_n)|},\\[1.5mm]
\alpha^{(q)}_{\chi_{(n-2,1,1)}} &= |\Hom(C_q, S_n)|^{-1}\sum_{\pi^q=1}\chi_{(n-2, 1, 1)}(\pi) & =&\hspace{2mm}
\frac{n(n-1)}{2}\frac{|\Hom(C_q, S_{n-2})|}{|\Hom(C_q, S_n)|}.
\end{align*}
\end{lemma}
We now use (\ref{Eq:HomQuotAsymp}) and \cite[Theorem 6]{MComb}
to compute the first terms of the asymptotic series for $\alpha_{\chi_{(n-2, 2)}}^{(q)}$.
We obtain
\[
\alpha_{\chi_{(n-2, 2)}}^{(q)} = n^{2/q}\big(1+R_q(n^{-1/q})\big)\bigg(1-2\sum_{\nu=1}^{q+3}
\tilde{Q}_\nu^{(q)} n^{-\nu/q}+ S_q(n^{-1/q})\bigg) + \mathcal{O}(n^{-\frac{q+2}{q}}),
\]
where the polynomials $R_q, S_q$ are given as follows.\\
\begin{center}
\begin{tabular}{c|l|l}
$q$ & $R_q(z)$ & $S_q(z)$\\[1mm] \hline
2 & $\frac{1}{2}z^2 -\frac{1}{4}z^3+\frac{1}{8}z^4+\frac{1}{32}z^5$\phantom{\Big|}
&$\frac{3}{4}z^2-\frac{7}{8}z^3+\frac{23}{64}z^4 + \frac{5}{128}z^5$ \\[1.5mm]
3 & $ \frac{1}{3}z^3 -\frac{2}{9}z^5 +\frac{5}{36}z^6$ &
$\frac{1}{3}z^4-\frac{1}{3}z^5+\frac{17}{108}z^6$\\[1.5mm]
4 & $\frac{1}{4}z^4 -\frac{1}{8}z^6$ & $\frac{3}{16}z^4+\frac{3}{8}z^5+\frac{25}{64}z^6+\frac{33}{64}z^7$\\[1.5mm]
5 & $\frac{1}{5}z^5$ & $\frac{3}{25}z^{8}$\\[1.5mm]
6 & $\frac{1}{6}z^6 -\frac{1}{12}z^9$ & $\frac{1}{12}z^6 + \frac{1}{6}z^7+\frac{1}{4}z^8-\frac{55}{216}z^9$\\[1.5mm]
7 & $\frac{1}{7}z^7$ & 0\\[1.5mm]
8, 12 & $\frac{1}{q}z^q$ & $\frac{3}{q^2}(z^q+z^{q+2}+z^{q+3})$\\[1.5mm]
9 & $\frac{1}{9}z^9$ & $\frac{1}{27}z^{12}$\\[1.5mm]
10, 18 & $\frac{1}{q}z^q$ & $\frac{3}{q^2}(z^q+z^{q+3})$\\[1.5mm]
14, 16 & $\frac{1}{q}z^q$ & $\frac{3}{q^2}z^q$\\[1.5mm]
$q\geq 20$, even & $\frac{1}{q}z^q$ & $\frac{3}{q^2}z^q$\\[1.5mm]
$q\geq 11$, odd & $\frac{1}{q}z^q$ & 0
\end{tabular}
\end{center}
Note that by Lemma \ref{Lem:MainLemFuchsExpl}, we have $\alpha_{\chi_{(n-2, 1, 1)}}^{(q)} =
1+\alpha_{\chi_{(n-2, 2)}}^{(q)}$, that is, the asymptotic series above contains all necessary information
for both characters.\\[3mm]
As an example, consider the triangle group
\[
\Gamma=\Big\langle x, y, z\,\big|\,x^2=y^3=z^7=xyz=1\Big\rangle.
\]
Here, $\mu(\Gamma)=\alpha(\Gamma)=\frac{1}{42}$ and $m_\Gamma=42$.
We first show how to determine the contribution of characters $\chi_\lambda$ with
$\Delta>2$. By Lemmas~\ref{Lem:CharPolynom} and \ref{Lem:CycleDistAsymp} (ii),
we see that, for fixed $\Delta$ and $q$ prime,
\[
\alpha_{\chi_\lambda}^{(q)} = \Big(1+\mathcal{O}(n^{-1/q})\Big)
\chi_{\lambda\setminus\lambda_1}(1)\frac{n^{\Delta/q}}{\Delta!}.
\]
Hence, for the triangle group $\Gamma$ we obtain
\[
\beta_{\chi_\lambda}^{(R)} = \Big(1+\mathcal{O}(n^{-1/7})\Big)
\frac{(\chi_{\lambda\setminus\lambda_1}(1))^3 n^{41\Delta/42}}{(\Delta!)^3(\chi_\lambda(1))^2}.
\]
Setting $\pi=1$ in (\ref{Eq:BetaDef}), we obtain
\begin{multline*}
|\Hom(\Gamma, S_n)| = \frac{|\Hom(C_2, S_n)||\Hom(C_3, S_n)||\Hom(C_7, S_n)|}{n!}\\[2mm]
\times\;\Bigg(\underset{\Delta\leq 2}{\sum_{\lambda}}\chi_\lambda(1) \beta_{\chi_\lambda}^{(R)} +
\underset{3\leq \Delta\leq 22}{\sum_{\lambda}}
\frac{(\chi_{\lambda\setminus\lambda_1}(1))^3n^{41\Delta/42}}{(\Delta!)^3\chi_\lambda(1)}
\hspace{.7mm}+\hspace{.7mm}\mathcal{O}(n^{-23/42})\Bigg).
\end{multline*}
Given our previous work in this section, the sum over $\Delta$ can be computed to whatever length
is required. We have cut the sum after the term $\Delta=22$, since this is the smallest precision
bringing to bear all phenomena occurring in such computations at arbitrary scale.
For $\Delta=0$, we have $\beta_{\chi_\lambda}^{(R)}=1$, whereas for $1\leq\Delta\leq 2$, we use the
asymptotic for $\alpha_{\chi_\lambda}^{(q)}$ computed above, to obtain
\[
\underset{\Delta\leq 2}{\sum_{\lambda}}\chi_\lambda(1)\beta_{\chi_\lambda}^{(R)} =
1+ n^{-1/42} + \frac{1}{2} n^{-1/21} + \frac{1}{2} n^{-11/21}\hspace{.7mm}+\hspace{.7mm}\mathcal{O}(n^{-23/42}).
\]
From Lemma~\ref{Lem:CharPolynom}, we obtain, for $\Delta$ fixed and $n\rightarrow\infty$, the asymptotic estimate
\[
\chi_\lambda(1) = \Big(1+\mathcal{O}(n^{-1})\Big)\frac{\chi_{\lambda\setminus\lambda_1}(1) n^{\Delta}}{\Delta!},
\]
which implies
\begin{eqnarray*}
\underset{3\leq \Delta\leq 22}{\sum_{\lambda}}
\frac{(\chi_{\lambda\setminus\lambda_1}(1))^3 n^{41\Delta/42}}{(\Delta!)^3\chi_\lambda(1)} &=&
\Big(1+\mathcal{O}(n^{-1})\Big) \underset{3\leq \Delta\leq 22}{\sum_{\lambda}}
\frac{(\chi_{\lambda\setminus\lambda_1}(1))^2 n^{-\Delta/42}}{(\Delta!)^2}\\
&=& \Big(1 + \mathcal{O}(n^{-1})\Big) \sum_{3\leq \Delta\leq 22}\frac{n^{-\Delta/42}}{\Delta!};
\end{eqnarray*}
and combining these estimates we obtain
\begin{multline*}
|\Hom(\Gamma, S_n)| = \frac{1}{n!}|\Hom(C_2, S_n)| |\Hom(C_3, S_n)| |\Hom(C_7, S_n)|\\[1mm]
\times\;
\bigg(1 + \sum_{\Delta=1}^{22}\frac{n^{-\Delta/42}}{\Delta!} + \frac{1}{2}n^{-11/21} +
\mathcal{O}(n^{-23/42})\bigg).
\end{multline*}
For $1\leq k\leq 22$, we compute $s_k(\Gamma)$ using the software
package GAP \cite{GAP}, and obtain
\[
s_k(\Gamma) = \begin{cases} 1, & k=1, 8, 9\\
 2, & k=7\\ 3, & k=15\\ 9, & k=14, 21\\ 13, & k=22\\ 0,& \mbox{otherwise,}
\end{cases}\qquad 1\leq k\leq 22.
\]
From these values, $h_k(\Gamma)$ and hence $d_k(\Gamma)$ are easily computed for
$1\leq k\leq 22$. The first three coefficients of the asymptotic series for $|\Hom(C_q, S_n)|$ are
given in \cite{MComb}, after Corollary 2. Putting together the various expansions, our final result is that
\begin{multline*}
s_n(\Gamma) = \frac{(2\pi)^{-\frac{53}{21}}e^{-\frac{1}{4}}}{\sqrt{42}} (n!)^{\frac{1}{42}}
n^{-\frac{11}{21}}\exp\big(n^{1/2}+n^{1/3}+n^{1/7}\big)\\[1mm]
\times\,\bigg\{1-\frac{2}{7}n^{-1/6}-\frac{1}{8}n^{-4/21}
-\frac{1}{9}n^{-3/14}-\frac{113}{147}n^{-1/3} -\frac{23}{140}n^{-5/14} + \frac{319}{8064}n^{8/21}\\[1mm]
+\frac{1}{72}n^{-17/42}+\frac{1}{162}n^{-3/7} + \frac{745}{8232}n^{-1/2}-\frac{28309}{64680}n^{11/21}
\hspace{.7mm}+\hspace{.7mm}\mathcal{O}(n^{-23/42})\bigg\}.
\end{multline*}
\subsection{One-relator groups}
\label{Subsec:OneRel}
The result of this subsection, apart from its inherent interest, also demonstrates
that certain Fuchsian groups $\Gamma$ with $\mu(\Gamma)>0$ and $\alpha(\Gamma)<0$
have a much faster growth than would be predicted by Theorem~\ref{Thm:FuchsAsymp}. Consider
a one-relator group
\begin{equation}
\label{Eq:OneRelType}
\Gamma = \Big\langle y_1, y_2,\ldots, y_s\hspace{.8mm}\big|\hspace{.8mm} y_1^{e_1} y_2^{e_2}\cdots y_s^{e_s} = 1\Big\rangle,
\end{equation}
and let $\bar{\Gamma}:= C_{e_1}\ast C_{e_2}\ast\cdots\ast C_{e_s}$.
\begin{theorem}
\label{Thm:OneRelAsymp}
Let $\Gamma$ be as in {\em (\ref{Eq:OneRelType}),} and suppose that
\[
\alpha(\Gamma) = -2\hspace{.7mm}+\hspace{.7mm}\sum_{1\leq j\leq s}\frac{2}{e_j} < 0.
\]
Then, as $n$ tends to infinity, we have
\[
s_n(\Gamma) \sim s_n(\bar{\Gamma}) \sim K (n!)^{\mu(\Gamma)-\alpha(\Gamma)/2}
\exp\Bigg(\sum_{j=1}^s
\underset{\nu<e_j}{\sum_{\nu\mid e_j}} \frac{n^{\nu/e_j}}{\nu}\hspace{.8mm}+\hspace{.8mm}
\frac{\alpha(\Gamma)-2\mu(\Gamma)+2}{4} \log n\Bigg),
\]
where
\[
K = \frac{\exp\Bigg(-\underset{2\mid e_j}{\sum\limits_j}
(2 e_j)^{-1}\Bigg)}{\left(2 \pi\right)^{\frac{2 + 2\mu(\Gamma)-\alpha(\Gamma)}{4}}
\sqrt{e_1 e_2\cdots e_s}}.
\]
\end{theorem}
\begin{proof}
We have
\begin{eqnarray}
\label{Eq:OneRelHom}
|\Hom(\Gamma, S_n)| &=& \Big|\Big\{(\pi_1,\ldots,\pi_s)\in S_n^s:
\hspace{1mm}\pi_1^{e_1} \pi_2^{e_2}\cdots \pi_s^{e_s} = 1\Big\}\Big|\nonumber\\[1.5mm]
&=& \sum_{{\bf c}_1,\ldots,{\bf c}_s} r_{e_1}
({\bf c}_1) r_{e_2}({\bf c}_2) \cdots r_{e_s}({\bf c}_s) N({\bf c}_1,\ldots, {\bf c}_s),
\end{eqnarray}
where
\[
N({\bf c}_1,\ldots, {\bf c}_s):= \big|\Big\{(\pi_1,\ldots,\pi_s)\in S_n^s:\hspace{1mm}
\pi_1\pi_2\cdots\pi_s=1,\, \pi_i\in{\bf c}_i\,(1\leq i\leq s)\Big\}\Big|.
\]
The contribution in (\ref{Eq:OneRelHom}) of the term corresponding to
${\bf c}_1 = {\bf c}_2 = \cdots = {\bf c}_s = 1$ is
\[
|\Hom(C_{e_1}, S_n)|\cdot|\Hom(C_{e_2}, S_n)|\cdots|\Hom(C_{e_s}, S_n)| = |\Hom(\bar{\Gamma}, S_n)|.
\]
We shall show that the sum over the remaining terms is of lesser order of magnitude. Since
$N({\bf c}_1,\ldots, {\bf c}_s)$ is invariant under permutation of its arguments, we may
assume that $|{\bf c}_1|\geq|{\bf c}_j|$ for all $j$. Hence,
\begin{eqnarray*}
N({\bf c}_1,\ldots,{\bf c}_s) &=& \Big|\Big\{(\pi_2,\ldots,\pi_s)\in S_n^{s-1}:
\hspace{1mm}\pi_2\pi_3\cdots\pi_s\in{\bf c}_1,\,\pi_j\in{\bf c}_j\,(2\leq j\leq s)\Big\}\Big|\\[1.5mm]
&\leq& |{\bf c}_2|\cdot|{\bf c}_3|\cdots|{\bf c}_s|\\[1.5mm]
&\leq& \prod_{1\leq j\leq s}|{\bf c}_j|^{1-\delta_j}
\end{eqnarray*}
for any choice of non-negative real numbers $\delta_1,\delta_2\ldots,\delta_s$ such that $\sum_j\delta_j=1$.
 For a non-empty set $J\subseteq[s]$, define $S_J:= \sum_{j\in J}1/e_j$. By our assumption,
 $S_J<1$, and, by definition, $\sum_{j\in J}1/(S_J e_j) = 1$. Using the above estimate for $N({\bf c}_1,\ldots,{\bf c}_s)$
 with $\delta_j=1/(S_J e_j)$, dividing Equation (\ref{Eq:OneRelHom}) by $|\Hom(\bar{\Gamma}, S_n)|$,
 and interchanging product and sum, we find that
\begin{eqnarray*}
0\leq\frac{|\Hom(\Gamma,S_n)|}{|\Hom(\bar{\Gamma},S_n)|} - 1 &\leq&
\underset{\emptyset\not= J\subseteq [s]}{\sum_J}
\hspace{.8mm}\underset{{\bf c}_j\neq1\Leftrightarrow j\in
J}{\sum_{{\bf c}_1,\ldots,{\bf c}_s}} \prod_{j\in J} \frac{r_{e_j}({\bf c}_j)
|{\bf c}_j|^{1-1/(S_J e_j)}}{|\Hom(C_{e_j}, S_n)|}\\[1.5mm]
&=&\underset{\emptyset\not= J\subseteq [s]}{\sum_J} \prod_{j\in J}
\sum_{{\bf c}\neq1} \frac{r_{e_j}({\bf c}) |{\bf c}|^{1-1/(S_J
e_j)}}{|\Hom(C_{e_j}, S_n)|}.
\end{eqnarray*}
Consider the factor in the last expression corresponding to $j\in J$ for a given set $J$.
Classifying the conjugacy classes according to the number $\ell$ of points moved by each of its elements,
this factor can be written as
\begin{multline}
\label{Eq:OneRelEq2}
\sum_{2\leq\ell\leq n} \underset{\mbox{\tiny
${\bf c}$ moves $\ell$ points}}{\sum_{\bf c}} \frac{r_{e_j}({\bf c}) |{\bf c}|^{1-1/(S_J
e_j)}}{|\Hom(C_{e_j}, S_n)|} = \\[1.5mm]
\sum_{2\leq\ell\leq n} \binom{n}{\ell}^{1-1/(S_J e_j)}
\frac{|\Hom(C_{e_j}, S_{n-\ell})|}{|\Hom(C_{e_j}, S_n)|}
{\sum_{\bf c}}^\ast r_{e_j}({\bf c}) |{\bf c}|^{1-1/(S_J e_j)},
\end{multline}
where the innermost sum extends over all fixed-point free
conjugacy classes of $S_\ell$. {}From \cite[Corollary~2]{MComb} we
deduce that
\[
\frac{|\Hom(C_{e_j}, S_{n-\ell})|}{|\Hom(C_{e_j}, S_n)|} \ll \bigg[\binom{n}{\ell} \ell!\bigg]^{-(1-1/e_j)}.
\]
Applying the Cauchy-Schwarz inequality, we find for the innermost sum in (\ref{Eq:OneRelEq2}) that
\begin{eqnarray*}
{\sum_{\bf c}}^\ast r_{e_j}({\bf c}) |{\bf c}|^{1-1/(S_J e_j)} &\leq&
\bigg(\sum_{\bf c} \big(r_{e_j}({\bf c})\big)^2 |{\bf c}|\bigg)^{1/2}
\bigg(\sum_{\bf c} |{\bf c}|^{1-2/(S_J e_j)}\bigg)^{1/2}\\[1.5mm]
&\leq& \bigg(\sum_{\chi\in\mbox{\tiny Irr}(S_\ell)} \big(\chi(1)\big)^{2-\frac{4}{e_j}+
\varepsilon}\bigg)^{1/2} \bigg(\sum_{\bf c} |{\bf c}|^{1-2/(S_J e_j)}\bigg)^{1/2},
\end{eqnarray*}
where we have used Proposition~\ref{Prop:RootNumberEst} (i) to estimate the first factor. If
$S_J e_j>2$, then we bound the second factor by $(\ell!)^{1-2/(S_J e_j)+\varepsilon}$; otherwise, each summand is $\leq 1$, and
\[
\sum_{\bf c} |{\bf c}|^{1-2/(S_J e_j)} \ll (\ell!)^\varepsilon.
\]
In both cases,
\begin{equation}
\label{Eq:OneRelEq3}
\sum_{\bf c} r_{e_j}({\bf c})
|{\bf c}|^{1-1/(S_J e_j)} \leq \max\Big(\big(\ell!\big)^{1-\frac{1}{e_j}-\frac{1}{S_J e_j}+\varepsilon},
\big(\ell!\big)^{\frac{1}{2}-\frac{1}{e_j}+\varepsilon}\Big).
\end{equation}
Putting (\ref{Eq:OneRelEq3}) back into (\ref{Eq:OneRelEq2}), we find that the left-hand side of
(\ref{Eq:OneRelEq2}) is bounded above by
\begin{multline*}
\sum_{2\leq\ell\leq n} \binom{n}{\ell}^{1-1/(S_J e_j)} \bigg[\binom{n}{\ell}
\ell!\bigg]^{-(1-1/e_j)}\max\Big(\big(\ell!\big)^{1-\frac{1}{e_j}-\frac{1}{S_J e_j}+\varepsilon},
\big(\ell!\big)^{\frac{1}{2}-\frac{1}{e_j}+\varepsilon}\Big) =\\[1mm]
\sum_{2\leq\ell\leq n} \binom{n}{\ell}^{\frac{1}{e_j}-\frac{1}{S_J e_j}}
\max\Big(\big(\ell!\big)^{-\frac{1}{S_J e_j}+\varepsilon}, \big(\ell!\big)^{-\frac{1}{2}+\varepsilon}\Big).
\end{multline*}
For sufficiently large $n$, increasing $\ell$ by $1$ decreases a summand by at least a factor $\frac{1}{2}$;
hence, as $n\rightarrow\infty$,
\[
\sum_{2\leq\ell\leq n}\hspace{.6mm}\underset{\mbox{\tiny ${\bf c}$ moves $\ell$ points}}{\sum_{\bf c}}
\frac{r_{e_j}({\bf c}) |{\bf c}|^{1-1/(S_J e_j)}}{|\Hom(C_{e_j}, S_n)|} \ll
n^{-2(\frac{1}{S_J e_j}-\frac{1}{e_j})},
\]
and, therefore,
\[
|\Hom(\Gamma, S_n)| \sim |\Hom(\bar{\Gamma}, S_n)|,\quad n\rightarrow\infty.
\]
Since $\chi(\bar{\Gamma}) = \sum_j 1/e_j - s +1<0$, and since the
proof of \cite[Proposition~1]{MInvent} only depends on an
asymptotic estimate, we infer in particular that
\begin{equation*}
\sum_{0< k<n} \binom{n}{k} \frac{|\Hom(\Gamma,S_k)| |\Hom(\Gamma,
S_{n-k})|}{|\Hom(\Gamma, S_n)|} \rightarrow
0\hspace{2mm}\mbox{as}\hspace{2mm}n\rightarrow\infty.
\end{equation*}
Combining this fact with the transformation formula (\ref{Eq:transform})
and \cite[Theorem~3]{Wright}, we find that
\begin{equation*}
s_n(\Gamma) \sim |\Hom(\Gamma, S_n)|/(n-1)! \sim
|\Hom(\bar{\Gamma}, S_n)|/(n-1)! \sim s_n(\bar{\Gamma}).
\end{equation*}
The explicit asymptotic formula given for $s_n(\bar{\Gamma})$
results from \cite[Theorem~1]{MInvent}.
\end{proof}
{\bf Remark.} Theorem~\ref{Thm:OneRelAsymp} was proved in \cite[Example~1(ii)]{MuPuEx}
for $\alpha(\Gamma)<-1$ (note that the invariant $\alpha(\Gamma)$ defined
in \cite[Sec.~2]{MuPuEx} is $\frac{\alpha(\Gamma)+2}{2}$ in our present notation).
\subsection{Demu\v skin groups}
\label{subsec:Demu}
Let $p$ be a prime. A pro-$p$-group $\Gamma$ is termed a Poincar\'e group
of dimension $n$, if $\Gamma$ has cohomological dimension $n$, and the algebra
$H^*(\Gamma)$ is finite dimensional and satisfies Poincar\'e duality. A Poincar\'e group
of dimension 2 is called a Demu\v skin group. These are one-relator groups
\[
\Gamma=\big\langle x_1, \ldots, x_m\,|\,R(x_1, \ldots, x_m)=1\big\rangle,\quad m=\dim H^1(\Gamma),
\]
and, for $p\neq 2$, the defining relation may be taken to be
\[
R=x_1^{p^h}[x_1, x_2][x_3, x_4]\cdots [x_{m-1}, x_m],\quad h\in\NN\cup\{\infty\},
\]
with the understanding that $x_1^{p^h}=1$ if $h=\infty$; cf \cite{Dem1}--\cite{Dem3} and
\cite{Labute}. For $m=2$, these groups are metacyclic, and their subgroup growth can
be computed using the methods of \cite{GSS}. The similarity of $R$ to a surface
group relation would allow us to estimate $s_n(\Gamma)$ asymptotically for $m\geq 6$,
using only tools from \cite{MuPuSurf}. The case $m=4$ however needs a more careful analysis.
As another application of our estimates for multiplicities of root number functions,
we prove the following result.
\begin{theorem}
\label{Thm:Demush}
For integers $q\geq 1$ and $d\geq 2,$ let
\[
\Gamma_{q, d} = \Big\langle x_1, y_1, \ldots, x_d, y_d\,\big|\,x_1^{q-1}[x_1, y_1] [x_2,y_2]\cdots[x_d, y_d]=1\Big\rangle.
\]
Then there exist explicitly computable constants $\gamma_\nu(\Gamma_{q, d}),$ such that
\[
s_n(\Gamma_{q, d}) \approx \delta n (n!)^{2d-2}\left\{1+\sum_{\nu=1}^\infty
\gamma_{\nu}(\Gamma_{q, d})\hspace{.3mm}n^{-\nu}\right\},\quad (n\rightarrow\infty),
\]
where
\[
\delta=\begin{cases}1, &q\mbox{\em\ even}\\2, &q\mbox{\em\ odd.}\end{cases}
\]
\end{theorem}
The proof runs parallel to the proof of Theorem~\ref{Thm:FuchsAsymp},
once we have established the following.
\begin{lemma}
Let $q\geq 2$ an integer. For a partition $\lambda\vdash n,$ define
the coefficient $l_{\chi_\lambda}^{(q)}$ by means of the equation
\[
\big|\big\{(\sigma, \tau)\in S_n^2: \sigma^{q-1}[\sigma, \tau] = \pi\big\}\big| = n!\sum_{\lambda\vdash n}l_{\chi_\lambda}^{(q)}
\chi_\lambda(\pi).
\]
Then we have $|l_{\chi_\lambda}^{(q)}|\leq \sqrt{\frac{m_{\chi_\lambda}^{(q)}}{\chi_\lambda(1)}}$.
Moreover, for a partition $\mu\vdash l,$ and a partition $\lambda\vdash n$
with $\lambda\setminus\lambda_1=\mu,$ the quantity $\chi_\lambda(1)\hspace{.3mm}l_{\chi_\lambda}^{(q)}$
is a constant depending only on $\mu,$ provided $n$ is sufficiently large.
\end{lemma}
\begin{proof}
Writing the equation $x^{k-1}[x, y] = \pi$ as $x^k (x^{-1})^y=\pi$, we
see that the number of solutions can be computed as
\[
\sum_\mathbf{c} |C_{S_n}(\mathbf{c})|\cdot\big|\big\{\sigma, \tau\in \mathbf{c}:\, \sigma^n \tau=\pi\big\}\big| =
\sum_\mathbf{c} |\mathbf{c}|\sum_{\lambda\vdash n}
\frac{\chi_\lambda(\mathbf{c})\chi_\lambda(\mathbf{c}^k)\chi_\lambda(\pi)}{\chi_\lambda(1)},
\]
that is,
\begin{eqnarray*}
l_{\chi_\lambda}^{(q)} & = & \frac{1}{n!}\sum_\mathbf{c}|\mathbf{c}|
\frac{\chi_\lambda(\mathbf{c})\chi_\lambda(\mathbf{c}^q)}{\chi_\lambda(1)}\\
 & \leq & \frac{1}{n!\chi_\lambda(1)}\left(\sum_\mathbf{c} |\mathbf{c}|(\chi_\lambda(\mathbf{c}))^2\right)^{1/2}
\left(\sum_\mathbf{c} |\mathbf{c}|(\chi_\lambda(\mathbf{c}^q))^2\right)^{1/2}\\
 & \leq & \frac{1}{\sqrt{\chi_\lambda(1)}}\left(\frac{1}{n!}
 \sum_\mathbf{c} |\mathbf{c}|\chi_\lambda(\mathbf{c}^q)\right)^{1/2}\\
 & = & \sqrt{\frac{m_{\chi_\lambda}^{(q)}}{\chi_\lambda(1)}}.
\end{eqnarray*}
For the second claim we have to compute
\[
\frac{1}{n!}\sum_{\pi\in S_n} \frac{\chi_\lambda(\pi)\chi_\lambda(\pi^k)}{\chi_\lambda(1)}.
\]
We express $\chi_\lambda(\pi)$ as a polynomial in the functions $s_i(\pi)$.
Then $\chi_\lambda(\pi)\chi_\lambda(\pi^k)$ is also a
polynomial in these functions, and our claim follows from Lemma~\ref{Lem:Poisson}.
\end{proof}
As an example, consider $\lambda=(n-1, 1)$. Then
\[
\chi_\lambda(\pi)\chi_\lambda(\pi^k)= (s_1(\pi)-1)\bigg(\sum_{t|q}ts_t(\pi)-1\bigg).
\]
The expected value of the first factor is 0, and it is stochastically
independent of $s_2, \ldots, s_q$. Hence, we have
\[
l_{\chi_\lambda}^{(q)} = \frac{1}{\chi_\lambda(1)\hspace{.3mm}n!}\sum_{\pi\in S_n}\big((s_1(\pi))^2 -
s_1(\pi)\big),
\]
and the computations leading to the second assertion in
Proposition~\ref{Prop:RootNumberEst} (ii) give
$l_{\chi_\lambda}^{(q)} = \frac{1}{\chi_\lambda(1)}=\frac{1}{n-1}$.\\[3mm]
For $\lambda=(n-2, 2)$, we find in the case $q$ even
\begin{multline*}
\chi_\lambda(\pi)\chi_\lambda(\pi^q) = \frac{1}{4}(s_1(\pi))^4-\frac{3}{2}(s_1(\pi))^3
+\frac{13}{4}(s_1(\pi))^2 - 3s_1(\pi) +1\\[2mm]
+s_2(\pi)\big(-2(s_1(\pi))^2+6s_1(\pi)-4\big) + (s_2(\pi))^2\big((s_1(\pi))^2-3s_1(\pi)+5\big)-2(s_2(\pi))^3\\[2mm]
+\underbrace{\Big(\frac{1}{2}(s_1(\pi))^2-\frac{3}{2}s_1(\pi)-s_2(\pi)+1\Big)}_{=:A}
\Bigg(\frac{1}{2}\underset{t\geq 3}{\sum_{t|q}} t^2 (s_t(\pi))^2 - \frac{3}{2}\underset{t\geq 3}{\sum_{t|q}} t s_t(\pi)
+\underset{t\geq 3}{\underset{t=2(t, q)}{\sum_{t|2q}}} t s_t(\pi)\Bigg).
\end{multline*}
Note that the expected value
of $A$ is 0, hence, since the second factor contains only terms $s_t(\pi)$ with $t\geq 3$,
it is stochastically independent of the first factor, and the expectation of the last summand vanishes;
using Lemma~\ref{Lem:Poisson}, we find that the remaining terms vanish as well. Dealing in a similar way with
the other cases, we obtain
\begin{center}
\begin{tabular}{c|c|c}
 & $q$ even & $q$ odd\\[1mm] \hline
$l_{\chi_{(n-2, 2)}}^{(q)}$ & 0 & $\displaystyle\frac{1}{n^2-3n}\phantom{\Bigg|}$\\[1mm]
$l_{\chi_{(n-2, 1, 1)}}^{(q)}$ & $\displaystyle\frac{13}{2(n^2-3n+2)}$ & $\displaystyle\frac{9}{2(n^2-3n+2)}$
\end{tabular}
\end{center}
As an example for the computation of the coefficients $\gamma_\nu(\Gamma_{q, d})$, we
consider the case $d=2$ and $q\geq 2$. From the values given above we deduce the estimate
\[
\frac{h_n(\Gamma)}{(n!)^2} = 1+\frac{1}{(n-1)^3} + \left.\begin{cases}\frac{26}{(n^2-3n+2)^3}, & q\mbox{ even}\\[1mm]
\frac{4}{(n^2-3n)^3} + \frac{18}{(n^2-3n+2)^3}, & q\mbox{ odd.}\end{cases}\right\}\hspace{.8mm}+\hspace{.8mm}\mathcal{O}(n^{-9}).
\]
For small values of $k$, we compute $h_k(\Gamma_{q, 2})$ as follows:
\[
h_1(\Gamma_{q, 2})=1,\quad h_2(\Gamma_{q, 2}) =\begin{cases} 8,& (q, 2)=1\\ 4,& (q, 2)=2,\end{cases}\quad
h_3(\Gamma_{q, 2})=\begin{cases} 72, & (q, 6)=1\\ 45, & (q, 6)=2\\ 63, & (q, 6)=3\\ 36, & (q, 6)=6, \end{cases}
\]
\[
h_4(\Gamma_{q,2}) = \begin{cases}1424,&(q,6)=1\\
                                 720,& (q,6)=2\\
                                 1280,& (q,6)=3\\
                                 576,& (q,6)=6,
                    \end{cases}\quad h_5(\Gamma_{q,2}) = \begin{cases}37192,&(q,30)=1\\
                                                                      21092,&(q,30)=2\\
                                                                      36040,&(q,30)=3\\
                                                                      35792,&(q,30)=5\\
                                                                      20840,&(q,30)=6\\
                                                                      19692,&(q,30)=10\\
                                                                      34640,&(q,30)=15\\
                                                                      19440,&(q,30)=30.
                                                         \end{cases}
\]
Proceeding as in Subsection~\ref{Subsec:FuchsGeneric}, we obtain
\[
s_n(\Gamma_{q, 2}) = \delta n(n!)^2R(n),
\]
where $\delta$ is as in Theorem~\ref{Thm:Demush}, and
\[
R(n) =
\begin{cases}
   1-n^{-1}-7n^{-2}-56n^{-3}-1237n^{-4}-33573n^{-5}+\mathcal{O}(n^{-6}),& (q, 30)=1\\
   1-n^{-1}-3n^{-2}-37n^{-3}-\phantom{1}623n^{-4}-19460n^{-5}+\mathcal{O}(n^{-6}),& (q, 30)=2\\
   1-n^{-1}-7n^{-2}-47n^{-3}-1111n^{-4}-32826n^{-5}+\mathcal{O}(n^{-6}),& (q, 30)=3\\
   1-n^{-1}-7n^{-2}-56n^{-3}-1237n^{-4}-32173n^{-5}+\mathcal{O}(n^{-6}),& (q, 30)=5\\
   1-n^{-1}-3n^{-2}-28n^{-3}-\phantom{1}497n^{-4}-19541n^{-5}+\mathcal{O}(n^{-6}),& (q, 30)=6\\
   1-n^{-1}-3n^{-2}-37n^{-3}-\phantom{1}623n^{-4}-18060n^{-5}+\mathcal{O}(n^{-6}),& (q, 30)=10\\
   1-n^{-1}-7n^{-2}-47n^{-3}-1111n^{-4}-31426n^{-5}+\mathcal{O}(n^{-6}),& (q, 30)=15\\
   1-n^{-1}-3n^{-2}-28n^{-3}-\phantom{1}497n^{-4}-18141n^{-5}+\mathcal{O}(n^{-6}),& (q, 30)=30.
\end{cases}
\]
Note that the series for $s_n(\Gamma_{q, 2})$ is far more dependent on $q$ -- and therefore
on $\Gamma_{q, 2}$ itself -- than the series for $h_n(\Gamma_{q, 2})$.
\section{Finiteness Results}
\label{Sec:FuchsianFinite}
Call two finitely generated groups
$\Gamma$ and  $\Delta$ {\it equivalent}, denoted $\Gamma\sim\Delta$, if
\[
s_n(\Gamma)=(1+o(1))s_n(\Delta),\qquad(n\rightarrow\infty).
\]
In \cite[Theorem 3]{MInvent} a characterization in terms of
structural invariants is given for the equivalence relation $\sim$
on the class of groups $\Gamma$ of the form
\[
\Gamma=G_1*G_2*\dots*G_s*F_r
\]
with $r, s\geq 0$ and finite groups $G_\sigma$, and it is shown
that  each $\sim$-class decomposes into finitely many isomorphism
classes. Here we are concerned with the analogous problems for
Fuchsian groups.
\begin{theorem}\label{Thm:FuchsWeakEquiv}
\begin{enumerate}
\item Let $\Gamma=C_{a_1}*\dots*C_{a_k}*F_r$ and $\Delta =
C_{b_1}*\dots*C_{b_l}*F_{r'}$ be free products of cyclic groups
such that $s_n(\Gamma)\asymp s_n(\Delta)$. Then
$r=r'$ and $\{a_1,\ldots,a_k\}=\{b_1,\ldots,b_l\}$ as multisets.\\[-2mm]
\item Let
\begin{multline*}
\Gamma = \Big\langle x_1, \ldots, x_r, y_1, \ldots, y_s, u_1, v_1, \ldots, u_t, v_t\,\big|\,\\
x_1^{a_1}=\cdots=x_r^{a_r} = x_1\cdots x_r y_1^{e_1}\cdots
y_s^{c_s} [u_1, v_1]\cdots [u_t, v_t] = 1\Big\rangle
\end{multline*}
and
\begin{multline*}
\Delta = \Big\langle x_1, \ldots, x_{r'}, y_1, \ldots, y_{s'}, u_1, v_1, \ldots, u_{t'}, v_{t'}\,\big|\,\\
x_1^{a'_1}=\cdots=x_{r'}^{a'_{r'}} = x_1\cdots x_{r'}
y_1^{e'_1}\cdots y_s^{e'_{s'}} [u_1, v_1]\cdots [u_{t'},
v_{t'}] = 1\Big\rangle
\end{multline*}
be Fuchsian groups, such that $\alpha(\Gamma), \alpha(\Delta)>0$.
Then $\Gamma\sim\Delta$ if and only if
\vspace{1mm}
\begin{enumerate}
\item The multisets $\{a_i:1\leq i\leq r\}$ and
$\{a'_i:1\leq i\leq r'\}$ coincide,\\[-3.5mm]
\item $\mu(\Gamma) = \mu(\Delta),$\\[-3.5mm]
\item $\delta = \delta'$.
\end{enumerate}
\item Let
\[
\Gamma = \Big\langle y_1,\ldots,y_s\,\big|\,y_1^{e_1} y_2^{e_2}\cdots y_s^{e_s} = 1\Big\rangle
\]
and
\[
\Delta = \Big\langle y_1,\ldots,y_{s'}\,\big|\,y_1^{e_1'} y_2^{e_2'}\cdots y_{s'}^{e_{s'}'} = 1\Big\rangle
\]
be two one-relator groups with $\alpha(\Gamma),\alpha(\Delta)<0$. Then the following assertions are equivalent:
\vspace{1.5mm}
\begin{enumerate}
\item $\Gamma \sim \Delta.$\\[-3.5mm]
\item $s=s'$ and $\{e_1,\ldots, e_s\}=\{e_1',\ldots,e_{s'}'\}$ as multisets.\\[-3mm]
\item $\hat{\Gamma} \cong \hat{\Delta},$ where the hat denotes profinite completion.
\end{enumerate}
\end{enumerate}
\end{theorem}
The proof of Theorem~\ref{Thm:FuchsWeakEquiv} requires the
following two auxiliary results.
\begin{lemma}
\label{divisorsums}
Let $A=\{a_1, \ldots, a_k\}$ and $B=\{b_1, \ldots, b_l\}$ be
multisets of integers, such that
\[
\sum_{d|a_i} \frac{1}{a_i} = \sum_{d|b_i}\frac{1}{b_i}
\]
for all $d\geq 2$. Then $A=B$.
\end{lemma}
\begin{proof}
We argue by induction on $n=k+l$. For $n\leq 1$, there is nothing
to show. Assume that our claim holds for all multisets $A', B'$
with $|A'|+|B'|\leq n-1$, and let $A, B$ be multisets as above.
Let $d$ be the greatest integer, such that $\sum_{d|a_i}1/a_i>0$.
Then $d=\max A$, and the value of the sum is $|\{i:a_i=\max
A\}|/\max A$. The same holds for $B$, hence the greatest element
of both multisets as well as the multiplicity of this maximum
coincide. Deleting these elements in both multisets yields a pair
of multisets $A', B'$ of smaller cardinality, which are equal by
the induction hypothesis. Hence we deduce $A=B$.
\end{proof}
\begin{lemma}
\label{Lem:GammaHat}
Given positive integers $k$ and $l,$ disjoint tuples of variables $\vec{x}_i = (x_{i1},\ldots,x_{ia_i})$
for $1\leq i\leq k,$ words $w_1(\vec{x}_1),\ldots, w_k(\vec{x}_k),$ and (possibly empty) words
$v_{ij}(\vec{x}_i)$ for $1\leq ~i\leq~k$ and $1\leq j\leq l,$ and a permutation $\sigma\in S_k,$ define
\[
\Gamma = \Big\langle \vec{x}_1, \ldots, \vec{x}_k\,\big|\,w_1(\vec{x}_1)\cdots
w_k(\vec{x}_k) = v_{ij}(\vec{x}_i) = 1\,(1\leq i\leq k,\,1\leq j\leq l)\Big\rangle
\]
and
\[
\Gamma_\sigma = \Big\langle \vec{x}_1, \ldots, \vec{x}_k\,\big|\,w_{\sigma(1)}(\vec{x}_{\sigma(1)})\cdots
w_{\sigma(k)}(\vec{x}_{\sigma(k)}) = v_{ij}(\vec{x}_i) = 1\,(1\leq i\leq k,\,1\leq j\leq l)\Big\rangle.
\]
Then $\Gamma$ and $\Gamma_\sigma$ have isomorphic profinite completions.
\end{lemma}
\begin{proof}
For $1\leq i\leq k$ and a finite group $G$, let
\[
N_i^{(G)}(g) := \Big|\Big\{\vec{x}_i\in G^{a_i}: w_i(\vec{x}_i) = g,
v_{i1}(\vec{x}_i)=\cdots=v_{il}(\vec{x}_i)=1\Big\}\Big|,\quad g\in G.
\]
Since $N_i^{(G)}$ is a class function, we can introduce Fourier coefficients $\alpha_{\chi, i}$ via
\[
N_i^{(G)}(g) = \sum_{\chi\in\mathrm{Irr}(G)}\alpha_{\chi, i}\hspace{.3mm}\chi(g),\quad g\in G.
\]
Then, using orthogonality, we have
\begin{eqnarray*}
|\Hom(\Gamma, G)| & = & \underset{g_1g_2\cdots g_k=1}{\sum_{g_1, \ldots, g_k}}
\prod_{1\leq i\leq k} N_i^{(G)}(g_i)\\[1.5mm]
 & = & \sum_{\mathbf{c}_1, \ldots, \mathbf{c}_k\subseteq G}\frac{|\mathbf{c}_1|\cdots|\mathbf{c}_k|}{|G|}
\sum_\chi\sum_{\chi_1, \ldots, \chi_k}\frac{\chi(\mathbf{c}_1)\cdots\chi(\mathbf{c}_k)}{(\chi(1))^{k-2}}
\hspace{.3mm}\chi_1(\mathbf{c}_1)\hspace{.2mm}\alpha_{\chi_1, 1}\cdots\chi_k(\mathbf{c}_k)
\hspace{.2mm}\alpha_{\chi_k, k}\\[1.5mm]
 & = & |G|^{k-1}\sum_\chi\frac{\alpha_{\overline{\chi}, 1}\cdots\alpha_{\overline{\chi}, k}}{(\chi(1))^{k-2}}.
\end{eqnarray*}
Since this computation leads to the same character formula when replacing $\Gamma$ by $\Gamma_\sigma$, we deduce
that, for each finite group $G$,
\[
|\Hom(\Gamma, G)| = |\Hom(\Gamma_\sigma, G)|;
\]
in particular, $s_n(\Gamma)=s_n(\Gamma_\sigma)$ for all $n$. Writing
\[
|\Hom(\Gamma, G)|=\sum_{U\leq G}|\mbox{Epi}(\Gamma, U)|
\]
and using M\"obius inversion in the subgroup lattice of $G$, this gives
\[
|\mbox{Epi}(\Gamma, G)| = \sum_{U\leq G}\mu(U, G)|\Hom(\Gamma, U)|,
\]
thus also
\begin{equation}
\label{Eq:EqualEpi}
|\mbox{Epi}(\Gamma, G)| = |\mbox{Epi}(\Gamma_\sigma, G)|.
\end{equation}
For $n\in\NN$, define finite groups $G_n$ and $G_n^\sigma$ via
\[
G_n=\Gamma\big/\bigcap\limits_{(\Gamma:\Delta)\leq n} \Delta,\quad
G_n^\sigma=\Gamma_\sigma\big/\bigcap\limits_{(\Gamma_\sigma:\Delta)\leq n} \Delta.
\]
From (\ref{Eq:EqualEpi}) we know in particular, that $G_n$ is a
homomorphic image of $\Gamma_\sigma$. Let $N$ be the kernel of
such a projection map $\phi:\Gamma_\sigma\rightarrow G_n$.
Since $s_\nu(\Gamma_\sigma)=s_\nu(G_n)$ for $\nu\leq n$, we have
\[
N\leq \bigcap_{(\Gamma_\sigma:\Delta)\leq n}\Delta,
\]
hence, $G_n^\sigma$ is a homomorphic image of $G_n$. By symmetry,
$G_n$ and $G_n^\sigma$ are isomorphic. By the universal properties
of $G_n$ and $G_n^\sigma$,
\[
\hat{\Gamma}\cong\lim\limits_\leftarrow G_n \cong \lim\limits_\leftarrow G_n^\sigma \cong \hat{\Gamma}_\sigma,
\]
as claimed.
\end{proof}
\begin{proofof}{Theorem {\rm \ref{Thm:FuchsWeakEquiv}}}
(i) By \cite[Theorem 1]{MInvent}, the assumption
$s_n(\Gamma)\asymp s_n(\Delta)$ is equivalent to the assertion
that
\begin{multline*}
(n!)^{-\chi(\Gamma)}\exp\left(\sum_{i=1}^k\sum_{d|a_i}\frac{n^{d/a_i}}{d}
+\frac{r+\chi(\Gamma)+1}{2}\log n\right)\asymp\\
(n!)^{-\chi(\Delta)}\exp\left(\sum_{i=1}^l\sum_{d|b_i}\frac{n^{d/b_i}}{d}
+\frac{r'+\chi(\Delta)+1}{2}\log n\right).
\end{multline*}
Comparing orders of magnitude as in the proof of \cite[Theorem
3]{MInvent}, we find first that $\chi(\Gamma)=\chi(\Delta)$, then,
successively, that
\[
\sum_{t|a_i} \frac{t}{a_i} = \sum_{t|b_i}\frac{t}{b_i}, \quad
t\geq 2,
\]
and, finally, that $r=r'$. Our claim follows now from Lemma
\ref{divisorsums}.\\[3mm]
(ii) By Theorem \ref{Thm:FuchsAsymp}, the assertion that
$\Gamma\sim\Delta$ is equivalent to
\[
\delta s_n(C_{a_1}*\dots*C_{a_r}*F_{s+2t}) \sim \delta'
s_n(C_{a'_1}*\dots*C_{a'_{r'}}*F_{s'+2t'}),\quad
(n\rightarrow\infty).
\]
By part (i), the latter assertion is equivalent to the conjunction
of
\[
C_{a_1}*\dots*C_{a_r}*F_{s+2t}\cong
C_{a'_1}*\dots*C_{a'_{r'}}*F_{s'+2t'}
\]
and $\delta=\delta'$, whence our claim.\\[3mm]
(iii) The equivalence of (a) and (b) follows from Theorem~\ref{Thm:OneRelAsymp} and part (i). Since (c)
obviously implies (a), it suffices to show that (b) implies (c); but this follows immediately from
Lemma \ref{Lem:GammaHat} upon setting $l=0$ and $w_i(\vec{x}_i) = x_i^{e_i}$.
\end{proofof}
Denote by $\F$ the class of all groups $\Gamma$ having a presentation of the form (\ref{Eq:GammaDef})
with $\alpha(\Gamma)>0$. We distinguish three refinements of the equivalence relation $\sim$ on $\F$:
 (i) the relation $\approx$ of strong equivalence defined via
\[
\Gamma\approx\Delta :\Leftrightarrow
s_n(\Gamma)=s_n(\Delta)(1+\mathcal{O}(n^{-A}))\mbox{ for every }A>0,
\]
(ii) isomorphy, and (iii) equality of the system of parameters
\[
(r,s,t;a_1,a_2\ldots,a_r,e_1,e_2,\ldots,e_s)
\]
in the Fuchsian presentation (\ref{Eq:GammaDef}), denoted $\Gamma=\Delta$. Clearly,
\begin{equation}
\label{Eq:Refine}
\Gamma = \Delta \Rightarrow \Gamma\cong\Delta \Rightarrow \Gamma\approx\Delta \Rightarrow \Gamma\sim\Delta.
\end{equation}
All these implications are in fact strict. To see this, define
\[
\Gamma_j = \Big\langle x,y,z,u\,\big|\,R_j(x,y,z,u)=1\Big\rangle,\quad 1\leq j\leq3,
\]
where
\[
R_j:= \begin{cases}[x,y] [z,u],& j=1\\
                   [x,y] z^2u^2,& j=2\\
                   x^2y^2z^2u^2,& j=3.
     \end{cases}
\]
Then $\Gamma_1$ and $\Gamma_2$ are isospectral, that is, $s_n(\Gamma_1)=s_n(\Gamma_2)$
for all $n$ (in particular, $\Gamma_1\approx \Gamma_2$), but $\Gamma_1\not\cong\Gamma_2$;
and $\Gamma_2\cong\Gamma_3$ but $\Gamma_2\neq\Gamma_3$.
Our next result implies that $\approx$ is a much finer equivalence relation than $\sim$. It appears
that the asymptotic series carries most of the structural information on Fuchsian groups
which can be detected via subgroup growth.
\begin{theorem}
\label{thm:FuchsStrongEquiv}
Each $\approx$-equivalence class of $\F$ decomposes into finitely many classes with respect to $=;$ that is,
each group $\Gamma\in\F$ has only finitely many presentations of the form {\em (\ref{Eq:GammaDef}),} and is $\approx$-equivalent
to at most finitely many non-isomorphic $\F$-groups.
\end{theorem}
\begin{corollary}
\label{Cor:Finite}
Let $\Gamma\in\mathcal{F}$ be given by a representation as in {\em (\ref{Eq:GammaDef})}.\vspace{-5mm}
\begin{enumerate}
\item The set $\{\Delta\in\mathcal{F}:\Delta\sim\Gamma\}/\cong$ is finite if and only if one of the following holds:\\[-3mm]
\begin{enumerate}
\item $s=t=0,$\\[-3mm]
\item $s=1,\, t=0,\, \sum_{i=1}^r\big(1-\frac{1}{a_i}\big) < 2,$\\[-3mm]
\item $s+2t=2,\, r=1$.\\[-3mm]
\end{enumerate}
\item The set $\{\Delta\in\mathcal{F}:\Delta\sim\Gamma\}/\cong$ is infinite, but
$\{\Delta\in\mathcal{F}: s_n(\Delta)=(1+\mathcal{O}(n^{-2\mu(\Gamma)}))s_n(\Gamma)\}/\cong$ is finite, if and only if
the following three conditions hold:\\[-3mm]
\begin{enumerate}
\item $s+2t+\sum_{i=1}^r\big(1-\frac{1}{a_i}\big) \geq 3,$\\[-3mm]
\item $a_i$ is odd for $1\leq i\leq r,$\\[-3mm]
\item $e_j=2$ for $1\leq j\leq s$ with at most one exception, and for the exceptional index $j_0$ (if it occurs)
we have $e_{j_0}=2^{p-1}$ for some prime $p$.\vspace{-10mm}
\end{enumerate}
\end{enumerate}
\end{corollary}

\begin{proofof}{Theorem~{\em \ref{thm:FuchsStrongEquiv}}}
Let $(r,s,t;a_1,a_2\ldots,a_r,e_1,e_2,\ldots,e_s)$
be a given set of parameters, let $\Gamma$ be the corresponding group and
assume that $\alpha(\Gamma)>0$. We have to show that there are only finitely
many tuples $(r',s',t';a_1',a_2'\ldots,a_r',e_1',e_2',\ldots,e_s')$, such that
for the corresponding group $\Delta$ we have $\Gamma\approx\Delta$. Before computing
the coefficients of the asymptotic series for $s_n(\Gamma)$ and $s_n(\Delta)$,
we show that we may assume without loss that $h_\nu(\Gamma)=h_\nu(\Delta)$ for $\nu=2, 3$.
In fact, from Theorem~\ref{Thm:FuchsWeakEquiv} (ii) we infer that $r+s+2t=r'+s'+2t'$, hence,
$|\Hom(\Delta, S_\nu)|\leq (\nu!)^{r+s+2t}$, that is, there are only finitely many choices
for $h_\nu(\Delta)$ for each fixed $\nu$. Hence, in the sequel we may
assume that $d_\nu(\Gamma)=d_\nu(\Delta)$ for $\nu=2, 3$, where the $d_\nu$ are given as in
(\ref{Eq:DkTransform}), and that $h_n(\Gamma)=\big(1+\mathcal{O}(n^{-3\mu(\Gamma)})\big)h_n(\Delta)$.
Theorem~\ref{Thm:FuchsWeakEquiv} (ii) already
implies that the multisets $\{a_i:1\leq i\leq r\}$ and $\{a_i':1\leq i\leq r'\}$ coincide.
From Propositions~\ref{prop:MainLemFuchs} and \ref{Prop:RootNumberEst} (i), and Equation
(\ref{hncharacter}), we see that
\begin{equation}
\label{Eq:HnCharForFinite}
h_n(\Gamma) = (n!)^{s+2t-2}\prod_{i=1}^r |\Hom(C_{a_i}, S_n)|
\left\{\underset{\Delta< 3\mu(\Gamma)/\alpha(\Gamma)}{\sum_{\lambda\vdash n}}
\frac{\prod_{i=1}^r\alpha_{\chi_\lambda}^{(a_i)} \prod_{j=1}^s m_{\chi_\lambda}^{(e_j)}}
{(\chi_\lambda(1))^{r+s+2t-2}} + \mathcal{O}(n^{-3\mu(\Gamma)})\right\}.
\end{equation}
In view of Proposition~\ref{Prop:RootNumberEst} (ii), the contribution
of partitions $\lambda$ with $3\leq\Delta\leq 3\mu(\Gamma)/\alpha(\Gamma)$
is of lesser order than the error term, hence, we can compute
$h_n(\Gamma)$ up to a relative error of order $n^{-3\mu(\Gamma)}$ using
only coefficients already computed in the previous sections.
Inserting
the values for $\alpha_{\chi_\lambda}^{(q)}$ computed in
Subsection~\ref{Subsec:CoeffComp}, and the multiplicities as given in
Proposition~\ref{Prop:RootNumberEst} (ii) into the right-hand side of
(\ref{Eq:HnCharForFinite}), we obtain
\begin{eqnarray}
h_n(\Gamma) & = & \delta (n!)^{s+2t-2}\prod_{i=1}^r |\Hom(C_{a_i}, S_n)|
\Bigg\{ 1 + (n-1)^{-(r+s+2t-2)}\prod_{i=1}^r H_{1, a_i}(n)
\prod_{j=1}^s\big(\tau(e_j)-1\big)\nonumber\\
&&+\left(\frac{n^2-3n+2}{2}\right)^{-(r+s+2t)}\underset{2|a_i}
{\prod_{1\leq i\leq r}} \frac{n}{n-1}\hspace{.4mm}H_{2, a_i}(n)
\underset{2\nmid a_i}{\prod_{1\leq i\leq r}}\left(H_{2, a_i}(n) -1\right)\nonumber\\
&&\qquad\times\;\prod_{j=1}^s\frac{1}{2}\big(\sigma(e_j)+(\tau(e_j))^2-3\tau(e_j)+
\tau_{\mathrm{odd}}(e_j)\big)\label{Eq:HnExplFinite}\\
&&+\left(\frac{n^2-3n}{2}\right)^{-(r+s+2t)}\underset{2|a_i}
{\prod_{1\leq i\leq r}}\bigg(\frac{n}{n-1}\hspace{.4mm}H_{2, a_i}(n)+1\bigg)
\underset{2\nmid a_i}{\prod_{1\leq i\leq r}}H_{2, a_i}(n)\nonumber\\
&&\qquad\times\;\prod_{j=1}^s\frac{1}{2}\big(\sigma(e_j)+(\tau(e_j))^2-3\tau(e_j)-\tau_{\mathrm{odd}}(e_j)+2\big)
\;+\;\mathcal{O}(n^{-3\mu(\Gamma)})\Bigg\},\nonumber
\end{eqnarray}
where
\[
H_{i, q}(n):= \binom{n}{i}\frac{|\Hom(C_q, S_{n-i})|}{|\Hom(C_q, S_n)|}.
\]
From \cite[Theorem 6]{MComb}, we see that
\[
(n-1)^{-(r+s+2t-2)}\prod_{i=1}^r H_{1, a_i}(n)\asymp n^{-\mu(\Gamma)},
\]
and all other contributions to the asymptotic series in (\ref{Eq:HnExplFinite})
are of lesser order, hence, expanding $h_n(\Delta)$ in the same way, we find
that $\Gamma\approx\Delta$ implies
\begin{equation}
\label{Eq:FiniteTauCond}
\prod_{j=1}^s\big(\tau(e_j)-1\big) = \prod_{j=1}^{s'}\big(\tau(e_j')-1\big).
\end{equation}
Moreover, the contribution of the character $\chi_{(n-1, 1)}$ to $h_n(\Gamma)$
and $h_n(\Delta)$ are identical. Next we consider terms of order $n^{-2\mu(\Gamma)}$.
Arguing as for terms of order $n^{-\mu(\Gamma)}$, we find that $h_n(\Gamma)=(1+o(n^{-2\mu(\Gamma)}))h_n(\Delta)$
is equivalent to (\ref{Eq:FiniteTauCond}) and
\begin{eqnarray}
&&\prod_{j=1}^s\Big(\sigma(e_j)+(\tau(e_j))^2-3\tau(e_j)+
\tau_{\mathrm{odd}}(e_j)\Big)\nonumber\\
&&\qquad\qquad+\;\prod_{j=1}^s\Big(\sigma(e_j)+(\tau(e_j))^2-3\tau(e_j)-\tau_{\mathrm{odd}}(e_j)+2\Big)\nonumber\\
 & = & \prod_{j=1}^{s'}\Big(\sigma(e_j')+(\tau(e_j'))^2-3\tau(e_j')+
\tau_{\mathrm{odd}}(e_j')\Big)\label{Eq:FiniteSigmaCond}\\
&&\qquad\qquad+\;\prod_{j=1}^{s'}\Big(\sigma(e_j')+(\tau(e_j'))^2-3\tau(e_j')-\tau_{\mathrm{odd}}(e_j')+2\Big).\nonumber
\end{eqnarray}
For a fixed tuple $(e_1, \ldots, e_s)$,
there are only finitely many tuples $(e_1', \ldots, e_{s'}')$ with $s'\leq s+2t,$ solving (\ref{Eq:FiniteSigmaCond}).
Indeed, the left-hand side is bounded by some constant,
whereas the right-hand side is bounded below by its greatest
factor, as all factors occurring in the last equation are $\geq 1$; thus $e_j'$ is bounded for all $j$.
\end{proofof}
\begin{proofof}{Corollary~{\em \ref{Cor:Finite}}}
(i) Using Theorem~\ref{Thm:FuchsAsymp}, one checks in each of these cases that $\{\Delta:\Delta\sim\Gamma\}$
is indeed finite. On the other hand, if none of the conditions (a)--(c) is satisfied, one easily computes that
the groups
\[
\Delta_e := \Big\langle x_1, \ldots, x_r, y_1, \ldots, y_{s+2t}\,\big|\,
x_1^{a_1}=\cdots=x_r^{a_r}=x_1\cdots x_r y_1^ey_2^2\cdots y_s^2=1\Big\rangle
\]
satisfy $\alpha(\Delta_e)\geq\frac{2}{e}$, and, if necessary, adjusting the parity of $e$, we have
$\Delta_e\sim\Gamma$ for infinitely many $e$. By Theorem~\ref{thm:FuchsStrongEquiv}, there is an infinite
sequence $(e_\nu)_{\nu\geq 1}$, such that $\Delta_{e_\nu}\sim\Gamma$, while $\Delta_{e_\nu}\not\approx\Delta_{e_\mu}$
for $\nu\neq\mu$, which implies our claim.\\[3mm]
(ii) In the proof of Theorem 7, we have seen that $s_n(\Delta)=(1+\mathcal{O}(n^{-2\mu(\Gamma)})) s_n(\Gamma)$
is equivalent to the conjunction of $\Delta\sim\Gamma$ and
\[
\prod_{j=1}^s (\tau(e_j)-1) = \prod_{j=1}^{s'} (\tau(e_j')-1).
\]
From this equation and part (i) it is easy to see that for a group $\Gamma$ satisfying (a)--(c) the described
sets are of the claimed cardinality. Now assume that $\Gamma$ is a group such that
$\{\Delta\in\mathcal{F}:\Delta\sim\Gamma\}$ is infinite, while $\{\Delta\in\mathcal{F}:
s_n(\Delta)=(1+\mathcal{O}(n^{-2\mu(\Gamma)})) s_n(\Gamma)\}$ is finite.
If $\delta(\Gamma)=1$, and $\{\Delta:\Delta\sim\Gamma\}$ is infinite, there are infinitely many
integers $e'$, such that $(e', 2, 2, \ldots, 2)$ solves this equation, and
$\delta(\Delta_{e'})=\delta(\Gamma)$. Define the integer $q=q(\Gamma)$ via
\[
q:= \left(\prod_{j=1}^s\big(\tau(e_j)-1\big)\right)+1.
\]
If $q$ is not 1 or prime, say $q=a\cdot b$ with $a, b\geq 2$, then
$s_n(\Delta_{2^ap^b})=(1+\mathcal{O}(n^{-2\mu(\Gamma)}))s_n(\Gamma)$ for all odd primes $p$,
and, by Theorem~\ref{thm:FuchsStrongEquiv}, there are infinitely many non-isomorphic groups among the groups $\Delta$
defined in this way. Hence, we may assume that all $e_j$ are even, $\tau(e_j)=2$ for all $j$ with at most one exception,
and for the exceptional index $j_0$, we have that $\tau(e_{j_0})$ is either 2 or $p+1$ for some prime $p$. This implies
our claim.
\end{proofof}

Thomas W. M\"uller, School of Mathematical Sciences, Queen Mary,
University of London, Mile End Road, E1\,4NS London, UK
(T.W.Muller@qmul.ac.uk)

Jan-Christoph Schlage-Puchta, Mathematisches Institut,
Albert-Ludwigs-Universit\"at, Eckerstr.~1, 79104 Freiburg, Germany
(jcp@math.uni-freiburg.de)
\end{document}